\documentclass[11pt]{amsart}
\usepackage{amsxtra,amssymb,amsmath,amsthm,amsbsy}
\usepackage[all]{xy}
\usepackage{enumerate}

\theoremstyle{plain}

\numberwithin{equation}{section}

\newtheorem{theorem}{Theorem}[section]

\newtheorem{definition-lemma}[theorem]{Definition-Lemma}

\theoremstyle{definition}
\newtheorem{example}[theorem]{Example}
\newtheorem{remark}[theorem]{Remark}

\newcounter{exer}
\counterwithin{exer}{section}

\newenvironment{exer} [1][]
{\refstepcounter{exer} \medskip 
\begin{center}
\begin{minipage}{5.0in} 
\footnotesize \noindent {\bf Exercise~\theexer. #1 }}
{\end{minipage}
\end{center}
\medskip}

\newcommand{\Id}         {{\mathrm {Id}}}

\newcommand{\Ker}        {{\mathrm {ker}}}

\newcommand{\Ad}         {{\mathrm {Ad}}}

\newcommand{\pr}         {{\mathrm{pr}}}

\newcommand{\Jac}        {\mathrm{Jac}}

\newcommand{\SP} [1]     {{\left\langle {{#1}} \right\rangle}}

\newcommand{\frakg}     {\mathfrak{g}}

\newcommand{\T}         {\mathbb{T}}

\newcommand{\TM}        {\T M}

\newcommand{\grd}         {\mathcal{G}}
\newcommand{\sour}        {\mathsf{s}}
\newcommand{\tar}         {{\mathsf{t}}}

\newcommand{\Cour}[1]      {[\![#1]\!]}

\newcommand{\Lie}        {\mathcal L}

\begin{document}
\title[]
{A brief introduction to Poisson geometry }

\author[]{Henrique Bursztyn}

\address{Instituto de Matem\'atica Pura e Aplicada,
Estrada Dona Castorina 110, Rio de Janeiro, 22460-320, Brasil }
\email{henrique@impa.br}

\date{}

\begin{abstract}
These notes are based on an introductory minicourse on Poisson geometry given at CRM, Barcelona, in July 2022.
They mostly contain foundational material, including motivating questions and key examples of Poisson structures, and highlight some recent tools and developments. 
The text is largely structured through exercises that lead readers to the desired conclusions, 
so it serves as a guided introduction to the basics of Poisson geometry.

\end{abstract}

\maketitle

\tableofcontents


\section{Introduction}\label{sec:intro}

These notes arose from a minicourse given at Centre de Reserca Matem\`atica, Barcelona, in July 2022, as part of the {\em 2022 Advanced Poisson School}. Versions of this minicourse were subsequently delivered at IMPA in 2023 and 2024, leading to a slight expansion of the original notes.

Keeping the original spirit of the minicourse, these notes focus on basic aspects of Poisson geometry. After a short introductory section below recalling the historical origins of Poisson structures, the notes are divided into four main parts: $\S$ \ref{sec:found} covers foundations, where we define the central objects of interest in Poisson geometry; 
$\S$ \ref{sec:example} presents some of the most relevant classes of examples of Poisson structures; $\S$ \ref{sec:results} concerns basic structural results and classical problems, mostly going back to Weinstein's seminal paper \cite{We83}, but with an effort to highlight refreshed perspectives and recent developments;
$\S$ \ref{sec:topics} briefly discusses additional topics, such as symplectic groupoids and Dirac structures,  that gained interest since \cite{We83} and led to noteworthy advances in the field.

Very little is actually proven in these notes; instead, many results and examples are presented 
through 132 exercises that permeate the main text with the purpose of guiding the reader toward the desired conclusions.

Despite the intention to offer a panoramic view of the subject,
the selection of material in these notes has an inherent bias, and several important topics and applications are omitted. Various textbooks on the subject, such as  \cite{CW,CFMbook,DZ,LGPV,MR}, provide additional material not covered here, as well as further details on topics that we only briefly mention or sketch.
In particular, these notes only treat Poisson geometry of $C^\infty$-manifolds, despite the importance of Poisson structures in the holomorphic and algebraic contexts (see e.g. \cite{Pym} and references therein).

Although not strictly necessary, readers would benefit from familiarity with basic concepts in symplectic geometry.

\bigskip


\noindent {\bf ``The'' Poisson bracket.}
Poisson geometry originated in the work of S. D. Poisson \cite{Poisson}, who introduced the so-called {\em Poisson bracket} as a tool for the geometric formulation of classical mechanics in the 19th century; see \cite{KSHistory} for more on the history. 
In its simplest form, the description of a classical mechanical system consists of a {\em phase space}
$$
\mathbb{R}^{2n} = \{ (q_i, p_i),\; i=1,\ldots, n \}, 
$$
with coordinates $q_i$ representing positions and $p_i$ representing momenta, along with a distinguished function $H\in C^\infty(\mathbb{R}^{2n})$, called the {\em Hamiltonian}, that gives rise to a {\em Hamiltonian vector field}:
$$
X_H :=
\left (\begin{matrix} 0 & \mathrm{Id} \\ -\mathrm{Id} & 0 \end{matrix} \right ) \left (\begin{matrix} \frac{\partial
H}{\partial q_i} \\[0.9mm] \frac{\partial
H}{\partial p_i} \end{matrix} \right)  = \sum_i\frac{\partial
H}{\partial p_i}\frac{\partial}{\partial q_i} - \frac{\partial
H}{\partial q_i}\frac{\partial}{\partial p_i}.
$$
The time evolution of the system is given by integral curves of $X_H$, i.e., curves $\sigma(t)=(q_i(t),p_i(t))$ in phase space satisfying
$\dot{\sigma}(t)=X_H(\sigma(t))$ or, equivalently,
$$
\dot{q}_i = \frac{\partial H}{\partial p_i}, \quad \dot{p}_i = -\frac{\partial H}{\partial q_i},
$$
which are {\em Hamilton's equations}.

In this context, Poisson introduced a skew-symmetric, $\mathbb{R}$-bilinear operation
$$
\{\cdot,\cdot\}: C^\infty(\mathbb{R}^{2n}) \times C^\infty(\mathbb{R}^{2n}) \to C^\infty(\mathbb{R}^{2n}),
$$ 
known as the {\em Poisson bracket}, given by
\begin{equation}\label{eq:original}
\{f,g\}   = \sum_i
\frac{\partial f}{\partial p_i}\frac{\partial g}{\partial q_i} -
\frac{\partial f}{\partial q_i}\frac{\partial g}{\partial p_i}.
\end{equation}

The Poisson bracket has a dynamical meaning as the rate of change of a given function with respect to the Hamiltonian vector field of another function: 
$$
\{H,f\} =\mathcal{L}_{X_H} f.
$$ 
It is immediate that any function is preserved by its own Hamiltonian flow,
\begin{equation}\label{eq:fpres}
    \Lie_{X_f}f = \{f,f\}=0.
\end{equation}

In order to gain insight on (and some times explicitly solve) Hamilton's equations, there was a special interest around Poisson's time in finding {\em conserved quantities} (also called {\em first integrals}) of the system, i.e., functions $f$ such that $\{H,f\} = \mathcal{L}_{X_H} f = 0$. 
(The link between the existence of conserved quantities and finding  solutions to Hamilton's equation has led to the important notion of {\em complete integrability}, see e.g. \cite[Chp.~10]{Arnoldbook}).
Poisson used his brackets to verify the following fact:
\begin{equation}\label{eq:poissonthm}
 \mathcal{L}_{X_H} f =0, \quad \mathcal{L}_{X_H} g =0 \, \implies  \mathcal{L}_{X_H} \{f,g\} =0,
\end{equation}
i.e., the Poisson bracket can be used to produce new conserved quantities from old ones. This result was clarified by C. G. Jacobi in 1842, who observed that Poisson's theorem was a consequence of a fundamental identity satisfied by the Poisson bracket, now known as the {\em Jacobi identity}:
$$
\{h,\{f,g\}\}+\{g,\{h,f\}\} + \{ f,\{g,h\} \}=0, \qquad f,g,h\in C^\infty(\mathbb{R}^{2n}).
$$
The notion of  ``Poisson structure'' appeared later in S. Lie's treatise on transformation groups of 1880; there Lie recognized the importance of the Jacobi identity when setting the foundations of the theory of Lie algebras and Lie groups. 

Despite their important role in physics (not only in classical mechanics but also in the passage to quantum mechanics \cite{Dirac}) 
Poisson brackets were not systematically studied in mathematics until the 1970s.
Two fundamental papers marking the beginning of the ``modern era'' of Poisson geometry are due to Lichnerowicz \cite{Lich} and Weinstein \cite{We83}. The general study of Poisson structures since then has been stimulated by such diverse areas as  geometric mechanics and field theories, integrable systems (in both finite and infinite dimensions), 
representation theory and quantum groups, 
singularity theory, 
noncommutative geometry and deformation quantization, among others;
more information on these topics and further references can be found in the textbooks
\cite{Arnoldbook,CW,DZ,GeShVa,Kac,KhesinWendt,KoroSoib,MR}.

\medskip

\noindent {\bf Acknowledgments}: 
I thank David Iglesias Ponte, Eva Miranda, Cedric Oms, and Roberto Rubio for their encouragement in having these notes written up. I am particularly indebted to David and Roberto for their feedback on  preliminary versions of the manuscript. While preparing these notes, I benefited from discussions with many people, including Daniel Alvarez, Alejandro Cabrera, Pedro H. Carvalho, Matías del Hoyo, and Hudson Lima. I am grateful to all the students from the courses on which these notes are based; I am especially thankful to Andrés Rodriguez and Tianhao Ye for carefully reading the manuscript and providing corrections and suggestions. {\em I am delighted to dedicate these notes to A. Weinstein, whose scientific leadership and fundamental contributions have shaped Poisson geometry as we know it.}


\section{Foundations}\label{sec:found}

\subsection{Poisson brackets}\label{subsub:brackets}

A {\em Poisson bracket} on a smooth manifold $M$ is an $\mathbb{R}$-bilinear operation $\{\cdot,\cdot\}: C^\infty(M)\times C^\infty(M) \to C^\infty(M)$ satisfying

\begin{enumerate}
\item $\{f,g\}=-\{g,f\}$,

\item $\{f,\{g,h\}\}+ \{h,\{f,g\}\} + \{g,\{h,f\}\} = 0$, \qquad\;\,\small{(Jacobi identity)}

\item $\{f,gh\}=\{f,g\}h + \{f,h\}g$, \qquad\qquad \qquad \qquad \qquad \, \small{(Leibniz rule)}
\end{enumerate}
for $f, g, h \in C^\infty(M)$. The pair $(M,\{\cdot,\cdot\})$ is called a {\em Poisson manifold}.

Conditions  (1) and (2) say that $(C^\infty(M),\{\cdot,\cdot\})$ is a Lie algebra, and the third condition is a compatibility of the Lie bracket with the commutative product on $C^\infty(M)$. One can use the same axioms to define a Poisson bracket $\{\cdot,\cdot\}$ on any commutative algebra $\mathcal{A}$, and the pair $(\mathcal{A},\{\cdot,\cdot\})$ is called a {\em Poisson algebra}. If only conditions (1) and (3) are satisfied, one refers to an {\em almost Poisson bracket}. 

A {{\em Poisson map}} from $(M_1,\{\cdot,\cdot\}_1)$ to $(M_2,\{\cdot,\cdot\}_2)$ 
is a smooth map $\varphi: M_1 \to M_2$ such that $\varphi^*: C^\infty(M_2)\to C^\infty(M_1)$ preserves brackets, i.e.,
$$
\{f,g\}_2\circ \varphi = \{f\circ \varphi, g\circ \varphi  \}_1, \qquad \forall f,g \in C^\infty(M_2).
$$
Poisson diffeomorphisms provide a natural notion of equivalence for Poisson manifolds, though there are other weaker notions of equivalence of interest  (see $\S \ref{subsec:sympreal}$).


For an almost Poisson bracket, the Leibniz rule implies that  any $f\in C^\infty(M)$ determines a unique vector field $X_f$ such that
$
\mathcal{L}_{X_f} g = \{f, g \}$, 
called the {{\em Hamiltonian vector field}} of $f$. 
Since
$$
\{f,g\}=dg(X_f)=-df(X_g),
$$
it is clear that, for any open subset $U\subseteq M$, $f|_U=0$ implies that $\{f,g\}|_U=0$ for all $g\in C^\infty(M)$;
hence almost Poisson brackets are local, in the sense that they can be restricted to open subsets.

By the skew symmetry of the bracket, $\mathcal{L}_{X_f}f=\{f,f\}=0$. 
In case of a Poisson bracket, the Jacobi identity ensures that \eqref{eq:poissonthm} holds.

We say that $f, g \in C^\infty(M)$ are in {\em involution} if $\{f,g\}=0$. A function $f$ is in involution with any other function $g$ if and only if $X_f=0$, in which case $f$ is called a {\em Casimir}.

We list a few initial examples.

\begin{example} Any manifold carries the trivial Poisson bracket $\{\cdot,\cdot\}\equiv 0$. A Poisson map from a Poisson manifold $(M,\{\cdot,\cdot\})$ to $(\mathbb{R}^k, 0)$ is same as a collection of functions $f_1,\ldots, f_k \in C^\infty(M)$ in pairwise involution, $\{f_i,f_j\}=0$. \hfill $\diamond$
\end{example}

\begin{example}
A direct calculation shows that the original Poisson bracket \eqref{eq:original} is a Poisson bracket on $M=\mathbb{R}^{2n}$
(see Example~\ref{ex:constant} below to avoid any calculation). 

More generally, let $(M,\omega)$ be a symplectic manifold, i.e., $\omega\in \Omega^2(M)$ is closed and nondegenerate (the map $\omega^\flat: TM\to T^*M$, $X\mapsto i_X\omega$, is an isomorphism). Any $f\in C^\infty(M)$ defines a Hamiltonian vector field $X_f$ by the condition $i_{X_f}\omega = df$, and 
\begin{equation}\label{eq:sympbrk}
\{f,g\}:= \omega(X_g,X_f) = dg(X_f)= \Lie_{X_f}g
\end{equation}
defines a Poisson bracket on $M$ (see $\S$ \ref{subsec:char}). This 
can be  verified by noticing that ``the'' Poisson bracket  \eqref{eq:original} on $\mathbb{R}^{2n}$ corresponds to the symplectic form
$\omega=\sum_i dq_i\wedge dp_i$, so by Darboux's theorem for symplectic forms, any bracket as in \eqref{eq:sympbrk} is locally equivalent to Poisson's original bracket. \hfill $\diamond$
\end{example}

\begin{exer}\label{exer:product} 
Consider Poisson manifolds $(M_1,\{\cdot,\cdot\}_1)$ and $(M_2,\{\cdot,\cdot\}_2)$, and let $M=M_1\times M_2$. Show that the formula
\begin{equation}\label{eq:product}
\{f,g\}(x_1,x_2)= \{f(\cdot,x_2),g(\cdot,x_2)\}_1(x_1) + \{f(x_1,\cdot),g(x_1,\cdot) \}_2(x_2)
\end{equation}
defines a Poisson structure on $M$, in such a way that the projections $p_i: M\to M_i$ are Poisson maps and $\{p_1^*C^\infty(M_1),p_2^*C^\infty(M_2)\}=0$. 

For a Poisson manifold $M$, note that the diagonal map $M\to M\times M$ is not Poisson (unless the Poisson bracket on $M$ is trivial); considering the category of Poisson manifolds and Poisson maps, conclude that the product of Poisson manifolds is not categorical\footnote{A product in a category satisfies the following universal property: given a product of objects $M_1\times M_2$ with projections $p_i: M_1\times M_2\to M_i$, $i=1,2$, then for any object $N$ and morphisms $\psi_i: N\to M_i$, $i=1,2$, there us a unique morphism $\psi: N\to M_1\times M_2$ such that $p_i\circ \psi=\psi_i$, $i=1,2$.} (cf. Exercise~\ref{exer:prod2} (b)). 
\end{exer}

\begin{example} \label{ex:family}
Let $S$ be a manifold with a smooth family of symplectic structures $\omega_t \in \Omega^2(S)$, $t\in \mathbb{R}$. Denote  by $\{\cdot,\cdot\}_t$ the Poisson bracket defined by $\omega_t$. Then there is a Poisson bracket on $M= S \times \mathbb{R}$ given by
$$
\{f,g\}(x,t) = \{f(\cdot,t), g(\cdot,t)\}_t(x).
$$
(More generally, in this example $\mathbb{R}$ could be replaced by any other manifold, and $\{\cdot,\cdot\}_t$ could be any smooth family of Poisson brackets, not necessarily symplectic.) \hfill $\diamond$
\end{example}

\subsection{Tensorial viewpoint: bivector fields}\label{subsub:tensor}
Given a manifold $M$, we will use the notation $\mathfrak{X}^k(M): =\Gamma(\wedge^k TM)$. In particular, the space of vector fields on $M$ will be denoted by $\mathfrak{X}^1(M)$.

For an  almost Poisson bracket $\{\cdot,\cdot\}$ on $M$, since 
$$
\{f,g\} = dg(X_f) = -df(X_g),
$$
we see that the  bracket of two functions only depends on their differentials. It follows that there is a unique bivector field $\pi\in \mathfrak{X}^2(M)$ such that
\begin{equation}\label{eq:bivector}
\{f,g\} = \pi(df, dg).
\end{equation}

By means of the previous formula one obtains a bijective correspondence between almost Poisson brackets on $C^\infty(M)$ and bivector fields $\pi$ on $M$. We will refer to a bivector field whose corresponding bracket satisfies the Jacobi identity
as a {\em Poisson bivector field}. One can therefore think of a {\em Poisson structure} on $M$ as either a  Poisson bracket $\{\cdot,\cdot\}$ on $C^\infty(M)$, or as a Poisson bivector field $\pi\in \mathfrak{X}^2(M)$.

In local coordinates $(x_1,\ldots,x_n)$ on $M$, a bivector field $\pi$ is written as
$$
\pi  = \frac{1}{2}\sum_{i,j} \pi_{ij}(x) \frac{\partial}{\partial x_i}\wedge \frac{\partial }{\partial x_j} = \sum_{i<j} \pi_{ij}(x) \frac{\partial}{\partial x_i}\wedge \frac{\partial }{\partial x_j},
$$
and the associated almost Poisson bracket takes the form
$$
\{f,g\}(x)  = \sum_{i,j}\pi_{ij}(x)\frac{\partial f}{\partial
x_i}\frac{\partial g}{\partial x_j}, \qquad \pi_{ij} = \{x_i,x_j\}.
$$

A bivector field $\pi$ on $M$ can be equivalently described by a vector-bundle map $\pi^\sharp: T^*M\to TM$ (over the identity map on $M$) satisfying $(\pi^\sharp)^*=-\pi^\sharp$ via 
$$
\pi(\alpha,\beta) = \beta(\pi^\sharp(\alpha)).
$$ 
Note that, for $f\in C^\infty(M)$,
$$
\pi^\sharp(df)=X_f.
$$ 

\begin{exer}\label{exer:pmaps}
Let $(M_1,\pi_1)$ and $(M_2,\pi_2)$ be Poisson manifolds. Show that the following are equivalent.
\begin{itemize}
\item[(a)] $\varphi: M_1\to M_2$ is a Poisson map; 
\item[(b)] $X_{\varphi^*f}$ is $\varphi$-related to $X_f$ for all $f\in C^\infty(M_2)$; 
\item[(c)] 
$\pi_2^\sharp|_{\varphi(x)}= d_x\varphi \circ \pi_1^\sharp|_x \circ (d_x\varphi)^*$ for all $x\in M_1$.  
\end{itemize}
\end{exer}

\begin{exer} \label{exer:prod2}
\begin{itemize}
\item[(a)]For Poisson manifolds $(M_1,\pi_1)$ and $(M_2,\pi_2)$,
show that the Poisson bracket on $M_1\times M_2$ in Exercise \ref{exer:product} corresponds to the natural direct product of bivector fields $\pi_1\times \pi_2 \in \mathfrak{X}^2(M_1\times M_2)$.
\item[(b)] Suppose that $\varphi_1: M\to M_1$ and $\varphi_2: M\to M_2$ are Poisson maps. Show that $M\stackrel{(\varphi_1,\varphi_2)}{\longrightarrow} M_1\times M_2$ is a Poisson map if and only if  $\{\varphi_1^* C^\infty(M_1), \varphi_2^* C^\infty(M_2)\}=0$.
\end{itemize}
\end{exer}

Given a bivector field $\pi$ on $M$, the failure of the Jacobi identity of the corresponding bracket $\{\cdot,\cdot\}$ is measured by the {\em Jacobiator}, i.e., 
the trilinear operation $\mathrm{Jac}_\pi: C^\infty(M)\times C^\infty(M)\times C^\infty(M) \to C^\infty(M)$,
$$
\mathrm{Jac}_\pi(f,g,h)=\{f,\{g,h\}\} + \{h,\{f,g\}\} + \{g,\{h,f\}\}.
$$
We will simply write $\mathrm{Jac}$ whenever $\pi$ is clear from the context.

\begin{exer}\label{exer:jac}  Verify that
$$
\mathrm{Jac}(f,g,h) = \mathcal{L}_{[X_f,X_g]}h- \mathcal{L}_{X_{\{f,g\}}}h = (\mathcal{L}_{X_f}\pi)(dg,dh).
$$
\end{exer}

It follows from the previous exercise that, for a Poisson structure, the map $C^\infty(M)\to \mathfrak{X}^1(M)$, $f \mapsto X_f$ is a Lie algebra homomorphism,
\begin{equation}\label{eq:brkpres}
[X_f, X_g]= X_{\{f,g\}}, \qquad \forall \, f,g \in C^\infty(M),
\end{equation}
as well as the fact that Poisson bivector fields are preserved by Hamiltonian flows. 
Another direct consequence of Exercise \ref{exer:jac} is that the Jacobiator only depends on the differentials of functions,  so there exists a unique trivector field $\Upsilon_\pi\in \mathfrak{X}^3(M)=\Gamma(\wedge^3 TM)$ such that, for all $f$, $g$, $h \in C^\infty(M)$,
\begin{equation}\label{eq:jac}
\mathrm{Jac}(f,g,h) = \Upsilon_\pi(df, dg, dh).
\end{equation}
Hence the vanishing of $\mathrm{Jac}$ can be verified by checking the vanishing of $\Upsilon_\pi$ on local frames of $T^*M$. In particular, in a coordinate chart $(x_1,\ldots, x_n)$, it suffices to check that $\mathrm{Jac}(x_i, x_j, x_k)=0$, \; $ \forall\, i$, $j$, $k$.

\begin{example}\label{ex:constant}
On any open subset of $\mathbb{R}^m$, a constant bivector field is Poisson. In coordinates $(x_1,\ldots,x_m)$, the fact that  $\pi_{ij}=\{x_i,x_j\}$ is constant implies that
$\{x_i, \{x_j,x_k\}\}= \Lie_{X_{x_i}}\{x_j,x_k\}=0$, so $\mathrm{Jac}(x_i,x_j,x_k)=0$ for all $i,j, k$.
This is the case for Poisson's original bracket \eqref{eq:original}, which corresponds to the bivector field
$$
\pi_{can} = \sum_i \frac{\partial}{\partial p_i}\wedge \frac{\partial}{\partial q_i} 
$$
in $\mathbb{R}^{2n}$. \hfill $\diamond$

\end{example}

\begin{example} \label{exer:2d} On a 2-dimensional manifold $M$ any trivector field is trivial, so any bivector field is Poisson. E.g. any $\varphi\in C^\infty(\mathbb{R}^2)$ defines a Poisson structure
$\pi = \varphi(x,y)\partial_{x}\wedge \partial_y$ on $\mathbb{R}^2$. \hfill $\diamond$
\end{example}

\begin{example}\label{ex:su2}
On $M=\mathbb{R}^3=\{\xi= (x,y,z)\}$, we have the Poisson bracket
$$
\{f,g\}(\xi)= \SP{\xi,(\nabla f|_\xi\times \nabla g|_\xi)},
$$
with corresponding bivector field $\pi=z \partial_x\wedge \partial_y + x \partial_y\wedge \partial_z + y \partial_z\wedge \partial_x$. In this case, 
integral curves of the Hamiltonian vector field of the function $H(x,y,z)=\frac{x^2}{2I_x} + \frac{y^2}{2I_y} +\frac{z^2}{2I_z}$ are the same as solutions of the Euler equations of a rigid body with moments of inertia $I_x$, $I_y$, $I_z$ (see e.g. \cite[$\S$ 15]{MR}). \hfill $\diamond$
\end{example}

\subsection{The Schouten bracket}\label{subsec:schouten}

For a bivector field $\pi$ on $M$, there is a more intrinsic way to express the trivector field $\Upsilon_\pi$ in terms of $\pi$ through the so-called {\em Schouten bracket}. For each $k=0,1,\ldots$, consider the space $\mathfrak{X}^k(M)=\Gamma(\wedge^k TM)$ of $k$-vector fields on $M$, and the (graded commutative) algebra 
$$
\mathfrak{X}^\bullet(M) = \oplus_{k=0}^{\mathrm{dim}(M)}\mathfrak{X}^k(M)
$$ 
of multivector fields (here $\mathfrak{X}^0(M)=C^\infty(M)$). The Schouten bracket extends the usual Lie bracket of vector fields as follows: 
it is the unique $\mathbb{R}$-bilinear bracket
$$
[\cdot,\cdot]: \mathfrak{X}^k(M)\times \mathfrak{X}^l(M)\to \mathfrak{X}^{k+l-1}(M)
$$ 
such that, for $X\in \mathfrak{X}^k(M)$, $Y\in \mathfrak{X}^l(M)$,
\begin{itemize}
    \item $[X,Y] = - (-1)^{(k-1)(l-1)}[Y,X]$,
    \item $[X,Y\wedge Z] = [X,Y]\wedge Z + (-1)^{(k-1)l}Y\wedge [X,Z]$,
    \item If $Z$ is a vector field, then $[Z,\cdot]=\mathcal{L}_Z$.
\end{itemize}
It follows from these properties that $[\cdot,\cdot]$ satisfies a graded version of the Jacobi identity, 
$$
(-1)^{(k-1)(m-1)} [X,[Y,Z]] + (-1)^{(m-1)(l-1)}[Z, [X,Y]] + (-1)^{(l-1)(k-1)}[Y, [Z,X]]=0,
$$
for $X\in \mathfrak{X}^k(M)$, $Y\in \mathfrak{X}^l(M)$, $Z\in \mathfrak{X}^m(M)$, so the Schouten bracket makes $\mathfrak{X}(M)$ into a Gerstenhaber algebra (i.e., $[\cdot,\cdot]$ is a graded Poisson bracket of degree -1).

\begin{exer}\label{exer:schouten}
Given $\pi\in \mathfrak{X}^2(M)$, check that
\begin{itemize}
\item[(a)]  $[\pi,f]=-X_f$,  for $f\in C^\infty(M)$;
\item[(b)] $[\pi, [\pi, \cdot]] = \frac{1}{2}[[\pi,\pi], \cdot]$.
\end{itemize}
\end{exer}

\begin{remark}
In local coordinates $(x_1,\ldots,x_n)$ and using the notation $\xi_i=\frac{\partial}{\partial x_i}$, one can write $X\in \mathfrak{X}^k(M)$ and $Y\in \mathfrak{X}^l(M)$ as
$$
X=\sum_{i_1<\ldots<i_k}a_{i_1\cdots i_k}(x) \xi_{i_1}\ldots \xi_{i_k}, \quad Y= \sum_{i_1<\ldots<i_l}b_{i_1\cdots i_l}(x) \xi_{i_1}\ldots \xi_{i_l},
$$
so we can regard them as ``functions'' of $x_i$ and $\xi_i$. With this notation, the Schouten bracket of $X$ and $Y$ has the following local expression:
$$
[X,Y] = \sum_i \frac{\partial X}{\partial \xi_i}\frac{\partial Y}{\partial x_i} - (-1)^{(k-1)(l-1)}
\frac{\partial X}{\partial x_i}\frac{\partial Y}{\partial \xi_i}.
$$
By formally interpreting $x_i$, $\xi_j$ as local coordinates, one can regard  the Schouten bracket as an analogue of Poisson's original bracket on the ``shifted'' cotangent bundle $T^*[1]M$ (the notation ``$[1]$'' indicates that fiber coordinates are assigned degree 1), see e.g. \cite[$\S$ 2]{RoyGraded}. \hfill $\diamond$
\end{remark}

One can verify (see e.g. \cite{DZ}) that
$$
\Upsilon_\pi = \frac{1}{2}[\pi,\pi].
$$
So a Poisson bivector field is characterized by the condition $[\pi,\pi]=0$.

\begin{exer}\label{exer:ld}
Consider vector fields $X,Y$ on $M$ and let $\pi=X\wedge Y \in \mathfrak{X}^2(M)$. Verify that $\pi$ is Poisson if and only if $X$, $Y$ and $[X,Y]$ are linearly dependent at all points. 
\end{exer}

%

\subsection{The characteristic distribution}\label{subsec:char}

By a {\em distribution} on a manifold $M$ we will mean a subset $D\subseteq TM$ such that $D_x:= D\cap T_xM$ is a vector subspace of $T_xM$. A distribution $D$ is called {\em smooth} if for any $x\in M$ and $v\in D_x$, there exists a vector field $X\in \mathfrak{X}^1(M)$ such that $X|_x=v$ and $X|_y\in D_y$ for all $y\in M$.
The dimension of $D_x$ is called the {\em rank} of $D$ at $x$, and a distribution with constant rank is called {\em regular}. A distribution that is smooth and regular is the same as a vector subbundle of $TM$.

\begin{exer}
Let $\omega\in \Omega^2(M)$ and $\omega^\flat: TM\to T^*M$, $\omega^\flat(X)=\omega(X,\cdot)$. Consider the distribution 
$$
K=\ker(\omega^\flat)=\{X\in TM\,|\, \omega(X,\cdot)=0\}\subseteq TM.
$$
Show that $K$ is smooth if and only if it is locally regular (i.e., regular on connected components). (It might be helpful to notice that the rank of a smooth distribution is a lower semi-continuous function on $M$.) 
\end{exer}

For a bivector field $\pi \in \mathfrak{X}^2(M)$,  the image of the map
\begin{equation}\label{eq:anchor}
\pi^\sharp: T^*M \to TM, \; \alpha \mapsto \pi(\alpha,\cdot),
\end{equation}
defines a distribution on $M$,
$$
R := \pi^\sharp(T^*M) \subseteq TM,
$$
called the
{\em characteristic distribution} of $\pi$.
Note that $R_x=\{ X_f|_x\;,\; f\in C^\infty(M)\}$.

\begin{exer}
For any bivector field $\pi$, show that $R$ is smooth.
\end{exer}

As the next exercise shows, each $R_x$ naturally acquires the structure of a symplectic vector space, and the collection of these symplectic vector spaces completely codifies the bivector field $\pi$. (In case of Poisson structures, much more will be true: the characteristic distribution will give rise to ``symplectic leaves'', see $\S$ \ref{subsec:integ}.)

\begin{exer}\label{exer:fibersympl} Let $V$ be a (real) vector space, and let $\pi \in \wedge^2 V$.
Consider $\pi^\sharp: V^*\to V$ defined by $\beta(\pi^\sharp(\alpha))=\pi(\alpha,\beta)$, and let $R=\pi^\sharp(V^*)\subseteq V$. Show that there is a unique nondegenerate skew-symmetric bilinear form $\Omega$ on $R$ given by $\Omega(u,v)=\pi(\beta,\alpha)$, for $u=\pi^\sharp(\alpha)$ and $v=\pi^\sharp(\beta)$. Conversely, show that given a pair $(R,\Omega)$, where $R\subseteq V$ is a subspace and $\Omega\in \wedge^2R^*$ is nondegenerate, there is a unique $\pi\in \wedge^2V$ such that $R=\pi^\sharp(V^*)$ and $\Omega$ is defined as before. 
\end{exer}

We define the {\em rank} of  $\pi$ at $x \in M$ as the rank of $R$ at $x$, and we say that $\pi$ is {\em regular} if so is $R$.  
We say that a bivector field $\pi$ is
 {\em nondegenerate} if  $R=TM$, or, equivalently, if $\pi^\sharp: T^*M\to TM$ is an isomorphism.

Any nondegenerate 2-form $\omega\in \Omega^2(M)$ defines a nondegenerate
bivector field $\pi$ by $\pi^\sharp=(\omega^\flat)^{-1}$, or 
$$
\pi(\beta,\alpha)= \omega((\omega^\flat)^{-1}(\alpha),(\omega^\flat)^{-1}(\beta)),
$$
and this establishes a bijective correspondence

\smallskip

\begin{center}{ {nondegenerate 2-forms \; $\rightleftharpoons$ \; {nondegenerate} bivector fields}}.
\end{center}

\smallskip

\begin{exer}\label{exer:poisymp} Suppose that $\pi$ is a nondegenerate bivector field with corresponding 2-form $\omega$, so that, for $f\in C^\infty(M)$, $X_f=\pi^\sharp(df)$ if and only if $i_{X_f}\omega=df$.  For $f,g,h \in C^\infty(M)$, compare $d\omega(X_f,X_g,X_h)$ and $ \Jac(f,g,h)$, and verify
that $\pi$ is Poisson if and only if $\omega$ is symplectic. 
\end{exer}

It follows from the previous exercise that symplectic structures are equivalent to nondegenerate Poisson structures:

\smallskip

\begin{center}{ {Symplectic structures \; $\rightleftharpoons$ \;  {nondegenerate} Poisson structures}}.
\end{center}

\smallskip

In spite of this equivalence, maps that preserve Poisson brackets are not the same as maps that preserve symplectic forms, as explained in the next exercise.


\begin{exer} Consider symplectic manifolds $(M_i,\omega_i)$, with corresponding Poisson brackets $\{\cdot,\cdot\}_i$, $i=1,2$, and let $\varphi: M_1\rightarrow M_2$ be a smooth map.
\begin{itemize}
\item[(a)] Show that if $\varphi$ preserves symplectic forms (i.e.,  $\varphi^*\omega_2=\omega_1$), then it must be an immersion, and that if $\varphi$ is a Poisson map, then it must be a submersion.

\item[(b)] Prove that, if $\varphi$ is a (local) diffeomorphism, then it is a Poisson map if and only if it preserves symplectic forms.

\item[(c)] Find examples of $M_1$, $M_2$ and $\varphi: M_1\rightarrow M_2$ such that (1) $\varphi$ is a Poisson map but does not preserve symplectic forms; (2) $\varphi$ preserves symplectic forms  but is not a Poisson map.

\end{itemize}

\end{exer}

\subsection{Poisson, coisotropic and cosymplectic submanifolds}
We briefly discuss some important types of submanifolds in Poisson geometry.
By a submanifold of a manifold $M$ we mean a manifold $N$ together with an injective immersion $\iota: N\hookrightarrow M$; we will often omit $\iota$ and identify $N$ with a subset of $M$. 

Let $\zeta$ be a $k$-vector field on $M$. We say that $\zeta$ is {\em tangent} to a submanifold $N \hookrightarrow M$ if $\zeta_x\in \wedge^k T_xN \subseteq \wedge^k T_xM$ for all $x\in N$. In this case $N$ inherits a $k$-vector field $\zeta_N$ by restriction of $\zeta$. 

\begin{exer}\label{exer:tangentmulti} 
Check that a $k$-vector field $\zeta$ is tangent to $N\hookrightarrow M$ if and only if, for any open subset $U\subseteq N$,
$(i_{df}\zeta)|_U=0$ for all $f\in C^\infty(M)$ such that $f|_U=0$.
\end{exer}

\begin{exer}\label{exer:tangent}
\begin{itemize}
\item[(a)] Show that a bivector field $\pi$ on $M$ is tangent to a submanifold $N$  if and only if $R|_N=\pi^\sharp(T^*M|_N)\subseteq TN$ if and only if $\pi^\sharp(\mathrm{Ann}(TN))=0$. In this case, check that the characteristic distribution of the induced bivector field $\pi_N$ on $N$ is $R|_N$.

\item[(b)] Consider the Jacobiator 3-vector field $\Upsilon_\pi$, see \eqref{eq:jac}. 
Show that if $\pi$ is tangent to $N$, then, for any open subset $U\subseteq N$,  $\mathrm{Jac}_\pi(f,g,h)|_U= \mathrm{Jac}_{\pi_N}(f|_U,g|_U,h|_U)$ for all $f,g,h \in C^\infty(M)$, and hence
$\Upsilon_\pi$ is also tangent to $N$ (see Exercise \ref{exer:tangentmulti}). 
\end{itemize}
\end{exer}

A submanifold $N$ of a Poisson manifold $M$ is a {\em Poisson submanifold} if it is equipped with a Poisson structure for which the inclusion map $N\hookrightarrow M$ is a Poisson map. (It is clear that there is at most one Poisson structure on a submanifold making it into a Poisson submanifold.)


\begin{exer}
Let $N$ be a submanifold of a Poisson manifold $(M,\pi)$.
Show that $N$ is a Poisson submanifold if and only if $\pi$ is tangent to $N$, if and only if every Hamiltonian vector field on $M$ is tangent to $N$.

\end{exer}

\begin{exer}\label{exer:Ipois}
Let $N$ be a submanifold and $I_N= \{ f \in C^\infty(M)\,|\, f|_N=0\} \subseteq C^\infty(M)$ its vanishing ideal. Check that if $N$ is Poisson then $I_N$ is a Lie-ideal, i.e., $\{C^\infty(M),I_N\}\subseteq I_N$, and that the converse holds if $N$ is embedded.
\end{exer}

\begin{exer}\label{exer:casimir}
Let $(M,\pi)$ be a Poisson manifold and $\Psi=(\psi_1,\ldots,\psi_k): M \to \mathbb{R}^k$ a map such that $r\in \mathbb{R}^k$ is a regular value, and consider the submanifold $N = \Psi^{-1}(r)$. Show that $N$ is Poisson if and only if  $X_{\psi_i}|_N=0$ for all $i=1,\ldots,k$. (In particular, level sets of Casimirs are Poisson submanifolds.)
\end{exer}

\begin{exer}\label{exer:S2}
Consider $\mathbb{R}^3$ with the Poisson structure  $\pi=z \partial_x\wedge \partial_y + x \partial_y\wedge \partial_z + y \partial_z\wedge \partial_x$. Show that the unit sphere $\mathbb{S}^2 \subseteq \mathbb{R}^3$ is a Poisson submanifold, and that the Poisson structure it inherits corresponds to the symplectic structure given by the negative of its area form.
\end{exer}

\begin{exer} \label{exer:spheres}
On $\mathbb{R}^{n+1}=\{ (x,y, z) \,|\, x,y \in \mathbb{R}, \, z=(z_1,\ldots,z_{n-1}) \in \mathbb{R}^{n-1} \}$, consider the bivector field $\pi$
defined by the bracket relations
$$
\{x,y\}=|z|^2,\; \{x,z_i\}= -yz_i,\; \{y,z_i\}=x z_i ,\; \{z_i,z_j\}=0, 
$$
where $|z|^2=z_1^2 + ... +z_{n-1}^2$.
Show that $\pi$ is a Poisson structure on $\mathbb{R}^{n+1}$ and the unit sphere $\mathbb{S}^n\subseteq \mathbb{R}^{n+1}$ is a Poisson submanifold.
\end{exer}

The Poisson submanifolds of a symplectic manifold are just its open subsets. We will now see ways to extend the notions of coisotropic and symplectic submanifolds to the Poisson setting. (In physics, these are known as {\em first-class} and {\em second-class} constraints, respectively.)

Given a submanifold $N$ of a Poisson manifold $(M,\pi)$, we define 
\begin{equation}\label{eq:TNpi}
TN^\pi= \pi^\sharp(\mathrm{Ann}(TN)),
\end{equation}
which is the Poisson analogue of symplectic orthogonals\footnote{If $W$ is a subspace of a symplectic vector space $(V,\Omega)$, its symplectic orthogonal is $W^\Omega:=\{v\in V\,|\, \Omega(v,w)=0 \, \forall w\in W\}=(\Omega^\flat)^{-1}(\mathrm{Ann}(W))$. Recall that a subspace $W\subseteq V$ is {\em coisotropic} if $W^\Omega\subseteq W$ and {\em symplectic} if $W\cap W^\Omega=\{0\}$ (equivalently, if $V=W\oplus W^\Omega$).} in symplectic geometry. With this notation, Poisson submanifolds are defined by the condition $TN^\pi=0$ (Exercise \ref{exer:tangent}).

\begin{exer}\label{exer:TNpibrk}
Check that $TN^\pi = \mathrm{Ann}((\pi^\sharp)^{-1}(TN))$.
Use this fact to show that, if $f$, $g\in C^\infty(M)$ are such that $df|_N, dg|_N \in \Gamma(\mathrm{Ann}(TN^\pi))$, 
 then $(d\{f,g\}) |_N \in \Gamma(\mathrm{Ann}(TN^\pi))$.
\end{exer}

\begin{exer}\label{exer:TNpi}
\begin{itemize}
\item[(a)] Show that the following holds pointwise: $TN^\pi = (TN\cap R)^\Omega\subseteq R$ (the symplectic orthogonal of $TN\cap R$ in $R$), where we regard  the characteristic distribution $R=\pi^\sharp(T^*M)$ as a family of symplectic vector spaces (see Exercise~\ref{exer:fibersympl}). In particular, $(TN^\pi)^\pi = TN\cap R$.

\item[(b)] Check that $\mathrm{rank}(TN^\pi)\leq \mathrm{dim}(M)-\mathrm{dim}(N)$, and hence $\mathrm{rank}(TN)+ \mathrm{rank}(TN^\pi)\leq \mathrm{dim}(M)$. 
\end{itemize}
\end{exer}

A submanifold $N$ is called {\em coisotropic} if $TN^\pi\subseteq TN$. 
An immediate consequence of the previous exercise is that $N$ is a coisotropic submanifold if and only if, at each point, $TN\cap R$ is a coisotropic subspace of $R$.

\begin{exer}\label{exer:coisoIn} Following the notation of Exercise \ref{exer:Ipois}, show that if $N$ is coisotropic then $\{I_N,I_N\} \subseteq I_N$, and that the converse holds as long as $N$ is embedded. 
\end{exer}

\begin{exer}
Consider the submanifold $N$ of $M$ as in Exercise~\ref{exer:casimir}.
Show that $N$ is coisotropic if and only if $\{\psi_i,\psi_j\}|_N=0$  for all $i,j=1,\ldots,k$.
\end{exer}

The following three exercises illustrate the role of coisotropic manifolds in Poisson geometry, see \cite{WeinCoiso}.

\begin{exer}
Show that a map $ M_1\to M_2$ between Poisson manifolds is a Poisson map if and only if its graph is a coisotropic submanifold of $M_1\times \overline{M_2}$, where $\overline{M_2}$ has minus the Poisson structure of $M_2$. 
\end{exer}

\begin{exer}\label{exer:invim}
Consider a Poisson map $\varphi: M_1\to M_2$ and a coisotropic submanifold $N\hookrightarrow M_2$ such that $\varphi^{-1}(N)$ is a submanifold of $M_1$ with $T(\varphi^{-1}(N))=(d\varphi)^{-1}(TN)$ (this is the case e.g. when $\varphi$ is transverse to $N$). Show that $\varphi^{-1}(N)$ is coisotropic.
\end{exer}


\begin{exer}
Let $\varphi: M \to B$ be a surjective submersion, and consider the submanifold $M\times_B M= \{(x,y) \in M\times M\,|\, \varphi(x)=\varphi(y)\}  \subseteq M\times M$. 

Suppose that $M$ is a Poisson manifold. Show that $B$ carries a Poisson structure for which $\varphi$ is a Poisson map if and only if $M\times_B M$ is a coisotropic submanifold of $M\times \overline{M}$. (Exercises \ref{exer:coisoIn} and \ref{exer:invim} can be helpful.)
\end{exer}

Symplectic submanifolds can be generalized to Poisson geometry in different ways; we will discuss one possibility now (see $\S$ \ref{subsec:dirac} for a more thorough discussion).

Given a Poisson manifold $(M,\pi)$, a submanifold $N$ is called {\em cosymplectic}  if
\begin{equation}\label{eq:cosym}
TN \oplus TN^\pi = TM|_N.
\end{equation}
By Exercise~\ref{exer:TNpi} (b), this last condition is equivalent to $TN + TN^\pi=TM|_N$.

\begin{exer}\label{exer:nondeg}
Let $N$ be a submanifold of $(M,\pi)$, and consider the (fiberwise) skew-symmetric bilinear form on $\mathrm{Ann}(TN)$ defined by the restriction of $\pi$. Show that $N$ is cosymplectic if and only if this bilinear form is nondegenerate on each fiber.
\end{exer}

It follows from the previous exercise that $\mathrm{Ann}(TN)$ is a symplectic vector bundle, and this motivates the name ``cosymplectic'' for this class of  submanifolds\footnote{The term ``cosymplectic'' in this context is not to be confused with the concept of ``cosymplectic structure'' due to Libermann, defined on an odd-dimensional manifold $M^{2n+1}$ by a closed 1-form $\eta$ and a closed 2-form $\omega$ such that $\eta\wedge \omega^n$ is a volume form; see \cite[Prop.~4.19]{CFMbook} for a description of how such structures are actually related to Poisson structures.}; cosymplectic submanifolds are alternatively called  {\em Poisson transversals} \cite{FrejMarcNormal} (see Example~\ref{example:PDsubm} (b) for further comments on terminology).


\begin{exer}\label{exer:cosympR} Let $N$ be a submanifold of $(M,\pi)$ and $R=\pi^\sharp(T^*M)$  (see Exercise~\ref{exer:TNpi} (a)). Show that
\begin{itemize}
\item[(a)] pointwise, $TN\cap TN^\pi = \{0\}$ if and only if $TN\cap R\subseteq R$ is a symplectic subspace;

\item[(b)] a submanifold $N$ is cosymplectic if and only if $TN\cap R\subseteq R$ is a symplectic subspace pointwise and $N$ is transverse to $R$, i.e., $TN+R|_N=TM|_N$.
\end{itemize}
\end{exer}

\begin{exer}\label{exer:slice}
Let $N$ be a submanifold of $(M,\pi)$. Show that if the cosymplectic condition \eqref{eq:cosym} holds at a given point in $N$, then it holds locally (in $N$) around this point. Use this fact to  show that, if $x\in N$ is such that $T_xN \oplus R_x=T_xM$, then there is a neighborhood  
  of $x$ in $N$ that is a cosymplectic submanifold.
\end{exer}

A cosymplectic submanifold $N$ is not a Poisson (or, more generally, coisotropic) submanifold unless $TN=TM|_N$, i.e., $N$ is an open subset of $M$. Nevertheless, a key property of cosymplectic submanifolds is that they always inherit a natural Poisson structure (a further discussion of how this Poisson structure relates to the ambient Poisson structure on $M$ is presented in $\S$ \ref{subsec:dirac}).

If $N$ is a cosymplectic submanifold of $(M,\pi)$, we can write 
$$
T^*M|_N= T^*N\oplus (TN^\pi)^*,
$$
using the identifications
$$
T^*N\stackrel{\sim}{\rightarrow} \mathrm{Ann}(TN^\pi)\subseteq T^*M|_N,\quad \; 
(TN^\pi)^*\stackrel{\sim}{\rightarrow} \mathrm{Ann}(TN) \subseteq T^*M|_N,
$$
given by the dual maps to the projections $TM|_N \twoheadrightarrow TN$ and $TM|_N\twoheadrightarrow TN^\pi$.
Since $\pi^\sharp(\mathrm{Ann}(TN))=TN^\pi$ and $\pi^\sharp(\mathrm{Ann}(TN^\pi))=(TN^\pi)^\pi=TN\cap R \subseteq TN$ (Exercise \ref{exer:TNpi} (a)), we see that $\pi^\sharp|_N: T^*M|_N\to TM|_N$ decomposes as maps
\begin{equation}\label{eq:piNchi}
T^*N\to TN, \quad \mathrm{ and } \quad (TN^\pi)^*\to TN^\pi.
\end{equation}
The first map defines a bivector field $\pi_N$ on $N$, while the second map defines an element $\chi \in \Gamma(\wedge^2 TN^\pi)$, so that 
$$
\pi|_N = \pi_N + \chi.
$$
Note that $\chi$ is nondegenerate, cf. Exercise~\ref{exer:nondeg}.

\begin{exer}\label{exer:charpiN}
Following Exercise~\ref{exer:cosympR}, check that the characteristic distribution of $\pi_N$ is $TN\cap R$, with pointwise symplectic form given by pullback of the one on $R$ . 
\end{exer}

The next exercise indicates how to verify that the bivector field $\pi_N$ is indeed a Poisson structure on $N$. We assume that $N$ is embedded for simplicity (otherwise one can work locally on $N$).

\begin{exer}\label{exer:cosymPoiss} 
\begin{itemize}
\item[(a)] Show that the bracket $\{\cdot,\cdot\}_N$ on $C^\infty(N)$ corresponding to $\pi_N$ can be calculated as follows:
$$
\{f,g\}_N= \{\widehat{f},\widehat{g}\}|_N,
$$
where $\widehat{f}, \widehat{g} \in C^\infty(M)$ are extensions of $f, g \in C^\infty(N)$ satisfying $d\widehat{f}|_{TN^\pi} = d\widehat{g}|_{TN^\pi}=0$. (Note that such extensions always exist.)

\item[(b)] Use (a) and Exercise~\ref{exer:TNpibrk} to show that 
$$
\mathrm{Jac}_{\pi_N}(f,g,h)= \mathrm{Jac}_{\pi}(\widehat{f},\widehat{g},\widehat{h})|_N,\qquad f,g,h \in C^\infty(N).
$$ 
Therefore
$\{\cdot,\cdot\}_N$ satisfies the Jacobi identity, so $\pi_N$ is a Poisson structure on $N$.
\end{itemize}
\end{exer}

The following exercise gives a general local formula for the bracket $\{\cdot,\cdot\}_N$ in terms of $\{\cdot,\cdot\}$, known as the {\em Dirac bracket}.

\begin{exer}[(Dirac bracket)]  Consider a Poisson manifold $(M,\pi)$,  a map $\Psi=(\psi_1,\ldots,\psi_k): M \to \mathbb{R}^k$ such that $r\in \mathbb{R}^k$ is a regular value, and the submanifold $N = \Psi^{-1}(r)$. Let $(c^{ij})$ be the matrix with entries $c^{ij}=\{\psi_i,\psi_j\}|_N$. (a) Show that $N$ is cosymplectic if and only if the matrix $(c^{ij})$ is invertible. Let $(c_{ij})$ be the inverse matrix. (b)  Prove that the Poisson bracket on $N$ is given by
$$
\{f,g\}_N = (\{\widetilde{f},\widetilde{g}\} -\sum_{i,j} \{\widetilde{f}, \psi_i\}c_{ij}\{\psi_j, \widetilde{g}\})|_N,
$$
where $\widetilde{f}$ and $\widetilde{g}$ are {\em arbitrary} extensions of $f$ and $g$ to $M$ (cf. Exercise~\ref{exer:cosymPoiss} (a)).
\end{exer}

\begin{exer}\label{exer:cosympfactor}
Let $M=S\times N$ be the product of two Poisson manifolds $S$ and $N$, with $S$ symplectic. 
\begin{itemize}
\item[(a)] Show that, for all $x\in S$, $\{x\}\times N$ is a cosymplectic submanifold of $M$ and the Poisson structure it inherits as such agrees with the original Poisson structure on $N$. 
\item[(b)] Check that, for $y\in N$, $S\times \{y\}$ is a Poisson submanifold of $M$ if and only if $\pi_N|_y=0$.
\end{itemize}
\end{exer}

\subsection{Some Poisson invariants}
Let $(M,\pi)$ be a Poisson manifold. We can use the Schouten bracket $[\cdot,\cdot]$ on multivector fields to define an operator
$$
d_\pi:=[\pi,\cdot]: \mathfrak{X}^\bullet(M)\to \mathfrak{X}^{\bullet +1}(M).
$$
By the graded Jacobi identity for the Schouten bracket, we have that $d_\pi^2=\frac{1}{2}[[\pi,\pi], \cdot]=0$ (see Exercise~\ref{exer:schouten}). The cohomology of the cochain complex $(\mathfrak{X}^\bullet(M), d_\pi)$, denoted by $H_\pi^\bullet(M)$,
is called the {\em Poisson cohomology} of $M$.

Recall the bundle map $\pi^\sharp: T^*M\to TM$, $\pi^\sharp(\alpha)=\pi(\alpha, \cdot)$, and consider the induced $C^\infty(M)$-linear homomorphism $(\pi^\sharp)^*: \Omega^\bullet(M)\to \mathfrak{X}^\bullet(M)$ given by
$$
(\pi^\sharp)^*\eta (\alpha_1,\ldots,\alpha_k)=\eta(\pi^\sharp(\alpha_1),\ldots,\pi^\sharp(\alpha_k)),
$$
for $\eta\in \Omega^k(M)$ and $\alpha_i\in\Omega^1(M)$, $i=1,\ldots,k$. 

\begin{exer}\label{exer:dRpoisson}
Show that, for any $\eta \in \Omega^k(M)$, 
$$
(\pi^\sharp)^* d \eta=  d_\pi (\pi^\sharp)^* \eta,
$$ 
i.e.,
$(\pi^\sharp)^*$ is a morphism from the de Rham complex to the Poisson complex. (Exercise~\ref{exer:schouten} checks this equality when $\eta$ has degree 0; the general result can be reduced to this case and the case where $\eta$ is an exact 1-form.)
\end{exer}

It follows that $(\pi^\sharp)^*: \Omega^\bullet(M)\to \mathfrak{X}^\bullet(M)$ induces a homomorphism in cohomology, 
$$
H^\bullet_{dR}(M)\to H^\bullet_\pi(M).
$$
When $\pi$ is nondegenerate (i.e., symplectic), this map is an isomorphism. In general, this map is neither surjective nor injective, and Poisson cohomology tends to be very different from de Rham cohomology.
For example, when $\pi=0$, we have $H_\pi^k(M)=\mathfrak{X}^k(M)$, so Poisson cohomology groups can be infinite-dimensional even when de Rham cohomology groups are trivial.

As with other types of cohomologies, Poisson cohomology groups have meaningful interpretations in low degrees. 

\begin{itemize}
    \item The zeroth Poisson cohomology group  consists of functions $f\in C^\infty(M)$ such that (see Exercise~\ref{exer:schouten})
    $$
    d_\pi f = -X_f =0,
    $$ 
    i.e., $H_\pi^0(M)$ is the space of Casimir functions.

   \item 
    A {\em Poisson vector field} on $(M,\pi)$ is a vector field $X$ such that 
   $$
   \mathcal{L}_X \pi = - d_{\pi} X = 0,
   $$ 
  i.e., $X$ is an infinitesimal automorphism of the Poisson structure (its flow is by  Poisson automorphisms). Equivalently, a Poisson vector field $X$ is a derivation of the Poisson bracket,
   $$
\mathcal{L}_X \{f,g\}= \{\mathcal{L}_Xf, g\} + \{f,\mathcal{L}_Xg\},
   $$
   for all $f,g\in C^\infty(M)$. The space of Poisson vector fields, denoted by $\mathfrak{X}_{\pi}(M)$, contains the space $\mathfrak{X}_{\mathrm{Ham}}(M)$ of Hamiltonian vector fields (cf. Exercise~\ref{exer:jac}), which are the inner derivations of the bracket.   
   
   \begin{exer}
   Note that $\mathfrak{X}_{\pi}(M)$ is closed under the Lie bracket of vector fields.
    Check that, if $X \in  \mathfrak{X}_\pi(M)$  and $X_f \in \mathfrak{X}_{\mathrm{Ham}}(M)$, then $[X,X_f]= X_{\mathcal{L}_Xf} \in \mathfrak{X}_{\mathrm{Ham}}(M)$, and hence   
    $\mathfrak{X}_{\mathrm{Ham}}(M) \subseteq \mathfrak{X}_\pi(M)$ is a Lie ideal.
   \end{exer}

   The first Poisson cohomology group is given by
   $$
   H^1_\pi(M)=\mathfrak{X}_\pi(M)/\mathfrak{X}_{\mathrm{Ham}}(M),
   $$
i.e., it is   the Lie algebra of outer derivations of the Poisson bracket.

   \item The second Poisson cohomology group on $(M,\pi)$ has an interpretation in terms of deformations of the Poisson structure $\pi$. For $\pi_1\in \mathfrak{X}^2(M)$ and denoting by $\varepsilon$ a formal parameter,
   the condition
   $$
d_\pi\pi_1=[\pi,\pi_1] =0
   $$
is equivalent to the property that $\pi+\varepsilon \pi_1$ satisfies the Jacobi identity up to order $\epsilon^2$, 
$$
[\pi+\varepsilon \pi_1, \pi+\varepsilon \pi_1]=0 \mod \varepsilon^2.
$$
If this holds we view $\pi_1$ as defining an infinitesimal deformation of $\pi$.
This infinitesimal deformation is called {\em trivial} if $\pi_1 = d_\pi X = -\mathcal{L}_X\pi$ for a vector field $X$, since in this case $\pi$ and $\pi+\varepsilon \pi_1$ are isomorphic up to order $\varepsilon^2$ via the time-$\varepsilon$ flow of $X$:
$$
(\varphi^\varepsilon_X)_*\pi = \pi + \varepsilon \pi_1 \mod \varepsilon^2.
$$
Therefore $H^2_\pi(M)$ is the space of infinitesimal deformations of $\pi$ modulo trivial deformations, so it could be intuitively regarded as the tangent space at $\pi$ to the moduli space of Poisson structures on $M$ up to Poisson automorphisms. This indicates the importance of degree 2 Poisson cohomology in the study of normal forms of Poisson structures. Note that $\pi$ itself defines a class in $H^2_\pi(M)$; when this class vanishes, the Poisson structure is called {\em exact} (see Exercise~\ref{exer:lineareuler} for an example).
\end{itemize}

Explicit computations of Poisson cohomology groups are notoriously difficult; see \cite[Ch.~2]{DZ} for many of the known results and calculation methods. See also $\S$~\ref{subsec:linear1}, $\S$~\ref{subsec:PL} and $\S$~\ref{subsec:log} for further references in specific examples.

\begin{exer}\label{exer:2dcohomology}
Consider $M=\mathbb{R}^2$ with Poisson structure $\pi= x\partial_x \wedge\partial_y$.
\begin{itemize}
\item[(a)] Verify that the restriction of a Poisson vector field to the $y$-axis is a constant multiple of $\frac{\partial}{\partial y}$, and that a 
Poisson vector field is Hamiltonian if and only if it vanishes on the $y$-axis.

\item[(b)] Show that $H_\pi^0(M)=\mathbb{R}$, $H^1_\pi(M)=\mathbb{R}$ (with generator given by the class of $\frac{\partial}{\partial y}$), and $H^2_\pi(M)=0$.
\end{itemize}
\end{exer}

A Poisson manifold $(M,\pi)$ has a canonical class in $H^1_\pi(M)$, known as the {\em modular class}. Let us assume for simplicity that $M$ is orientable. Given a volume form $\eta \in \Omega^{top}(M)$, we recall that the {\em divergence} of a vector field $Y$ on $M$ with respect to $\eta$ is the function $\mathrm{div}_{\eta}(Y) \in C^\infty(M)$ defined by
$$
\mathcal{L}_Y\eta = \mathrm{div}_{\eta}(Y) \eta.
$$


\begin{exer}\label{exer:modvf}
Check that the map $C^\infty(M) \to C^\infty(M)$, $f\mapsto \mathrm{div}_{\eta}(X_f)$, is a derivation of the pointwise product of functions.
\end{exer}

Therefore each volume form $\eta$ defines a vector field $X_\eta$ on $(M,\pi)$ by the condition
$$
\mathcal{L}_{X_\eta}f = \mathrm{div}_{\eta}(X_f), \ \ \forall f\in C^\infty(M).
$$
We call $X_\eta$ the {\em modular vector field} with respect to $\eta$.

\begin{exer}\label{exer:modprop}
\begin{itemize}
\item[(a)] Show that $X_\eta$ is a Poisson vector field.
\item[(b)] Let $\eta'=f \eta$ be another volume form, for $f \in C^\infty(M)$ nowhere vanishing. Show that $X_\eta - X_{\eta'}= X_{\log |f|}$.
\end{itemize}
\end{exer}

By the previous exercise, $X_\eta$ defines a 1-cocycle in the Poisson complex of $(M,\pi)$, and its cohomology class
$$
[X_\eta] \in H_\pi^1(M)
$$
is independent of the choice of volume form. We refer to $[X_\eta]$ as the {\em modular class} of $(M,\pi)$. A Poisson manifold is called {\em unimodular} when its modular class is trivial.

\begin{exer}\label{exer:unim}
Show that a Poisson manifold is unimodular if and only if there exists a volume form that is invariant by all Hamiltonian vector fields. Conclude that a symplectic manifold is unimodular.
\end{exer}

\begin{exer}\label{exer:modlog}
Show that the modular class of the Poisson manifold in Exercise~\ref{exer:2dcohomology} is nontrivial, given by $[-\frac{\partial}{\partial y}]$.
\end{exer}

See also Exercise~\ref{exer:modchar}.

\begin{remark} We may drop the orientability requirement on $M$ by replacing volume forms by smooth densities in the definition of modular vector fields. \hfill $\diamond$
\end{remark}

The modular class of a Poisson manifold goes back to the work of Koszul, see \cite{KSModular} for original references; it was rediscovered in \cite{WeinMod}, where the terminology was introduced motivated by analogies with the modular theory of von Neumann algebras.

\section{Examples} \label{sec:example}

We have seen some initial examples of Poisson manifolds in $\S$ \ref{subsub:brackets} and $\S$ \ref{subsub:tensor}. We now discuss some broader classes of examples.

\subsection{{Quotients by symmetries}}\label{subsec:symmetry}

Let $(M,\pi)$ be a Poisson manifold carrying an action of a Lie group $G$ by Poisson diffeomorphisms. Then the space of invariant functions on $M$,
$$
C^\infty(M)^G:=\{ f \in C^\infty(M)\,|\, \sigma^*f=f \; \forall\,\sigma\in G\} \subseteq C^\infty(M),
$$ 
is a Poisson subalgebra, i.e.,
$$
\{C^\infty(M)^G, C^\infty(M)^G\}\subseteq C^\infty(M)^G,
$$ 
as shown by the following (more general) exercise.

\begin{exer}\label{exer:symmetry}
Let $K$ and $M$ be Poisson manifolds. An immediate consequence of \eqref{eq:product} is that a map $\psi: K\times M \to M$ is Poisson if and only if
$$
\{f,g\}_M(\psi(k,x)) = \{f(\psi(k,\cdot)),g(\psi(k,\cdot))\}_M(x) +
\{f(\psi(\cdot,x)),g(\psi(\cdot,x))\}_K(k).
$$
(When $\{\cdot,\cdot\}_K\equiv 0$, it follows that $\psi$ is Poisson if and only if $\psi(k,\cdot): M\to M$ is a Poisson map for each $k\in K$.) Show that if $f\circ \psi(k,\cdot) = f$ and $g\circ \psi(k,\cdot) = g$ for all $k\in K$, then $\{f,g\}_M \circ \psi(k,\cdot)= \{f,g\}_M$ for all $k\in K$. 
\end{exer}

Suppose that the $G$-action on $M$ is free and proper, so that $B=M/G$ is a smooth manifold and the quotient map $p: M \to B$ is a surjective submersion. Then $C^ \infty(B)$ inherits a Poisson bracket from the identification $p^*: C^\infty(B)\stackrel{\sim}{\to} C^\infty(M)^G$,
and in this way $B$ acquires a Poisson structure with the property that $p:M\to B$ is a Poisson map. 

A special case of this construction is when $(M,\omega)$ is symplectic and $G$ acts on $M$ by symplectomorphisms. 
For instance,  whenever a Lie group $G$ acts on a manifold $Q$, its cotangent lift\footnote{The cotangent lift of a diffeomorphism $\psi: Q\to Q$ is $(d\psi^{-1})^*: T^*Q\to T^*Q$.} induces a $G$-action on $M=T^*Q$ preserving the canonical symplectic form (cf. Exercise~\ref{exer:linarcanonical} (a)). If the $G$-action on $Q$ is free and proper, so is the lifted action on $T^*Q$, and hence $T^*Q/G$ inherits a Poisson structure. As an example, consider the action of a Lie group $G$ on itself by right multiplication, $\sigma \mapsto (r_{\sigma^{-1}}:G\to G)$, and the lift of this action to $T^*G$, $\sigma \mapsto ((dr_{\sigma})^*: T^*G\to T^*G))$. In this case the map $T^*G\to \mathfrak{g}^*$ defined by right translations gives rise to an identification
of $T^*G/G$ with $\mathfrak{g}^*$.  

\begin{exer}\label{exer:liepois}
Check that the induced Poisson structure on $\mathfrak{g}^*=T^*G/G$ is given by
$$
\{f,g\}(\xi):=  \xi([df|_\xi, dg|_\xi]),
$$
for $f, g \in C^\infty(\mathfrak{g}^*)$ and $\xi \in \mathfrak{g}^*$ (you may get a sign depending on your conventions). (Here $df|_\xi$ and $dg|_\xi$ are viewed as elements in $\mathfrak{g}$  through the identification $T^*_\xi \mathfrak{g}^* = (\mathfrak{g}^*)^*=\mathfrak{g}$.) 
\end{exer}

\subsection{{Linear Poisson structures I: Lie algebras}}\label{subsec:linear1}
Let $V$ be a (real, finite-dimen\-sio\-nal) vector space. Consider $V^*=C^\infty_{lin}(V) \subseteq C^\infty(V)$, the subspace of linear functions on $V$. A Poisson structure on $V$ is called {\em linear} if
$$
\{C^\infty_{lin}(V),C^\infty_{lin}(V)\} \subseteq C^\infty_{lin}(V).
$$
In this case the restriction of $\{\cdot,\cdot\}$  makes $V^*=C^\infty_{lin}(V)$ into a Lie algebra.

Conversely,
let $(\mathfrak{g}, [\cdot,\cdot])$ be a Lie algebra. As a vector space, we can identify it with the space of linear functions on $\mathfrak{g}^*$, 
$$
\mathfrak{g}=C^\infty_{lin}(\mathfrak{g}^*) \subseteq C^\infty(\mathfrak{g}^*).
$$

\begin{exer}
Check that
(a) there exists at most one almost Poisson bracket on $C^\infty(\mathfrak{g}^*)$ such that $\{\cdot,\cdot\}|_{\mathfrak{g}}=[\cdot,\cdot]$, and
(b) an  almost Poisson bracket on $C^\infty(\mathfrak{g}^*)$ with this property is Poisson.
\end{exer}

By the previous exercise, the bracket on $C^\infty(\mathfrak{g}^*)$  given in Exercise \ref{exer:liepois},
$$
\{f,g\}(\xi):=  \xi([df|_\xi, dg|_\xi]),
$$
is the unique Poisson structure on $\mathfrak{g}^*$ satisfying $\{\cdot,\cdot\}|_{\mathfrak{g}}=[\cdot,\cdot]$. This bracket is characterized by the fact that, for $u\in \mathfrak{g}$ (viewed as a linear function on $\mathfrak{g}^*$), the corresponding Hamiltonian vector field on $\mathfrak{g}^*$ is
\begin{equation}\label{eq:Xu}
X_u|_\xi = \xi(\mathrm{ad}_u(\cdot)) = \mathrm{ad}^*_u(\xi) \; \in \mathfrak{g}^*=T_\xi\mathfrak{g}^*.
\end{equation}

In conclusion, we have a natural bijection between linear Poisson structures on a vector space $V=\mathfrak{g}^*$ and Lie algebra structures on its dual $V^*=\mathfrak{g}$.

With respect to linear coordinates $(\xi_1,\ldots,\xi_n)$ on $V$ (corresponding to a basis on $V^*$), a linear Poisson structure is determined by linear functions $\{\xi_i,\xi_j\}=\sum_k c_{ijk}\xi_k$, so it has the form
$$
\pi = \sum_{i<j} (\sum_k c_{ijk} \xi_k) \frac{\partial}{\partial \xi_i}\wedge \frac{\partial}{\partial \xi_j},
$$
and $c_{ijk}$ are the structure constants of the corresponding Lie bracket on $V^*$.

Note that Example~\ref{ex:su2} is the linear Poisson structure on the dual of the Lie algebra $\mathfrak{so}(3)$.

\begin{exer}\label{exer:addition}
Let $\pi$ be a Poisson structure on a vector space $V$. Show that $\pi$ is linear if and only if the addition map $+: V\times V \to V$ is Poisson.
\end{exer}

\begin{exer}\label{exer:liedual} Let $\mathfrak{g}$ and  $\mathfrak{h}$ be Lie algebras. Show that a  linear map $\mathfrak{g}\to \mathfrak{h}$ is a Lie algebra homomorphism if and only if the dual map $\mathfrak{h}^*\to \mathfrak{g}^*$ is Poisson.
\end{exer}

\begin{exer}\label{exer:Lieann}
Let $\mathfrak{g}$ be a Lie algebra and $\mathfrak{h}\subseteq \mathfrak{g}$  be a subspace. Show that $\mathfrak{h}$ is an ideal (resp. Lie subalgebra)  if and only if $\mathrm{Ann}(\mathfrak{h})\subseteq \mathfrak{g}^*$ is a Poisson (resp.  coisotropic)  submanifold. Moreover, when $\mathfrak{h}$ is an ideal, check that the  Poisson structure on $\mathrm{Ann}(\mathfrak{h})$  coincides with the one on the dual of the Lie algebra $\mathfrak{g}/\mathfrak{h}$ under the natural identification $\mathrm{Ann}(\mathfrak{h})=(\mathfrak{g}/\mathfrak{h})^*$.
\end{exer}

\begin{exer}\label{exer:modchar}
Let $\mathfrak{g}$ be a Lie algebra. Show that the modular vector field on $\mathfrak{g}^*$ with respect to any constant volume form
is the constant vector field 
$$
\mathfrak{g}^*\ni \xi \mapsto \chi \in T_\xi\mathfrak{g}^*=\mathfrak{g}^*,
$$ 
where $\chi$ is the {\em modular character} of $\mathfrak{g}$, given by $\chi(u)=\mathrm{Tr}(\mathrm{ad}_u)$. 
\end{exer}

The Poisson cohomology of linear Poisson structures associated to compact Lie algebras is described in \cite{GinzWein}. For a description of the Poisson cohomology of all three-dimensional linear Poisson structures, see \cite{Zeiser}.



The next two exercises concern a generalization of linear Poisson structures.

\begin{exer}\label{exer:2cocycle}
\begin{itemize}
\item[(a)] Let $\mathfrak{g}$ be a Lie algebra  with linear Poisson structure $\pi_{\mathfrak{g}^*}$ on $\mathfrak{g}^*$. Let $\lambda$ be a constant Poisson structure on $\mathfrak{g}^*$ (viewed as an element in $\wedge^2 \mathfrak{g}^*$).
 Show that $\pi_{\mathfrak{g}^*}+\lambda$ is a Poisson structure on $\mathfrak{g}^*$ if and only if
$$
\lambda(u_1, [u_2,u_3]) + \lambda(u_3, [u_1,u_2]) + \lambda(u_2, [u_3,u_1]) =0, \quad \; \forall \, u_1, u_2, u_3\in \mathfrak{g},
$$
i.e., $\lambda$ is a {\em 2-cocycle} on $\mathfrak{g}$.

\item[(b)] A 2-cocycle $\lambda$ defines a Lie algebra structure on $\mathfrak{g}\oplus \mathbb{R}$ via $[(u,t),(v,s)]=([u,v] ,\lambda(u,v))$ (a central extension of $\mathfrak{g}$ by $\mathbb{R}$). Check that $(\mathfrak{g}^*,\pi_{\mathfrak{g}^*}+\lambda)$ sits in $(\mathfrak{g}\oplus \mathbb{R})^* = \mathfrak{g}^*\oplus \mathbb{R}$ as a Poisson submanifold via $\xi \mapsto (\xi,1)$.
\end{itemize}
\end{exer}

\begin{exer}\label{exer:aff}
Let $C^\infty(V)_{aff}\subseteq C^\infty(V)$ denote the subspace of affine functions on a vector space $V$. Show that a Poisson structure on $V$ is {\em affine}, in the sense that $\{C^\infty_{aff}(V),C^\infty_{aff}(V)\} \subseteq C^\infty_{aff}(V)$, if and only if it is of the form $\pi + \lambda$, where $\pi$ is a linear Poisson structure on $V$ and $\lambda\in \wedge^2V$ is a 2-cocycle on the Lie algebra $V^*$. 
\end{exer}

Many important examples of linear and affine Poisson structures on duals of Lie algebras occur in infinite dimensions, providing the Hamiltonian description of various PDEs related e.g. to fluid dynamics, see e.g. \cite{KhesinWendt}.

Linear (and affine) Poisson structures play a central role in the theory of Hamiltonian actions.
Let $M$ be a Poisson manifold and $\mathfrak{g}$ a Lie algebra.
A $\mathfrak{g}$-action (i.e. Lie-algebra homomorphism) $\psi: \mathfrak{g}\to \mathfrak{X}^1(M)$  is called {\em weakly Hamiltonian}
if there exists a linear map 
$\mu^*: \mathfrak{g}\to C^ \infty(M)$ 
such that 
\begin{equation}\label{eq:hammu}
\psi(u)= X_{\mu^*u}, \; \mbox{for all} \; u\in \mathfrak{g}. 
\end{equation}
The $\mathfrak{g}$-action is said to be {\em Hamiltonian} when there is such a $\mu^*$ that is, in addition, a Lie algebra homomorphism. 
 
Any linear map $\mu^*: \mathfrak{g}\to C^\infty(M)$ is equivalent to a smooth map $\mu: M\to \mathfrak{g}^*$ via 
$$
\SP{\mu(x),u}=\mu^*(u)(x),
$$ 
for $x\in M$ and $u\in \mathfrak{g}$. For a weakly Hamiltonian $\mathfrak{g}$-action $\psi: \mathfrak{g}\to \mathfrak{X}^1(M)$, a map $\mu: M\to \mathfrak{g}^*$
for which \eqref{eq:hammu} holds is called a {\em moment map}.

\begin{exer}\label{exer:mmap}
Check that a linear map $\mu^*:\mathfrak{g}\to C^\infty(M)$ is a Lie algebra homomorphism if and only if $\mu: M\to \mathfrak{g}^*$ is a Poisson map.
\end{exer}

An example of a Hamiltonian action is the coadjoint $\mathfrak{g}$-action on $\mathfrak{g}^*$; in this case the identity map $\mathfrak{g}^*\to \mathfrak{g}^*$ is a moment map (the corresponding Lie-algebra homomorphism $\mathfrak{g}\to C^\infty(\mathfrak{g}^*)$ is the natural inclusion of $\mathfrak{g}$ as linear functions on $\mathfrak{g}^*$), and generating vector fields are given in \eqref{eq:Xu}.

\begin{exer}\label{exer:equiv} If $\mu: M\to \mathfrak{g}^*$ is a moment map for a weakly Hamiltonian $\mathfrak{g}$-action,  check that $\mu$ is  a Poisson map if and only if it is $\mathfrak{g}$-equivariant with respect to the coadjoint action (cf. Exercise~\ref{exer:pmaps} (b)).
\end{exer}

For Hamiltonian actions, moment maps will always be assumed to be equivariant (or equivalently, by the previous exercise, to be Poisson maps).
Given a Hamiltonian $\mathfrak{g}$-action $\psi$ on $(M,\pi)$ with moment map $\mu:M\to \mathfrak{g}^*$, the data $(M,\pi,\psi, \mu)$ is called a  {\em Hamiltonian $\mathfrak{g}$-space}. Note that the action $\psi$ is determined by $\mu$ via \eqref{eq:hammu}; on the other hand, any Poisson map $\mu: M\to \mathfrak{g}^*$ defines a Hamiltonian action on $M$ by the same formula \eqref{eq:hammu} (see \eqref{eq:brkpres} and Exercise~\ref{exer:mmap}).
It follows that Hamiltonian $\mathfrak{g}$-spaces are the same as Poisson maps into $\mathfrak{g}^*$,
\smallskip

\begin{center}{ {Hamiltonian $\mathfrak{g}$-spaces \; $\rightleftharpoons$ \; Poisson maps into $\mathfrak{g}^*$}.}
\end{center}

\smallskip

More generally, given a Poisson manifold $M$, a Lie algebra $\mathfrak{g}$, and a smooth map $\mu: M\to \mathfrak{g}^*$, define the skew-symmetric bilinear map $\lambda: \mathfrak{g}\times \mathfrak{g}\to C^\infty(M)$,
$$
\lambda(u,v)= \{\mu^*u,\mu^*v\} - \mu^*[u,v].
$$


\begin{exer}
Show that $\mu$ is the moment map for a weakly Hamiltonian $\mathfrak{g}$-action on $M$ if and only if $\lambda$ takes values in Casimirs of $M$. If $M$ is symplectic and connected, conclude that $\lambda$ defines an element in $\wedge^2 \mathfrak{g}^*$ and verify that it is a 2-cocycle on $\mathfrak{g}$ (cf. Exercise~\ref{exer:2cocycle}).

\end{exer}

Suppose that $(M,\omega)$ is a connected symplectic manifold carrying a weakly Hamiltonian $\mathfrak{g}$-action with moment map $\mu: M\to \mathfrak{g}^*$ and associated 2-cocycle $\lambda \in \wedge^2\mathfrak{g}^*$, as in the previous exercise. Then $\mathfrak{g}^*$ carries an affine Poisson structure $\pi_{\mathfrak{g}^*} + \lambda$, cf. Exercises~\ref{exer:2cocycle} and \ref{exer:aff}, and  an affine $\mathfrak{g}$-action 
$\psi_\lambda: \mathfrak{g}{\to} \mathfrak{X}^1(\mathfrak{g}^*)$, where
$$
\psi_\lambda(u)|_\xi = \mathrm{ad}^*_u(\xi) + i_u\lambda \; \in \, \mathfrak{g}^*=T_\xi \mathfrak{g}^*.
$$
(The affine action $\psi_\lambda$ is simply the Hamiltonian action on $(\mathfrak{g}^*, \pi_{\mathfrak{g}^*} + \lambda)$ with  moment map given by the identity map $\mathfrak{g}^* {\to} \mathfrak{g}^*$, cf. \eqref{eq:Xu}.)

\begin{exer}
Verify that the moment map $\mu$ is a Poisson map $(M,\omega) \to (\mathfrak{g}^*,\pi_{\mathfrak{g}^*}+\lambda)$, noticing that $\mu$ is Poisson if and only if it is $\mathfrak{g}$-equivariant (with respect to the affine action on $\mathfrak{g}^*$). 
\end{exer}

\subsection{{Linear Poisson structures II: Lie algebroids}}\label{subsec:linear2}

The notion of linear Poisson structure makes sense, more generally, on vector bundles. 

For a vector bundle $E\stackrel{q}{\to} M$, denote by $C_{bas}^\infty(E)=q^*C^\infty(M)\subseteq C^\infty(E)$ the subalgebra of basic functions, and by
$C^\infty_{lin}(E)\subseteq C^\infty(E)$ the subspace of fiberwise linear functions.
We use the fact that $C_{bas}^\infty(E)\cdot C^\infty_{lin}(E)\subseteq C^\infty_{lin}(E)$ to regard $C^\infty_{lin}(E)$ as a $C^\infty(M)$-module, so that we have a natural identification $\Gamma(E^*)=C^\infty_{lin}(E)$ (as $C^\infty(M)$-modules).

A Poisson structure on the total space $E$ is called {\em linear} if 
$$
\{C^\infty_{lin}(E),C^\infty_{lin}(E)\} \subseteq C^\infty_{lin}(E).
$$
We denote by
\begin{equation}\label{eq:duallie}
[\cdot,\cdot]: \Gamma(E^*)\times \Gamma(E^*)\to \Gamma(E^*)
\end{equation}
the Lie bracket on $\Gamma(E^*)=C^\infty_{lin}(E)$ obtained by restriction of $\{\cdot,\cdot\}$.

\begin{exer}\label{exer:duallie}
Check that (by the Leibniz rule) a linear Poisson structure on $E\to M$ satisfies
$$
\{C^\infty_{lin}(E),C_{bas}^\infty(M)\}\subseteq C^\infty_{bas}(E)\qquad \{C_{bas}^\infty(E), C_{bas}^\infty(E)\}=0.
$$ 
Show also that there is a bundle map $\rho: E^*\to TM$ such that  
\begin{equation}\label{eq:rho}
\{\xi, q^*f\}=q^*\Lie_{\rho(\xi)}f,
\end{equation} 
for $\xi\in \Gamma(E^*)$ and $f\in C^\infty(M)$, and that
$\rho$ and the Lie bracket $[\cdot,\cdot]$ on $\Gamma(E^*)$ satisfy
$[\xi, f \eta]= (\Lie_{\rho(\xi)}f)\eta + f[\xi,\eta]$, for all $\xi,\eta\in \Gamma(E^*)$ and $f\in C^\infty(M)$.
\end{exer}

From the previous exercise, in local bundle coordinates $(x_1,\ldots,x_n,\xi_1,\ldots,\xi_r)$ on $E$ (defined by a local frame on $E^*$), a linear Poisson structure satisfies
$$
\{x_i,x_j\}=0,\;\; \{x_i,\xi_j\}=-\Lie_{\rho(\xi_j)}x_i =-\rho_{ij}(x),\;\; \{\xi_i,\xi_j\}=\sum_k c_{ijk}(x)\xi_k,
$$
so the corresponding bivector field $\pi\in \mathfrak{X}^2(E)$ locally has the form
$$
\pi = \sum_{i<j}(\sum_k c_{ijk}(x)\xi_k)\frac{\partial}{\partial \xi_i}\wedge \frac{\partial}{\partial \xi_j} - \sum_{i,j}\rho_{ij}(x)\frac{\partial}{\partial x_i}\wedge \frac{\partial}{\partial \xi_j}.
$$

\begin{exer}\label{exer:lineareuler}
Let $\kappa_t: E\to E$ denote the fiberwise scalar multiplication by $t\in \mathbb{R}$, and denote by $\mathcal{E}\in \mathfrak{X}^1(E)$ the Euler vector field (in local bundle coordinates $(x_i,\xi_j)$, $\mathcal{E}=\sum_j \xi_j\frac{\partial}{\partial \xi_j}$).
Show that a Poisson structure $\pi$ on $E$ is linear if and only if $(\kappa_t)_*\pi=t\pi$ for all $t>0$, if and only if $\mathcal{L}_{\mathcal{E}}\pi=-\pi$. 
\end{exer}

A {\em Lie algebroid} is a vector bundle $A\to M$ equipped with a Lie bracket $[\cdot,\cdot]$ on $\Gamma(A)$ and a vector bundle morphism $\rho: A\to TM$ (over the identity map in $M$), called the {\em anchor}, such that
$$
[u,fv]= (\Lie_{\rho(u)}f)v + f[u,v]
$$
for all $u,v \in \Gamma(A)$ and $f\in C^\infty(M)$. One refers to this last condition as the {\em Leibniz identity}. 

\begin{exer}
Show that, on a Lie algebroid, the anchor map $\rho: \Gamma(A)\to \mathfrak{X}^1(M)$ preserves Lie brackets.
\end{exer}

Examples of Lie algebroids include Lie algebras (when $M$ is a point) and tangent bundles $A=TM$ (with Lie bracket given by the usual commutator bracket of vector fields and $\rho$ the identity map) or, more generally, involutive subbundles  $F\subseteq TM$. We will see other examples later.

For now, the point to be made is that, given a linear Poisson structure on $E\to M$, the Lie bracket \eqref{eq:duallie} on $\Gamma(E^*)$ and the map $\rho$ defined in \eqref{eq:rho} make  (by Exercise~\ref{exer:duallie}) $E^*$ into a Lie algebroid. Conversely, a Lie algebroid structure on $A\to M$ gives rise to a linear Poisson structure on $A^*$, uniquely determined by the conditions
\begin{equation}\label{eq:linbrk}
\{u,v\}=[u,v],\qquad \{u,q^*f\}=q^*(\Lie_{\rho(u)}f)
\end{equation}
for $u, v \in \Gamma(A)=C_{lin}^\infty(A^*)$, $f\in C^\infty(M)$, and $q: A^*{\to} M$ the bundle projection.
These constructions are mutually inverse and establish a bijection between linear Poisson structures on vector bundles $E=A^*$ and Lie algebroid structures on $E^*=A$. 

\begin{exer}
Check that conditions \eqref{eq:linbrk} indeed define a unique linear Poisson bracket on $A^*$.
\end{exer}

\begin{exer}\label{exer:dualpoismap}
\begin{itemize}
\item[(a)] Let $A\to M$ and $B\to M$ be Lie algebroids. Show that a vector bundle map $\varphi: A\to B$ (over the identity map) is a Lie algebroid morphism (i.e., it intertwines anchor maps and the induced map $\Gamma(A)\to \Gamma(B)$ preserves brackets) if and only if $\varphi^*: B^*\to A^*$ is a Poisson map.

\item[(b)] Let $B\to M$ be a vector subbundle of a Lie algebroid $A\to M$.
Show that $B$ is a Lie subalgebroid (i.e., $\Gamma(B)\subseteq \Gamma(A)$ is a Lie subalgebra) if and only if $\mathrm{Ann}(B)\subseteq A^*$ is a coisotropic submanifold.
\end{itemize}
\end{exer}

Cf. Exercises~\ref{exer:liedual} and \ref{exer:Lieann}.

\begin{exer}\label{exer:linarcanonical}
\begin{itemize}
\item[(a)] Check that the Poisson bracket defined by the canonical symplectic form on $T^*M$ is linear, and that the dual Lie algebroid structure on $TM$ is the canonical one, i.e., with $\rho=\mathrm{Id}$.
(It follows that $\mathfrak{X}^1(M)$ sits in $C^\infty(T^*M)$ as a Lie subalgebra, and hence any action $\mathfrak{g}\to \mathfrak{X}^1(M) \subseteq C^\infty(T^*M)$ automatically defines a Hamiltonian $\mathfrak{g}$-action on $T^*M$, 
known as its {\em cotangent lift}.) 

\item[(b)] Let $\pi$ be a linear Poisson structure on $E\to M$.
Show that $\pi$ is nondegenerate if and only if the anchor $E^*\to TM$ is an isomorphism, if and only if $E$ is isomorphic to $T^*M$ equipped with its canonical Poisson structure.
\end{itemize}
\end{exer}

\begin{remark}\label{rem:schouten}
For any Lie algebroid $A\to M$ there is a Schouten bracket on $\Gamma(\wedge^\bullet A)$, defined analogously to the one on $\S$ \ref{subsec:schouten} when $A=TM$, extending the Lie bracket on $\Gamma(A)$. When $A = \mathfrak{g}$ is a Lie algebra, the bracket on $\wedge^\bullet \mathfrak{g}$ agrees with the one induced from the Schouten bracket of left-invariant multivector fields on a Lie group with Lie algebra $\mathfrak{g}$. \hfill $\diamond$
\end{remark}

\subsection{{Poisson-Lie groups and their homogeneous spaces}}\label{subsec:PL}
 Let $G$ be a Lie group. A bivector field $\pi$ on $G$ is called {\em multiplicative} if 
$$
\pi_{\sigma_1\sigma_2}=(l_{\sigma_1})_* \pi_{\sigma_2} + (r_{\sigma_2})_* \pi_{\sigma_1}, \;\;\forall \sigma_1, \sigma_2\in G,
$$
where $l_\sigma, r_\sigma: G\to G$ denote left and right multiplication by $\sigma\in G$.
It follows that $\pi_e=0$, so any nontrivial example of multiplicative bivector field has varying rank. 

\begin{exer}
Check that a Poisson structure $\pi$ on $G$ is multiplicative if and only if 
the multiplication map $m: G\times G \to G$ is a Poisson map and that, in this case, the inversion map is anti-Poisson (i.e., it is a Poisson map from $(G,\pi)$ to $(G,-\pi))$). 
\end{exer}

 A Lie group equipped with a multiplicative Poisson structure is called a {\em Poisson Lie group}.

\begin{exer}
A bivector field $\pi$ on a Lie group $G$ is called {\em affine} if 
$$
\pi_{\sigma_1\sigma_2}=(l_{\sigma_1})_* \pi_{\sigma_2} + (r_{\sigma_2})_* \pi_{\sigma_1} - (r_{\sigma_2})_*(l_{\sigma_1})_* \pi_e, \;\;\forall \sigma_1, \sigma_2\in G,
$$
so that $\pi$ is multiplicative if and only if it is affine and $\pi_e=0$. Prove that if $\pi$ is affine, then $\Lie_{u^r}\pi$ is right invariant for all $u\in \mathfrak{g}$ (equivalently, $\Lie_{u^l}\pi$ is left invariant for all $u\in \mathfrak{g}$), and the converse holds if $G$ is connected.

\end{exer}

Any Lie group is a Poisson Lie group with the trivial Poisson structure. 
Multiplicative Poisson structures on vector spaces, viewed as abelian groups, are the same as linear Poisson structures (see Exercise~\ref{exer:addition}).

\begin{exer}\label{exer:semidir}
Let $G=\mathbb{R}\ltimes \mathbb{R}$ be the semi-direct product Lie group with respect to the action of $\mathbb{R}$ on itself via $y\mapsto e^x y$; i.e., $G$ is $\mathbb{R}^2=\{(x,y)\}$ as a manifold with group structure $(x_1,y_1)\cdot (x_2,y_2)= (x_1+x_2, y_1+ e^{x_1}y_2) $. For $\varphi\in C^\infty(\mathbb{R}^2)$, check that the Poisson structure $\varphi(x,y)\partial_x\wedge \partial_y$ (see Exercise~\ref{exer:2d}) is multiplicative if and only if $\varphi(x_1+x_2, y_1+e^{x_1}y_2)= \varphi(x_1,y_1)+ e^{x_1}\varphi(x_2,y_2)$. In particular, each of the Poisson structures $(e^x-1)\partial_x\wedge \partial_y$ and $y \partial_x\wedge \partial_y$  makes $\mathbb{R}\ltimes \mathbb{R}$ into a Poisson Lie group.

\end{exer}

An important class of examples of Poisson Lie groups arises as follows.
Let $G$ be a Lie group with Lie algebra $\mathfrak{g}$.
For $\mathfrak{r} \in \wedge^2\frakg$, consider the
bivector field on $\pi$ on $G$ given by
\begin{equation}\label{eq:r}
\pi_\sigma = (r_\sigma)_*\mathfrak{r} - (l_\sigma)_*\mathfrak{r} = (\mathfrak{r}^r - \mathfrak{r}^l)|_\sigma.
\end{equation}

The next exercise explains when such $\pi$ makes $G$ into a Poisson Lie group. We will make use of the Schouten bracket $[\cdot,\cdot]$ on $\wedge^\bullet \mathfrak{g}$, see Remark~\ref{rem:schouten}

\begin{exer}
\begin{itemize}
\item[(a)] Check that any bivector field $\pi$ of the form \eqref{eq:r} is multiplicative.
\item[(b)] Verify that $[\pi,\pi]=0$ if and only if
$[\mathfrak{r},\mathfrak{r}]^r=[\mathfrak{r},\mathfrak{r}]^l$, i.e.,
$[\mathfrak{r},\mathfrak{r}]$ is $\Ad$-invariant.
\end{itemize}
\end{exer}

When $\pi=\mathfrak{r}^r-\mathfrak{r}^l$ is a Poisson structure, i.e.,
$[\mathfrak{r},\mathfrak{r}]$ is $\Ad$-invariant, $\mathfrak{r}$ is called
an {\em r-matrix}  and the corresponding
Poisson-Lie group is called {\em exact} or {\em coboundary}; when $\mathfrak{r}$ satisfies
$[\mathfrak{r},\mathfrak{r}]=0$ (known as the classical Yang-Baxter equation), $\mathfrak{r}$ is called a {\em triangular r-matrix}. If $G$ is connected and semisimple, or if $G$ is compact, any multiplicative Poisson structure is defined by an $r$-matrix (see \cite[Thm.~1.11]{LuWe}).

\begin{example}\label{ex:SU2}
Let $G$ be the Lie group  
$$
SU(2)= \left \{ \left (\begin{matrix}
a &{b}\\
 -\overline{b} & \overline{a}
\end{matrix}
\right ),\; a= x+iy,\, b= z+i w \in \mathbb{C}, \, |a|^2+|b|^2=1 \right \}.
$$
Its Lie algebra $\frakg=\mathfrak{su}(2)$ has a basis $\{e_1, e_2, e_3\}$, where
\begin{equation}\label{eq:e123}
e_1=\frac{1}{2}\left (
\begin{matrix}
i & 0\\
0 & -i
\end{matrix}
\right ),\;\; e_2=\frac{1}{2}\left (
\begin{matrix}
0 & 1\\
-1 & 0
\end{matrix}
\right ),\;\; e_3=\frac{1}{2}\left (
\begin{matrix}
0 & i\\
i & 0
\end{matrix}
\right ),
\end{equation}
satisfying $[e_1,e_2]=e_3$, $[e_3,e_1]=e_2$, $[e_2,e_3]=e_1$. Notice that $e_1\wedge e_2\wedge e_3$ is $\Ad$-invariant; since it generates the 1-dimensional  space $\wedge^3 \frakg$,
any $\mathfrak{r} \in \frakg\wedge\frakg$ is such that $[\mathfrak{r},\mathfrak{r}]$ is a multiple of $e_1\wedge e_2\wedge
e_3$, and hence $\Ad$-invariant as well. The choice
$\mathfrak{r} = 2 e_2\wedge e_3$ defines the so-called ``standard'' Poisson structure on $SU(2)$.

\hfill $\diamond$
\end{example}

The next exercise relates the standard Poisson structure on $SU(2)$ with the Poisson structure on spheres given in Exercise~\ref{exer:spheres}. 

\begin{exer} \label{exer:SU2}
Consider $\mathbb{S}^3 = \{(x,y,z,w) \in \mathbb{R}^4\, |\, x^2+y^2+z^2+w^2=1\}$ with the Poisson structure $\pi$ defined by the bracket relations
$$
\{x,y\}= z^2 + w^2,\; \{x,z\}= -yz,\; \{x,w\}= -yw, \; \{y,z\}= xz,\; \{y,w\}=wx, \; \{z,w\}=0,
$$
as in Exercise \ref{exer:spheres}, and Lie group structure coming from the quaternionic multiplication in $\mathbb{R}^4=\{x + y\mathbf{i} + z\mathbf{j} + w \mathbf{k}\}$.
With respect to the usual identification of Lie groups $SU(2)=\mathbb{S}^3$, 
$$
\left (\begin{matrix}
x+iy &{z+iw}\\
 -z+iw & x-iy
\end{matrix}
\right ) \mapsto  x + y\mathbf{i} + z\mathbf{j} + w \mathbf{k},
$$
we have that $e_1=\mathbf{i}/2$, $e_2=\mathbf{j}/2$ and $e_3=\mathbf{k}/2$ (see \eqref{eq:e123}).
Verify that
$$
2(e_2^r\wedge e_3^r - e_2^l\wedge e_3^l) = \pi,
$$
i.e., $\pi$ agrees with the standard Poisson structure.
\end{exer}

A multiplicative Poisson structure on a Lie group $G$ admits an infinitesimal des\-cription in terms of the Lie algebra $\mathfrak{g}$. Any bivector field $\pi$ on $G$ such that $\pi_e=0$ gives rise to a map
\begin{equation}\label{eq:cobrk}
\delta: \mathfrak{g}\to \wedge^2\mathfrak{g},\;\;\; u \mapsto (\Lie_X\pi)_e,
\end{equation}
where $X$ is any vector field on $G$ satisfying $X|_e=u$ (check that this is well defined).
Let $[\cdot,\cdot]_*: \wedge^2\mathfrak{g}^*\to \mathfrak{g}^*$ be the dual map. One can prove that (see \cite[Thm.~1.6]{LuWe})
\begin{itemize}
\item if $\pi$ is multiplicative, then 
\begin{equation}\label{eq:compat}
\delta([u,v])=[\delta u,v]+[u,\delta v], \qquad \forall \, u,v\in \mathfrak{g};
\end{equation}
\item if $\pi$ is Poisson, then $[\cdot,\cdot]_*$ defines a Lie bracket on $\mathfrak{g}^*$.
\end{itemize}

\begin{exer}\label{exer:linearcobr}
Let $V$ be a vector space equipped with a linear Poisson structure $\pi$, regarded as an abelian Poisson Lie group. Check that the map $\wedge^2 V^*\to V^*$ dual to \eqref{eq:cobrk} coincides with the Lie bracket on $V^*$ corresponding to $\pi$ (in the sense of $\S$ \ref{subsec:linear1}).
\end{exer}

A {\em Lie bialgebra} is a pair of Lie algebras $(\mathfrak{g},[\cdot,\cdot])$ and $(\mathfrak{g}^*,[\cdot,\cdot]_*)$ such that the map $\delta$ dual to $[\cdot,\cdot]_*$, referred to as the {\em cobracket}, satisfies the compatibility condition \eqref{eq:compat}. We will use the notation  $(\mathfrak{g},\mathfrak{g}^*)$ or $(\mathfrak{g},\delta)$ for a Lie bialgebra. A  {\em Lie-bialgebra morphism} from $(\mathfrak{g}_1,\mathfrak{g}_1^*)$ to 
$(\mathfrak{g}_2,\mathfrak{g}_2^*)$ is a Lie algebra morphism $T: \mathfrak{g}_1\to \mathfrak{g}_2$ satisfying one (and hence all) of the following equivalent conditions: $T$ intertwines $\delta_1$ and $\delta_2$;  $T^*: \mathfrak{g}_2^*\to \mathfrak{g}_1^*$ is a Lie algebra morphism;  $T$ is a Poisson map (with respect to the Poisson structures on $\mathfrak{g}_1$ and $\mathfrak{g}_2$ dual to the Lie brackets on $\mathfrak{g}^*_1$ and $\mathfrak{g}^*_2$, respectively).

As explained above, any Poisson Lie group $G$ gives rise to a Lie bialgebra via \eqref{eq:cobrk}. Conversely, any Lie bialgebra $(\mathfrak{g},\delta)$ arises in this way  from a unique multiplicative Poisson structure on a connected, simply-connected Lie group $G$ integrating $\mathfrak{g}$ \cite[Thm.~1.8]{LuWe}.

A key fact about Lie bialgebras is that if $(\mathfrak{g},\delta)$ is a Lie bialgebra, then so is $(\mathfrak{g}^*,\delta_*)$ (where $\delta_*$ is the cobracket corresponding to the Lie bracket on $\mathfrak{g}$), called the Lie bialgebra {\em dual} to $(\mathfrak{g},\delta)$. For a Poisson Lie group $G$ with Lie bialgebra $(\mathfrak{g}, \delta)$, 
the connected, simply-connected Lie group integrating the dual Lie bialgebra $(\mathfrak{g}^*,\delta_*)$ is a Poisson Lie group $G^*$ called the {\em dual} of $G$. (More generally, we refer to any Poisson Lie group integrating the dual Lie bialgebra $(\mathfrak{g}^*,\delta_*)$ as a {\em dual} to $G$.)

As an example, for a Lie group $G$, regarded as a Poisson Lie group with the trivial Poisson structure, the dual Poisson Lie group is $\mathfrak{g}^*$ (see Exercise~\ref{exer:linearcobr}).

\begin{exer}
Consider the Lie group $G=\mathbb{R}\ltimes \mathbb{R}$ of Exercise \ref{exer:semidir}. Check that the Poisson Lie groups $(G,y\partial_x\wedge \partial_y)$ and $(G, (e^x-1)\partial_x\wedge \partial_y)$ are both self-dual, i.e., in each case there is an isomorphism $(\mathfrak{g},\mathfrak{g}^*)\cong (\mathfrak{g}^*,\mathfrak{g})$.
\end{exer}

\begin{exer}\label{exer:susb}
Consider $SU(2)=\mathbb{S}^3$ with the standard Poisson structure (see Exercise~\ref{exer:SU2}). Show that the dual Lie algebra $\mathfrak{su}(2)^*$ is isomorphic to $\mathbb{R}^3$ with Lie bracket defined by
$$
[e_1,e_2]= e_2 ,\;\;[e_1,e_3]= e_3,\;\; [e_2,e_3]= 0,
$$
called the ``book'' Lie algebra $\mathfrak{sb}(2,\mathbb{C})$ (since its regular coadjoint orbits look like the pages of an open book, see Exercise \ref{exer:g*leaves} (c)). In particular, as a Lie group, the dual of $SU(2)$ can be identified with $SB(2,\mathbb{C})=\left \{ \left (\begin{matrix}
a &{b + ic}\\
 0 & a^{-1}
\end{matrix}
\right ),  a>0,\, b,c \in \mathbb{R}\right\}$.
\end{exer}

If $(G,\pi_G)$ is a Poisson Lie group and $(M,\pi)$ is a Poisson manifold, we say that an action of $G$ on $M$ is {\em Poisson} if the action map $\Psi: G\times M\to M$, $(\sigma,x)\mapsto \sigma.x$, is a Poisson map; equivalently,
\begin{equation}\label{eq:poissact}
\pi_{\sigma . x}=(\Psi_\sigma)_*\pi_x + (\Psi_x)_*(\pi_G)_\sigma,
\end{equation}
 where $\Psi_\sigma=\Psi(\sigma,\cdot): M \to M$ and $\Psi_x=\Psi(\cdot,x): G\to M$, for all $\sigma \in G$, $x\in M$.
Note that $\pi$ is generally not $G$-invariant (cf. Exercise \ref{exer:notinv} below), but it directly follows from \eqref{eq:poissact} that 
if a $G$-action on $M$ is Poisson for two Poisson structures $\pi_1$ and $\pi_2$, then $\pi_1-\pi_2$ is a $G$-invariant bivector field.  

Even though a Poisson action is not an action by Poisson diffeomorphisms, 
it still follows from Exercise \ref{exer:symmetry} that the space of $G$-invariant functions $C^\infty(M)^G\subseteq C^\infty(M)$ is a Poisson subalgebra. Hence if the action is free and proper, the orbit space $M/G$ acquires a Poisson structure for which the quotient map $M\to M/G$ is Poisson (generalizing the discussion in $\S$ \ref{subsec:symmetry}). The theory of Hamiltonian actions and moment maps extends to Poisson actions of Poisson Lie groups -- but in this more general setting moment maps are Poisson maps into dual Poisson Lie groups \cite{Lumoment}.


\begin{exer}\label{exer:notinv}
Let $G$ be a Poisson Lie group with Lie bialgebra $(\mathfrak{g},\delta)$. Suppose that $G$ acts on a Poisson manifold $(M,\pi_M)$ with infinitesimal action denoted by $\psi: \mathfrak{g}\to \mathfrak{X}(M)$. Verify that if the $G$-action on $M$ is Poisson, then
$\Lie_{\psi(u)}\pi_M = \psi(\delta(u))$ for all $u\in \mathfrak{g}$, and that the converse holds if $G$ is connected.
\end{exer}




A class of examples of Poisson manifolds closely related to Poisson Lie groups is given by their Poisson homogeneous spaces \cite{Drinfeld1}: For a Poisson Lie group $G$, a {\em Poisson homogeneous $G$-space} is a Poisson manifold equipped with a transitive Poisson $G$-action. 

Any manifold $M$ carrying a transitive action of a Lie group $G$ is of the form $G/H$ for a closed subgroup $H$ (we consider $H$ acting on $G$ by right multiplication). So a Poisson homogeneous $G$-space is equivalently described by a closed subgroup $H\subseteq G$ and a Poisson structure on $G/H$ for which the map
$$
G \times G/H \to G/H, \qquad (\sigma_1, \sigma_2 H)\mapsto \sigma_1\sigma_2 H
$$
is a Poisson map. A special situation is when $H$ is a Poisson subgroup: in this case the $H$-action on $G$ is Poisson and the Poisson structure on $G$ descends to a homogeneous Poisson structure on $G/H$. In general, if a Poisson structure $\pi$ makes $G/H$ into a Poisson homogeneous space, it is {\em not} true that the quotient map $G\to G/H$ is a Poisson map. The next exercise explains when this happens.

\begin{exer}\label{exer:GmodH}
Let $H$ be a closed subbgroup of a Poisson Lie group $G$. 
\begin{itemize}
\item[(a)] Suppose that $H$ has the property that the Poisson bracket of any two $H$-invariant functions on $G$ vanishes on $H$. Use the multiplicativity of the Poisson structure on $G$ to show that if $f$ and $g$ are $H$-invariant, then so is $\{f,g\}$.

\item[(b)] Show that $G/H$ carries a Poisson structure for which $G\to G/H$ is a Poisson map if and only if $H$ is coisotropic. Moreover, in this case $G/H$ is a Poisson homogeneous space.
\end{itemize}
\end{exer}

\begin{exer}\label{exer:Hcoiso}
Let $G$ be a Poisson Lie group with Lie bialgebra $(\mathfrak{g},\mathfrak{g}^*)$, let
$H \subseteq G$ be a Lie subgroup. Show that if $H$ is coisotropic, then $\mathrm{Ann}(\mathfrak{h})\subseteq \mathfrak{g}^*$ is a Lie subalgebra; similarly, if $H$ is Poisson then $\mathrm{Ann}(\mathfrak{h})$ is a Lie ideal of $\mathfrak{g}^*$ (cf. Exercise~\ref{exer:Lieann}). Moreover, the converses of both statements hold when $H$ is connected.
\end{exer}

\begin{exer}\label{exer:affinehomog}
For a Poisson Lie group $(G,\pi_G)$, let $\pi$ be a Poisson homogeneous structure on $G$ itself, with respect to the action by left multiplication  ($H=\{e\}$). Show that
$\pi-\pi_G= (\pi_e)^l$, and that $\pi$ is an affine Poisson structure on $G$ (see Exercise~\ref{exer:aff}). 

(More is true: given any affine Poisson structure $\pi$ on $G$, $\pi-(\pi_e)^l$ is a multiplicative Poisson structure on $G$ with respect to which $(G,\pi)$ is a Poisson homogeneous space.)
\end{exer}


Hence affine Poisson structures and quotients of Poisson Lie groups by closed coisotropic subgroups provide many examples of Poisson homogeneous spaces.

\begin{remark}[Drinfeld's classification]\label{rem:double}
Given a Lie bialgebra $(\mathfrak{g},\mathfrak{g}^*)$, there is a natural Lie (bi)algebra structure on $\mathfrak{d}=\mathfrak{g}\oplus \mathfrak{g}^*$, known as the {\em Drinfeld double}; the Lie bracket on $\mathfrak{d}$ is characterized by the properties that its natural symmetric pairing is ad-invariant and $\mathfrak{g}$ and $\mathfrak{g}^*$ are Lie subalgebras. For a closed Lie subgroup $H$ of a Poisson Lie group $G$, Poisson homogeneous structures on $G/H$ are classified by suitable lagrangian subalgebras of $\mathfrak{d}$ \cite{Drinfeld2}, see also \cite{MeinLu,MeinDirac}; for coisotropic subgroups, the corresponding lagrangian subalgebras are of the form $\mathfrak{h}\oplus \mathrm{Ann}(\mathfrak{h})$, see Exercise~\ref{exer:Hcoiso}. \hfill $\diamond$
\end{remark}

As an explicit example of Poisson homogeneous space, consider the standard Poisson structure $\pi$ on $SU(2)=\mathbb{S}^3$ described in Exercise~\ref{exer:SU2}. The subgroup $H=\{(x,y,0,0)\,|\, x^2+y^2=1\} \subseteq \mathbb{S}^3$ is such that $\pi |_H=0$, so $H$ is a Poisson subgroup. The quotient space of $\mathbb{S}^3$ by the action of $H=\mathbb{S}^1$ by right multiplication is $\mathbb{C}P^1=\mathbb{S}^2 \subseteq \mathbb{R}^3=\{(x_1,x_2,x_3)\}$, with quotient map $\mathbb{S}^3\to \mathbb{S}^2$ given by the Hopf map
$$
x_1 = 2(xz-yw),\quad x_2= 2(yz+xw), \quad x_3= x^2+y^2-z^2-w^2.
$$
As explained above,  $\mathbb{S}^2$ inherits a quotient Poisson structure that makes it into a Poisson homogeneous space.
We describe it in the  next exercise.

\begin{exer}
Show that 
$$
\{x_1,x_2\}= 2(x_3-1)x_3 ,\quad \{x_2,x_3\}=2(x_3-1) x_1, \quad \{x_1,x_3\}= -2(x_3-1)x_2.
$$
\end{exer}

Following Exercise~\ref{exer:S2}, we see that this Poisson structure on $\mathbb{S}^2$ agrees with 
\begin{equation}\label{eq:bruhS2}
2(1-x_3)\pi_{area},
\end{equation}
where $\pi_{area}$ is the Poisson structure corresponding to the area form. So it is nondegenerate at all points except for $(0,0,1)$, where it vanishes. The next exercise provides suitable coordinates around this singular point.

\begin{exer}
\begin{itemize}
\item[(a)] Consider the stereographic projection (from the south pole $\mathbf{s}=(0,0,-1)$), 
$$
\mathbb{S}^2\setminus \mathbf{s} \to \mathbb{R}^2, \; (x_1,x_2,x_3) \mapsto \left (y_1=\frac{x_1}{1+x_3}, y_2=\frac{x_2}{1+x_3} \right),
$$ 
with inverse 
$$
x_1=\frac{2y_1}{(1+y_1^2+y_2^2)}, \; x_2=\frac{2y_2}{(1+y_1^2+y_2^2)}, \; x_3=\frac{(1-(y_1^2+y_2^2))}{(1+y_1^2+y_2^2))}.
$$ 
Show that $\{y_2,y_1\}=(y_1^2+y_2^2)(1+y_1^2+y_2^2)$.

\item[(b)] With the additional change of coordinates given by $z_1= y_1/\sqrt{1+y_1^2+y_2^2}$ and  $z_2= y_2/\sqrt{1+y_1^2+y_2^2}$, show that the Poisson structure becomes
$$
(z_1^2+ z_2^2)\partial_{z_2}\wedge \partial_{z_1}.
$$
\end{itemize}
\end{exer}

\begin{exer} \label{exer:SU2cov}
Show that any Poisson structure on $\mathbb{S}^2=\{(x_1,x_2,x_3)\,|\, x_1^2+x_2^2+x_3^2=1\}$ that makes it into a Poisson homogeneous $SU(2)$-space is of the form $2(\lambda-x_3)\pi_{area}$, for $\lambda\in \mathbb{R}$.
\end{exer}

The Poisson cohomology of the family of Poisson structures in the previous exercise is computed in \cite[Thm.~3.2]{RoySU}.

The Poisson structure \eqref{eq:bruhS2} on $\mathbb{S}^2$, obtained as a quotient of the standard Poisson structure on $SU(2)$, is an example of a {\em Bruhat} Poisson structure. Much more generally \cite[Thm.~4.7]{LuWe},  any connected, compact, semi-simple Lie group carries a nontrivial ``standard'' Poisson structure making it into a Poisson Lie group, and each of its coadjoint orbits carries a ``Bruhat'' Poisson structure making it into a Poisson homogeneous space.

Further details on this subsection can be found e.g. in \cite[Chp.~1]{KoroSoib}.

\subsection{Log-symplectic manifolds and symplectic Lie algebroids}\label{subsec:log} \
Log-symplec\-tic structures  are special types of Poisson structures that are very close to being symplectic, in the sense that they are generically nondegenerate. Their systematic study began in \cite{GMP,Radko}.

Let $M$ be a $2n$-dimensional manifold. Any bivector field $\pi$ defines a section $\pi^n:=\wedge^n\pi$ of the line bundle $\wedge^{2n} TM$; note that $\pi^n$ is nonzero at points where $\pi$ is nondegenerate, and it vanishes at points where the rank of $\pi$ is less than $2n$.

A Poisson structure $\pi$ on $M$ is called {\em log-symplectic} (or {\em b-symplectic}) if $\pi^n$ vanishes transversally, i.e., $\pi^n$ is transverse to the zero section of $\wedge^{2n} TM$.
It follows that the zero set $Z:=(\pi^n)^{-1}(0)$ is a hypersurface (i.e., a codimension 1 submanifold, possibly disconnected) of $M$, called the {\em singular} or {\em degeneracy locus} of $\pi$.

In local coordinates $(x_1,\ldots,x_{2n})$ on $U\subseteq M$, we have $\pi^n= f \partial_{x_1}\wedge\ldots\wedge \partial_{x_{2n}}$ for a smooth function $f$ on $U$, and $f^{-1}(0)=Z\cap U$. In these coordinates the transversality condition is equivalent to $0$ being a regular value of $f$.


The next exercise shows how 2-dimensional log-symplectic structures look like around points in the degeneracy locus.

\begin{exer}\label{exer:log2dim}
\begin{itemize}
\item [(a)] Consider a nondegenerate Poisson structure on $\mathbb{R}^2$. 
Let $f\in C^\infty(\mathbb{R}^2)$ be such that $f(m)=0$ and $df|_m\neq 0$. Prove that, in a neighborhood of $m$, one can find local coordinates $(y_1,y_2)$, such that $y_1=f$ and $\{y_1,y_2\}=1$. (Use the straightening theorem for $X_f$.) 

\item [(b)] Let $\pi$ be a log-symplectic structure on $\mathbb{R}^2$ with degeneracy locus $Z$.
 Show that any $m\in Z$ admits a neighborhood with coordinates $(y_1,y_2)$ such that $\pi = y_1\partial_{y_1}\wedge \partial_{y_2}$. (Write $\pi= f \overline{\pi}$, where $\overline{\pi}$ is nondegenerate, and use (a).)

\end{itemize}

\end{exer}

More generally, as we will see (Exercise~\ref{exer:log}), given a log-symplectic structure $\pi$ on a $2n$-dimensional manifold $M$ with degeneracy locus $Z\subseteq M$, any point in $Z$ admits a neighborhood $U$ with coordinates $(y_1,y_2,q_1,\ldots,q_{n-1}, p_1,\ldots,p_{n-1})$ such that 
\begin{equation}\label{eq:pilog}
\pi|_U = y_1\partial_{y_1}\wedge\partial_{y_2} + \sum_{i=1}^{n-1}\partial_{p_i}\wedge \partial_{q_i}.
\end{equation}

Since  $U\cap Z=\{ (0, y_2, q_1,\ldots,q_{n-1}, p_1,\ldots,p_{n-1})\}$, it is clear that $\pi$ is tangent to $Z$, i.e., $Z$ is a Poisson submanifold, and the rank of the restricted Poisson structure on $Z$ is $2n-2$. Away from the degeneracy locus, $\pi$ is symplectic; in the coordinates above, the corresponding symplectic form is
$$
\frac{1}{y_1}dy_2\wedge dy_1 + \sum_{i=1}^{n-1}dq_i\wedge dp_i
= dy_2\wedge d(\log|y_1|) + \sum_{i=1}^{n-1}dq_i\wedge dp_i,
$$
which is a 2-form with logarithmic singularities along $Z$.

Log-symplectic manifolds have been the object of much study is recent years. Existence results, constructions  and obstructions for log-symplectic structures are discussed e.g. in \cite{Gil,FrMiMa}.  A full set of invariants classifying them in dimension 2 was described in \cite{Radko} (partially generalized in \cite{GMP} to higher dimensions). Their Poisson cohomology is described in \cite{Radko}  and \cite{GMP,MarOso}.



A fruitful approach to Poisson structures that are ``generically symplectic'' relies on viewing them as symplectic structures on suitable Lie algebroids \cite{klaasse-divisortype,MirSc}. We will briefly recall this viewpoint to log-symplectic structures, following \cite{GMP}.

As seen in $\S$ \ref{subsec:linear2}, tangent bundles of manifolds are special cases of Lie algebroids. Let us consider a Lie algebroid $A\to M$, with anchor $\rho:A\to TM$ and Lie bracket $[\cdot,\cdot]$ on $\Gamma(A)$, and think of it as a replacement for $TM$. Then sections of  $\wedge^k A$ are regarded as analogues of $k$-vector fields. The space of {\em $A$-multivector fields} $\Gamma(\wedge^\bullet A)$ inherits a Schouten bracket extending the Lie bracket on $\Gamma(A)$ (see Remark~\ref{rem:schouten}) and the action of $A$-vector fields $u \in \Gamma(A)$ on functions $f\in C^\infty(M)$ via $\Lie_{\rho(u)}f$. Similarly, the space of $A$-forms $\Gamma(\wedge^\bullet A^*)$ carries a differential 
\begin{equation}\label{eq:dA}
d_A: \Gamma(\wedge^\bullet A^*) \to \Gamma(\wedge^{\bullet+1} A^*)
\end{equation}
satisfying 
$$
d_Af (u) = \Lie_{\rho(u)}f, \qquad d_A\xi(u,v)= \Lie_{\rho(u)}\xi(v)-\Lie_{\rho(v)}\xi(u) - \xi([u,v]),
$$
for $f\in C^\infty(M)$, $u,v \in \Gamma(A)$, and $\xi\in \Gamma(A^*)$. The corresponding cohomology is called the {\em Lie algebroid cohomology} of $A$.

An {\em $A$-Poisson structure} is an $A$-bivector field $\pi_A\in \Gamma(\wedge^2A)$ such that $[\pi_A,\pi_A]=0$. Since the anchor map preserves Lie brackets, it extends to a bracket-preserving map $\rho: \Gamma(\wedge^\bullet A)\to \mathfrak{X}^\bullet(M)$; as a consequence, 
\begin{equation}\label{eq:rhopi}
\pi:=\rho(\pi_A) \in \mathfrak{X}^2(M)
\end{equation}
is an ordinary Poisson structure on $M$. A situation of interest is when $\pi_A$ is nondegenerate, in the sense that it induces an isomorphism $A^*\to A$, or, equivalently, $\pi_A^n$ never vanishes. In analogy with the discussion in $\S$ \ref{subsec:char}, $\pi_A$ is nondegenerate if and only if it is the inverse of an {\em $A$-symplectic structure}, i.e.,  a $d_A$-closed 2-form $\omega_A\in \Gamma(\wedge^2 A^*)$ that is nondegenerate, in the sense that  $\omega_A^\flat: A\to A^*$ is an isomorphism. 
Log-symplectic structures are examples of Poisson structures defined from $A$-symplectic structures via \eqref{eq:rhopi}, for an appropriate choice of Lie algebroid $A$.

Given a manifold $M$ and a hypersurface $Z \hookrightarrow M$,  denote by $\mathfrak{X}^1_Z(M)$ the space of vector fields on $M$ that are tangent to $Z$. 
In a neighborhood $U$ of a point in $Z$ with adapted coordinates $(z, x_2,\ldots,x_l)$ (so that $U\cap Z$ is defined by $z=0$), a vector field tangent to $Z$ has the form
$X= a\partial_z + \sum_i b_i \partial_{x_i}$, for $a, b_i \in C^\infty(U)$ and $a(0,x_1,\ldots,x_l) =0$. Then $a = z b$, for $b\in C^\infty(U)$, and $X= b z\partial_z + \sum_i b_i \partial_{x_i}$, showing that the vector fields 
$$
e_1=z\partial_z, \;\; e_2=\partial_{x_2}, \; \ldots, \;  e_l=\partial_{x_l}
$$ 
form a basis  for the $C^\infty(U)$-module of vector fields on $U$ that are restrictions of vector fields in $\mathfrak{X}^1_Z(M)$. In other words, as a sheaf of $C^\infty(M)$-modules, $\mathfrak{X}^1_Z(M)$ is locally freely generated, and therefore there is a vector bundle $T_ZM$ over $M$ such that $\mathfrak{X}^1_Z(M)= \Gamma(T_ZM)$ \cite{Melrose}.

\begin{exer}
Take $M=\mathbb{S}^1$ and $Z$ to be a point. Check that $T_ZM$ is a nontrivial line bundle
(hence isomorphic to the M\"obius line bundle over $\mathbb{S}^1$), as opposed to $TM$, which is  trivial. What if $Z$ is the union of 2 points?
\end{exer}

The Lie bracket of vector fields defines a Lie bracket on $\mathfrak{X}^1_Z(M)$, and the inclusion  $\mathfrak{X}^1_Z(M) \to \mathfrak{X}^1(M)$ induces an anchor map $\rho: T_ZM \to TM$ making $T_ZM$ into a Lie algebroid, called the {\em log-tangent bundle} (or {\em $b$-tangent bundle}) associated with $Z\subseteq M$, see \cite[$\S 17.4$]{CW}. 
Note that the anchor map is an isomorphism away from $Z$, but its kernel has rank 1
(locally spanned by $e_1$) over $Z$.



\begin{exer}
\begin{itemize}
\item [(a)] Let $A=T_ZM$ be the log-tangent bundle associated with a hypersurface $Z\subseteq M$, and let $\pi_A$ be an $A$-Poisson structure. Let 
$\pi=\rho(\pi_A)$ be the induced Poisson structure on $M$. Show that $\pi_A$ is nondegenerate if and only if $\pi$ is log-symplectic.

\item [(b)] Let $\pi$ be a log-symplectic structure with degenerate locus $Z\subseteq M$, and let $A=T_ZM$. Show that there is a (unique) $A$-Poisson structure $\pi_A$ such that $\pi=\rho(\pi_A)$ (it may help to use that $Z$ is a Poisson submanifold, see \eqref{eq:pilog}).
\end{itemize}
\end{exer}

In conclusion, log-symplectic structures on a manifold $M$ with degeneracy locus $Z$ are equivalent to symplectic structures on the Lie algebroid $T_ZM$ (see \cite[Prop.~20]{GMP}).

For more general discussions on how to view other types of generically symplectic Poisson structures by means of symplectic structures on appropriate Lie algebroids, see e.g. \cite{klaasse-divisortype,MirSc} and references therein.



\section{Local structure and symplectic foliation} \label{sec:results}



Let $M$ be a manifold and $D\subseteq TM$ a smooth distribution (see $\S$ \ref{subsec:char}). An {\em integral submanifold} of $D$ is a connected submanifold $S \hookrightarrow M$ such that $D|_S=TS$. We say that $D$ is {\em integrable} if any point in $M$ is contained in an integral submanifold. An integrable distribution $D$ gives rise to a partition of $M$ into  maximal integral submanifolds (which are immersed, not necessarily embedded), called {\em leaves}; the collection of leaves will be referred to as a {\em foliation}, see e.g. \cite[Sec.~1.5]{DZ}. When $D$
has constant rank, the classical Frobenius theorem asserts that $D$ is
integrable if and only if it is involutive; the corresponding foliation in this
case is called {\em regular}.

\subsection{Integrability of the characteristic distribution and symplectic leaves}\label{subsec:integ}

Let $\pi$ be a bivector field on $M$, and let $R=\pi^\sharp(T^*M)\subseteq TM$ be its characteristic distribution. Suppose that $S\subseteq  M$ is an integral submanifold of $R$,
$$
R|_S=TS.
$$
Then  
\begin{itemize}
\item $\pi$ is tangent to $S$, and the induced bivector field $\pi_S$ on $S$ is nondegenerate (Exercise~\ref{exer:tangent} (a));
\item $\pi_S$ is equivalent to a nondegenerate 2-form $\omega_S\in \Omega^2(S)$, see $\S$ \ref{subsec:char}.
\end{itemize}

It follows that a bivector field $\pi$ whose characteristic distribution $R$ is integrable gives rise to a foliation in which each leaf $\iota: \mathcal{O}\hookrightarrow M$ carries a nondegenerate 2-form $\omega_\mathcal{O}$ such that
\begin{equation}\label{eq:leaf2form}
\{f,g\}(x)= \omega_\mathcal{O}(X_{\iota^*g}, X_{\iota^*f}) (x), \qquad\ \forall \; x\in \mathcal{O}.
\end{equation}
Note that $\pi$ is uniquely determined by the foliation and leafwise 2-forms.

\begin{exer}\label{exer:sympfol}
Let $\pi$ be a bivector field on $M$ with integrable characteristic distribution.
Check that $\pi$ is Poisson if and only if each leafwise 2-form $\omega_\mathcal{O}$ is closed, hence symplectic (see Exercises~\ref{exer:poisymp} and \ref{exer:tangent} (b)). In particular, $\pi$ is Poisson if leaves have dimension at most $2$.
\end{exer}

As we will see in $\S$ \ref{subsec:split} below, when $\pi$ is Poisson the integrability of $R$ always holds, 
as a consequence of Weinstein's splitting theorem (see also Remark~\ref{rem:Liealgprop}). 
In the regular case, this can be checked more directly using Frobenius theorem.

\begin{exer} \label{exer:preg}
Let $(M,\pi)$ be a regular Poisson manifold. Use \eqref{eq:brkpres} to show that the characteristic distribution $R\subseteq TM$ is involutive, and hence integrable.
\end{exer}

In summary, when $\pi$ is a Poisson structure, its characteristic distribution is integrable (see $\S$ \ref{subsec:split}) and each leaf $\mathcal{O}\hookrightarrow M$  is equipped with a symplectic form $\omega_\mathcal{O}$ uniquely determined by the property that the inclusion map is Poisson; the collection of all such symplectic leaves is called the {\em symplectic foliation} of $\pi$, and $\pi$ is uniquely characterized by its symplectic foliation.

\begin{exer}
Consider $M= \mathbb{T}^3=\mathbb{S}^1\times \mathbb{S}^1\times \mathbb{S}^1=\{(\theta_1,\theta_2,\theta_3)\}$, and for each $\lambda\in \mathbb{R}$, let
$\pi_\lambda= \frac{\partial}{\partial\theta_1}\wedge (\frac{\partial}{\partial \theta_2} + \lambda \frac{\partial}{\partial \theta_3})$. Verify that its characteristic distribution has constant rank equal to $2$ and is involutive, so $\pi_\lambda$ is Poisson (by Exercises~\ref{exer:sympfol} and \ref{exer:preg}). Check that all its symplectic leaves are embedded (and isomorphic to $\mathbb{T}^2$) when $\lambda \in \mathbb{Q}$; if $\lambda$ is irrational, leaves are isomorphic to cylinders and are not embedded.
\end{exer}

\begin{exer}\label{exer:cosympleaves}
Use Exercises~\ref{exer:cosympR} and \ref{exer:charpiN} to check that a submanifold $N$ of a Poisson manifold $M$ is cosymplectic if and only if it intersects each symplectic leaf $\mathcal{O}\hookrightarrow M$ transversally in a symplectic submanifold of $\mathcal{O}$, and, in this case, the symplectic foliation of the induced Poisson structure on $N$ is given by the connected components of the intersections of $N$ with the symplectic leaves of $M$.
\end{exer}

\begin{exer}\label{exer:leafdiffeo}
Let $M_1$ and $M_2$ be Poisson manifolds. Verify that a diffeomorphism $\varphi: M_1\to M_2$ is Poisson if and only if, for every $x\in M_1$, $\varphi$ restricts to a symplectomorphism from the symplectic leaf through $x$ to the symplectic leaf containing $\varphi(x)\in M_2$.
\end{exer}

As indicated by the next exercise, Casimir functions give useful information about symplectic leaves.

\begin{exer}\label{exer:casimirleaves} Let $(M,\pi)$ be a regular $m$-dimensional Poisson manifold with $\mathrm{rank}(\pi)=k$. Suppose that $f_1,\ldots, f_{m-k}$ are Casimir functions such that $df_1, \ldots, df_{m-k}$ are linearly independent at all points. Show that the symplectic leaves are given by connected components of the level sets of $(f_1,\ldots, f_{m-k}): M \to \mathbb{R}^{m-k}$.
\end{exer}

The next exercise illustrates a way to obtain bivector fields with integrable characteristic distributions that are not Poisson.

\begin{exer}\label{exer:Bfield}
 Let $(M,\pi)$ be a Poisson manifold and $B\in \Omega^2(M)$. Let $B^\flat: TM\to T^*M$, $B^\flat(X)=i_XB$, and consider the map $\Id + B^\flat \circ \pi^\sharp: T^*M\to T^*M$.
\begin{itemize}
\item[(a)] Check that $\Id + B^\flat \circ \pi^\sharp$ is an isomorphism if and only if, on each leaf $\mathcal{O}$ of $(M,\pi)$, the 2-form $\omega_\mathcal{O} + B_\mathcal{O}$ is nondegenerate, where $B_\mathcal{O}$ denotes the pullback of $B$ to $\mathcal{O}$.

\item[(b)] If $\Id + B^\flat \circ \pi^\sharp$ is an isomorphism, note that there is a bivector field $\pi_B$ on $M$ such that $\pi_B^\sharp = \pi^\sharp \circ (\Id + B^\flat \circ \pi^\sharp)^{-1}$, and whose characteristic distribution coincides with that of $\pi$, so they have the same leaves. Verify that, for each leaf $\mathcal{O}$, the 2-form on $\mathcal{O}$ induced by $\pi_B$ is given by $\omega_\mathcal{O} + B_\mathcal{O}$. Conclude (see Exercise~\ref{exer:sympfol})  that $\pi_B$ is Poisson if and only if the pullback of $dB$ to each leaf vanishes.
\end{itemize}
\end{exer}

For a closed 2-form $B$, we call the Poisson structure $\pi_B$ of the previous exercise the {\em gauge transformation} of $\pi$ by $B$, and say that the two Poisson structures are {\em gauge equivalent}.

As an application, we can use gauge transformations and the decomposition of a Poisson manifold into symplectic leaves to extend the well-known  Moser's trick from symplectic geometry to the Poisson setting, see e.g. \cite{AMLin}. We briefly recall the setup in symplectic geometry. For a smooth 1-parameter family of symplectic forms $\omega_t \in \Omega^2(M)$, let $B_t=\omega_t-\omega_0$, and suppose that there exists a 1-parameter family of 1-forms $\alpha_t\in \Omega^1(M)$ such that $\frac{d B_t}{dt}= - d\alpha_t$. Let $X_t$ be the time dependent vector field defined by $i_{X_t}\omega_t=\alpha_t$, and let $\phi_t$ be its flow. Then
$$
\phi_t^*\omega_t=\omega_0,
$$
which follows from the equality
$$
\frac{d }{dt} \phi_t^*\omega_t=\phi_t^*\left ( \mathcal{L}_{X_t}\omega_t +\frac{d}{dt}\omega_t  \right ) = \phi_t^*d \left ( i_{X_t}\omega_t - \alpha_t\right ).
$$

\begin{exer}\label{exer:poissonmoser}
Let $\pi_t$ be a 1-parameter family of Poisson structure on $M$ for which there exists a family of closed 2-forms $B_t\in \Omega^2(M)$ with $B_0=0$ such that $\pi_t=(\pi_0)_{B_t}$ (see Exercise~\ref{exer:Bfield}, part (b), for the notation). Suppose that there exists a 1-parameter family of 1-forms $\alpha_t \in \Omega^2(M)$ such that $\frac{d B_t}{dt}= - d\alpha_t$. Consider the time dependent vector field $X_t=\pi_t^\sharp(\alpha_t)$, and let $\phi_t$ be its flow. Then
$$
(\phi_t)_*\pi_0=\pi_t.
$$
\end{exer}


\subsection{The splitting theorem} \label{subsec:split}

The central result describing the local structure of a Poisson manifold is
{\em Weinstein's splitting theorem} \cite{We83}: {\em for a Poisson manifold $(M,\pi)$, any point $m$ admits a neighborhood with coordinates $(q_1,\ldots,q_k,p_1,\ldots,p_k,y_1,\ldots,y_l)$ (centered at $m$) where $\pi$ takes the form
\begin{equation}\label{eq:splitting}
\pi = {\sum_{i=1}^k \frac{\partial}{\partial p_i}\wedge \frac{\partial}{\partial q_i}} + \frac{1}{2}\sum_{i,j= 1}^l {\varphi_{ij}(y)\frac{\partial}{\partial y_i}\wedge \frac{\partial}{\partial y_j}}, \qquad \mbox{and} \quad \; \varphi_{ij}(0)=0.
\end{equation}
}
Since $\varphi_{ij}$ vanishes at $m$, the rank of $\pi$ at $m$ is $2k$.

The following are consequences of this result.

\begin{itemize}
\item Since the functions $\varphi_{ij}$ only depend on the coordinates $y_1$, $\ldots$, $y_l$, the theorem gives a local decomposition  of $M$ as the direct product of two Poisson manifolds, 
$$
(S,\pi_S)\times (N,\pi_N), 
$$
where $S$, defined by the coordinates $(q_i,p_i)$, has a nondegenerate Poisson structure $\pi_S= \sum_i \frac{\partial}{\partial p_i}\wedge \frac{\partial}{\partial q_i}$, and $N$, defined by the coordinates $(y_1,\ldots,y_l)$, has Poisson structure $\pi_N= \sum_{i<j} {\varphi_{ij}(y)\frac{\partial}{\partial y_i}\wedge \frac{\partial}{\partial y_j}}$ whose rank vanishes at the point $y=0$.

\item When $\pi$ has constant rank $2k$ around $m$, it locally looks like Poisson's original bracket \eqref{eq:original},
$$
\pi= \sum_{i=1}^k \frac{\partial}{\partial p_i}\wedge \frac{\partial}{\partial q_i}.
$$
When $\pi$ is nondegenerate, this local normal form recovers the Darboux theorem for symplectic structures.

\item The submanifold $S \hookrightarrow M$, $(q_i,p_i)\mapsto ((q_i,p_i), 0)$ satisfies $TS=R|_S$, i.e., it is an integral submanifold of the characteristic distribution $R$ containing $m=(0,0)$. It follows that $R$ is integrable, 
and its leaves define the symplectic foliation of $\pi$, as explained in $\S \ref{subsec:integ}$. In Weinstein's splitting coordinates, the leafwise symplectic forms have the canonical expression $\sum_i dq_i\wedge dp_i$.

\item Leaves of $(M,\pi)$ can be characterized by the following equivalence relation: $x, x' \in M$ belong to the same leaf if and only if there are functions $f_1, \ldots, f_r$ on $M$ such that
\begin{equation}\label{eq:flows}
x' = \phi^{t_1}_{X_{f_1}}\circ \ldots \circ \phi^{t_r}_{X_{f_r}}(x), 
\end{equation}
where $\phi^t_X$ denotes the flow of a vector field $X$ at time $t$. (This can be directly verified if $x$, $x'$ belong to a neighborhood where the splitting theorem holds, and in general by covering any path from $x$ to $x'$ with finitely many of those.)

\end{itemize}

\begin{exer}\label{exer:complete}
Let $\varphi: M_1\to M_2$ be a Poisson map that is {\em complete}, in the sense that whenever a Hamiltonian vector field $X_f$ on $M_2$ is complete, $X_{\varphi^*f}$ is also complete. 
Show that if $y=\varphi(x)$, then the entire symplectic leaf through $y$ is contained in the image of $\varphi$. Hence the image of a complete Poisson map is a union of symplectic leaves. 
\end{exer}

The following steps outline the proof of the splitting theorem.

\begin{enumerate}

\item Assuming that the rank of $\pi$ at $m$ is not zero, let $p_1$ be a function such that $X_{p_1}|_m\neq 0$. By the straightening theorem one can find a function $q_1$ around $m$ such that $\{p_1,q_1\}=1$.

\item The vector fields $X_{p_1}$ and $X_{q_1}$ are linearly independent around $m$ and commute: $[X_{p_1},X_{q_1}]= X_{\{p_1,q_1\}}= X_1=0$. It follows that there are coordinates $(y_1,\ldots,y_n)$ around $m$ satisfying $X_{q_1}=\partial_{y_{n-1}}$, $X_{p_1}= \partial_{y_n}$, in such a way that $dy_j(X_{q_1})=0$ and $dy_j(X_{p_1})=0$ for $j=1,\ldots,n-2$.

\item The differentials $dq_1$, $dp_1$, $dy_1, \ldots dy_{n-2}$ are linearly independent, so the functions $q_1$, $p_1$, $y_1,\ldots y_{n-2}$ form a local coordinate system around $m$, satisfying 
$$
\{p_1,q_1\}=1, \quad \{p_1,y_j\}=\{q_1, y_j\}=0.
$$ 
Hence $X_{p_1}=\partial_{q_1}$ and $X_{q_1}=-\partial_{p_1}$.

\item By the Jacobi identity,
$$
\partial_{p_1}\{y_i,y_j\}=-\{q_1, \{y_i,y_j\}\}=0, \quad\;\; \partial_{q_1}\{y_i,y_j\}=\{p_1, \{y_i,y_j\}\}=0,
$$
so the functions $\varphi_{ij}=\{y_i,y_j\}$ only depend on $y_1, \ldots,y_{n-2}$.
Hence
$$
\pi = \frac{\partial}{\partial p_1}\wedge \frac{\partial}{\partial q_1} + \frac{1}{2}\sum_{i,j}\varphi_{ij}(y)\frac{\partial}{\partial y_i}\wedge \frac{\partial}{\partial y_j}.
$$
If $\varphi_{ij}(0)=0$, the theorem is proven. If not, one repeats the argument for the Poisson structure $\frac{1}{2}\sum_{i,j}\varphi_{ij}(y)\frac{\partial}{\partial y_i}\wedge \frac{\partial}{\partial y_j}$.
\end{enumerate}

Different proofs of the splitting theorem, leading to more general results, can be found in \cite{BLM,FrejMarcNormal} 


\subsection{Symplectic leaves in examples}

\subsubsection{Poisson surfaces}For a Poisson structure $\pi$ on a 2-dimensional manifold $M$ (Example~\ref{exer:2d}), zero-dimensional leaves are given by points where the rank of $\pi$ is zero, and 2-dimensional leaves are given by the connected components of the open subset where $\pi$ is nondegenerate

\subsubsection{Duals of Lie algebras}
For a Lie algebra $\mathfrak{g}$, the infinitesimal coadjoint action of $\mathfrak{g}$ on $\mathfrak{g}^*$ sends $u\in \mathfrak{g}$ to the vector field $\mathrm{ad}^*_u$ on $\mathfrak{g}^*$ given by
$$
\mathfrak{g}^*\ni \xi \mapsto \mathrm{ad}_u^*(\xi)=\xi([u,\cdot]) \in \mathfrak{g}^*=T_\xi \mathfrak{g}^*.
$$
It follows from \eqref{eq:Xu} that the characteristic distribution of the linear Poisson structure on $\mathfrak{g}^*$ is tangent to the orbits of the coadjoint $G$-action, where $G$ is any connected Lie group with Lie algebra $\mathfrak{g}$.  Following \eqref{eq:leaf2form},
the symplectic form on a coadjoint orbit $\mathcal{O}$ is given, at $\xi \in \mathcal{O}$, by
$$
\omega_{\mathcal{O}}( \mathrm{ad}_u^*(\xi),  \mathrm{ad}_v^*(\xi))=\xi([v,u]).
$$


In the next exercise, we will consider $\mathbb{R}^3=\{(x,y,z)\}$ with linear Poisson structures corresponding to each of the 3-dimensional Lie algebras $\mathfrak{su}(2)$, $\mathfrak{sl}(2,\mathbb{R})$ and $\mathfrak{sb}(2,\mathbb{C})$ (see Exercise~\ref{exer:susb}), respectively.

\begin{exer}\label{exer:g*leaves} Verify the following (Exercise~\ref{exer:casimirleaves} should be useful).
\begin{itemize}
    \item[(a)] The only zero-dimensional leaf of $\pi= z\partial_x\wedge\partial_y + x \partial_y\wedge\partial_z + y \partial_z\wedge\partial_x$ is the origin, and its 2-dimensional leaves are given by the fibers of the map $f: \mathbb{R}^3\setminus \{(0,0,0)\} \to \mathbb{R}$, $f(x,y,z)=x^2+y^2+z^2$.
     
    \item[(b)] The only zero-dimensional leaf of $\pi= -z\partial_x\wedge \partial_y + x \partial_y\wedge \partial_z + y \partial_z\wedge \partial_x$ is the origin, and its 2-dimensional leaves are the (connected components of the) fibers of the map $f: \mathbb{R}^3\setminus \{(0,0,0)\} \to \mathbb{R}$, $f(x,y,z)=x^2+y^2-z^2$.

    
    \item[(c)] The zero-dimensional leaves of $\pi=x\partial_x\wedge \partial_z + y\partial_y\wedge \partial_z$ are points along the $z$-axis, and its 2-dimensional leaves are the fibers of the map 
    $$
    f: \mathbb{R}^3 \setminus \{(0,0,z), z\in \mathbb{R}\}\to \mathbb{S}^1, \;\;\;f(x,y,z)= 
    \left(\frac{x}{\sqrt{x^2+y^2}}, \frac{y}{\sqrt{x^2+y^2}}\right).
    $$
    
\end{itemize}
\end{exer}


In the more general case of Poisson-Lie groups, symplectic leaves 
are given by orbits of the so-called {\em dressing action}, see e.g. \cite[$\S$ 3.4]{KoroSoib} and \cite[Sec.~2]{LuWe}, that generalize the coadjoint action on the dual of a Lie algebra.

The next exercise describes some leaves of duals of Lie algebroids.

\begin{exer}\label{exer:leafA*1}
Let $A\to M$ be a Lie algebroid with surjective anchor map $\rho:A \twoheadrightarrow TM$. Show 
that, if $M$ is connected, then  $\rho^*$ identifies $(T^*M, \omega_{can})$ with  a symplectic leaf of $A^*$ (cf. Exercises \ref{exer:dualpoismap} and \ref{exer:complete}). More generally, a Lie algebroid $A$ defines a  foliation on $M$ whose tangent distribution is $\rho(A)\subseteq TM$ (see  Remark~\ref{rem:Liealgprop}); show that, for each leaf $\mathcal{O}\hookrightarrow M$, the cotangent bundle $(T^*\mathcal{O},\omega_{can})$ is identified with a symplectic leaf of $A^*$ via $\rho^*|_{\mathcal{O}}: T^*\mathcal{O}\to A^*|_\mathcal{O}\subseteq A^*$.
\end{exer}

\subsubsection{Quotients of symplectic manifolds}\label{subsub:quot}

We will now discuss the symplectic foliation of  Poisson manifolds obtained as quotients of symplectic manifolds. We start with some general considerations.

\begin{exer}\label{exer:lib1}
Let $(S,\omega)$ be a symplectic manifold and $\varphi: S\to M$ a surjective submersion. Consider the ``vertical'' bundle $\mathcal{V}=\Ker(d\varphi)\subseteq TS$, and let $\mathcal{V}^\omega \subseteq TS$ be the symplectic orthogonal subbundle.  Let $C^\infty_\mathcal{V}(S)=\{g\in C^\infty(S)\,|\, dg\in \Gamma(\mathrm{Ann}(\mathcal{V})) \}$, i.e., the space of functions that are locally constant on $\varphi$-fibers.
Notice that
\begin{itemize}
\item  $g\in C^\infty_\mathcal{V}(S)$ if and only if $X_g\in \Gamma(\mathcal{V}^\omega)$, 
\item  $(d_x\varphi)^*: T_{\varphi(x)}^*M \stackrel{\sim}{\to} \mathrm{Ann}(\mathcal{V}_x)$, and 
\begin{equation}\label{eq:Vom}
\mathcal{V}^\omega_x=\{X_{\varphi^*f}|_x \; , f\in C^\infty(M)\},  \;\;\; \forall x\in S.
\end{equation}
\end{itemize}
Verify that
\begin{itemize}
\item[(a)] If $\varphi^*C^\infty(M)\subseteq C^\infty(S)$ is a Poisson subalgebra (equivalently, there is a Poisson structure on $M$ for which $\varphi$ is a Poisson map), then $\mathcal{V}^\omega$ is involutive.
\item[(b)] If $\mathcal{V}^\omega$ is involutive, then $C^\infty_\mathcal{V}(S) \subseteq C^\infty(S)$ is a Poisson subalgebra.
\end{itemize}
\end{exer}

Let $(M,\pi)$ be a Poisson manifold, and suppose that there is a symplectic manifold $(S,\omega)$ and a Poisson surjective submersion $
\varphi: S\to M$. We will use $\varphi$ to give a description of the symplectic leaves of $M$.

Consider the vector subbundle $\mathcal{V}^\omega \subseteq TS$, the symplectic orthogonal of the distribution $\mathcal{V}$ tangent to $\varphi$-fibers. By the previous exercise (part (a)), $\mathcal{V}^\omega$ is involutive, so it is the tangent distribution to a  regular foliation on $S$. Since the Hamiltonian vector fields $X_{\varphi^*f}$ and $X_f$ are $\varphi$-related for all $f\in C^\infty(M)$ (see Exercise~\ref{exer:pmaps}), \eqref{eq:Vom} implies that that
\begin{equation*}\label{eq:distphi}
d\varphi( \mathcal{V}^\omega_x)= R_{\varphi(x)}, \quad \forall x\in S,
\end{equation*}
where $R$ is the characteristic distribution of $\pi$. As a consequence, denoting by $\mathcal{L}\stackrel{\iota_\mathcal{L}}{\hookrightarrow} S$  the leaf through a given point $x\in S$ and $\mathcal{O}$ the symplectic leaf of $\pi$ containing $\varphi(x)$, then
$\varphi(\mathcal{L}) \subseteq \mathcal{O}$; moreover, the map $\varphi_\mathcal{L}:=\varphi|_{\mathcal{L}}: \mathcal{L} \to \mathcal{O}$ is a submersion. 

By Exercise~\ref{exer:complete}, if we additionally assume the Poisson map $\varphi$ to be {\em complete}, then
$$
\varphi(\mathcal{L}) = \mathcal{O},
$$
showing that the symplectic leaves in $M$ are the images of the leaves of $\mathcal{V}^\omega\subseteq TS$ under $\varphi$. 
The discussion is summarized by the following commutative diagram of horizontal inclusions and vertical surjective submersions,
\begin{equation}\label{eq:diag}
\xymatrix {
{\mathcal{L}} \ar[d]_{\varphi_{\mathcal{L}}} \ar[r]^{\iota_{\mathcal{L}}} &
{(S,\omega)} \ar[d]^\varphi\\
(\mathcal{O},\omega_\mathcal{O}) \ar[r] & (M, \pi).
}
\end{equation}
The next exercise shows that the symplectic forms $\omega$ and $\omega_\mathcal{O}$ are related by
\begin{equation}\label{eq:redforms}
\varphi_{\mathcal{L}}^* \omega_{\mathcal{O}} = \iota_{\mathcal{L}}^*\omega.
\end{equation}

\begin{exer}\label{exer:redomega}
For $f,g \in C^\infty(M)$, check that 
$$
\varphi_{\mathcal{L}}^*\omega_{\mathcal{O}}(X_{\varphi^*g}(x),X_{\varphi^*f}(x))=\{f,g\}(\varphi(x)),
\quad\;\;\; \iota_{\mathcal{L}}^*\omega(X_{\varphi^*g}(x),X_{\varphi^*f}(x))=\{\varphi^*f, \varphi^*g\}(x).
$$
\end{exer}

A special case of the previous setup was considered in $\S$ \ref{subsec:symmetry}: $(S,\omega)$ is a symplectic manifold carrying an action of a Lie group $G$ by symplectomorphisms that is
free and proper, and $M=S/G$ is equipped with the Poisson structure such that the quotient map $\varphi: S\to S/G$ is  Poisson. We keep the notation $\mathcal{V}\subseteq TS$ for the distribution tangent to the $G$-orbits.

\begin{exer}\label{exer:completequot}
Show that the Poisson map $\varphi: S\to M=S/G$ is complete.
\end{exer}

Suppose that the $G$-action on $S$ is Hamiltonian (cf. $\S$ \ref{subsec:linear1}) with moment map $\mu: S\to \mathfrak{g}^*$; recall that this means that $\mu$ is $G$-equivariant (with respect to the coadjoint action on $\mathfrak{g}^*$) and satisfies 
$$
i_{\psi(u)}\omega=d(\mu^*u), \qquad \forall u\in \mathfrak{g},
$$ 
where $\psi: \mathfrak{g}\to \mathfrak{X}^1(M)$ is the infinitesimal action and $\mu^*u(x)=\langle \mu(x), u\rangle$.

For $\xi \in \mathfrak{g}^*$, denote by $G_\xi \subseteq G$ the stabilizer group of $\xi$.

\begin{exer}\label{exer:reduction}
Verify that
\begin{itemize}
\item[(a)] $\mu:S \to \mathfrak{g}^*$ is a submersion,
\item[(b)] $\Ker(d\mu) = \mathcal{V}^\omega$,
\item[(c)] for $\xi \in \mathfrak{g}^*$ and $x\in \mu^{-1}(\xi)$, $\sigma x \in \mu^{-1}(\xi)$ if and only if $\sigma \in G_\xi$.
\end{itemize}
\end{exer}

It follows from  (a) and (b) of the previous exercise that, when the $G$-action on $S$ is Hamiltonian, the  leaves of $\mathcal{V}^\omega$ are the connected components of the level sets of $\mu$.
From part (c) of the exercise, along with \eqref{eq:diag} and \eqref{eq:redforms}, it follows that the symplectic leaves of $S/G$ are identified with the connected components of the symplectic reduced spaces $\mu^{-1}(\xi)/G_\xi$,
$$
\xymatrix {
{\mathcal{L}} \ar[d]_{} \ar[r]^{} & \mu^{-1}(\xi) \ar[d] \ar[r] &
{(S,\omega)} \ar[d]^\varphi\\
\mathcal{O} \ar[r] &  \mu^{-1}(\xi)/G_\xi \ar[r]   & (S/G, \pi).
}
$$

\subsubsection{Standard Poisson structure on $\mathbb{S}^3$}

Let us consider the standard Poisson structure on $\mathbb{S}^3 = \{(x,y,z,w) \in \mathbb{R}^4\, |\, x^2+y^2+z^2+w^2=1\}$, given by
\begin{align*}
\pi&= (z^2+w^2){\partial_x}\wedge {\partial_y} + \left ( -y {\partial_x} + x {\partial_y}  \right )\wedge \left ( z {\partial_z} + w {\partial_w} \right )\\
&= (-z\partial_x+x\partial_z)\wedge (-z\partial_y+y\partial_z) + (-w\partial_x+x\partial_w)\wedge (-w\partial_y+y\partial_w),
\end{align*}
see Exercise~\ref{exer:SU2} (note that the second expression for $\pi$ makes apparent the fact that it is tangent to $\mathbb{S}^3$). The zero-dimensional leaves of $(\mathbb{S}^3,\pi)$ are the points in the unit circle $C=\{(x,y,0,0)\;|\; x^2+y^2=1 \}$. Hence the restriction of $\pi$ to the open subset $\mathbb{S}^3\setminus C$ is a regular Poisson structure of rank 2. 

Consider $\mathbb{S}^1 \subseteq \mathbb{R}^2$ and the map 
$$
F: \mathbb{S}^3\setminus C {\, \to \,} \mathbb{S}^1, \; \quad \; F(x,y,z,w) = \left(\frac{z}{\sqrt{z^2+w^2}}, \frac{w}{\sqrt{z^2+w^2}} \right).
$$
Denoting by $\mathbb{D}^2\subseteq \mathbb{R}^2$ the open unit disk, we have a  diffeomorphism 
\begin{equation}\label{eq:diffS3}
\mathbb{S}^3\setminus C \stackrel{\sim}{\to} \mathbb{D}^2\times \mathbb{S}^1,\;\;\; (x,y,z,w) \mapsto ((x,y), F(x,y,z,w)),
\end{equation}
with inverse map given by
$((x,y), c)\mapsto (x,y,\sqrt{1-x^2-y^2} \, c)$.

\begin{exer} 
\begin{itemize}
\item[(a)] Check that the functions $F_1(x,y,z,w)= \frac{z}{\sqrt{z^2+w^2}}$ and $F_2(x,y,z,w)=\frac{w}{\sqrt{z^2+w^2}}$ are Casimirs for $\pi$ on $\mathbb{S}^3\setminus C$.
\item[(b)]  Show that the diffeomorphism \eqref{eq:diffS3} identifies $(\mathbb{S}^3\setminus C,\pi)$ with the direct product of Poisson manifolds  $(\mathbb{D}^2,\pi_{\mathbb{D}^2}) \times (\mathbb{S}^1,0)$, where $\pi_{\mathbb{D}^2}=(1-x^2-y^2)\partial_x\wedge \partial_y$.
\end{itemize}
\end{exer}

It follows from the previous exercise that the 2-dimensional leaves of $\pi$ agree with the fibers of the map $F$. In more detail, the set of 2-dimensional leaves of $(\mathbb{S}^3,\pi)$ is in bijection with $\mathbb{S}^1$: for each $c\in \mathbb{S}^1$ the corresponding leaf is the image of 
$$
\mathbb{D}^2\to \mathbb{S}^3, \quad \;\;\; (x,y)\mapsto (x,y,\sqrt{1-x^2-y^2} \, c),
$$
with symplectic form in these coordinates given by $\frac{1}{x^2+y^2-1}dx\wedge dy$. 

\begin{exer}
Give an analogous description of the symplectic leaves of the higher-dimensional Poisson spheres of Exercise~\ref{exer:spheres}.
\end{exer}


Recall from $\S$ \ref{subsec:PL} that $\mathbb{S}^2$ has a natural Poisson homogeneous structure such that the Hopf map $\mathbb{S}^3\to \mathbb{S}^2$ is a Poisson map (see \eqref{eq:bruhS2}). 
This Poisson structure has two symplectic leaves: one singular point (the image of $C\subseteq \mathbb{S}^3$), and an open symplectic leaf given by the complement of the singular point.
The Hopf map takes each 2-dimensional symplectic leaf in $\mathbb{S}^3\setminus C$ symplectomorphically onto the open symplectic leaf of $\mathbb{S}^2$.

\begin{exer}
Describe the symplectic leaves of the Poisson structures on $\mathbb{S}^2$ of Exercise~\ref{exer:SU2cov}.
\end{exer}

\subsubsection{Log-symplectic structures}
The following exercise describes the symplectic leaves of log-symplectic manifolds ($\S$ \ref{subsec:log}).

\begin{exer} \label{exer:log}
Let $(M^{2n},\pi)$ be a log-symplectic manifold with degeneracy locus $Z$, and $m\in Z$. Check that, in Weinstein's splitting coordinates $(q_1,\ldots,q_k, p_1,\ldots,p_k, y_1,\ldots, y_l)$ around $m$, one must have $l=2$. Use Exercise~\ref{exer:log2dim} to conclude that there are local coordinates around $m$ such that
$$
\pi = y_1 \frac{\partial}{\partial y_1}\wedge \frac{\partial}{\partial y_2} + \sum_{i=1}^{n-1}\frac{\partial}{\partial p_i}\wedge \frac{\partial}{\partial q_i}.
$$
As a Poisson submanifold, $Z$ inherits a regular Poisson structure of corank one. So the symplectic leaves of $(M,\pi)$ are either open subsets given by connected components of $M\setminus Z$, or are contained in $Z$ and have dimension $2n-2$.
\end{exer}


\subsection{Transverse structure and linearization}\label{subsec:transverse}

Given a Poisson manifold $(M,\pi)$, recall that Weinstein's local splitting theorem ($\S$ \ref{subsec:split}) asserts that any point in $M$ admits a neighborhood that can be written as a direct product of Poisson manifolds, 
$$
(S,\pi_S) \times (N,\pi_N),
$$ 
where $\pi_S$ is nondegenerate and $\pi_N$ vanishes at the given point. The symplectic factor $S$ is an open subset of a symplectic leaf $\mathcal{O}\subseteq M$. The degenerate factor $N$ sits in $M$ as a cosymplectic submanifold (see Exercise~\ref{exer:cosympfactor}) and $\pi_N$
is called the {\em transverse Poisson structure} to the symplectic leaf $\mathcal{O}$;
although $(N,\pi_N)$ is not intrinsically defined, the terminology is justified by the following uniqueness result \cite{We83}: {\em let $m_0, m_1 \in \mathcal{O}$ and, for each $i=0,1$, let $N_i$ be a cosymplectic submanifold containing $m_i$ such that 
\begin{equation}\label{eq:slice}
T_{m_i}\mathcal{O}\oplus T_{m_i}N_i=T_{m_i}M;
\end{equation}
then there are open neighborhoods $U_0$ of $m_0$ in $N_0$ and $U_1$ of $m_1$ in $N_1$ with a Poisson diffeomorphism $(U_0,\pi_{N_0})\cong (U_1,\pi_{N_1})$ taking $m_0$ to $m_1$.}

In the previous result, one can think of the transverse cosymplectic submanifolds as obtained from local Weinstein splittings around points in the symplectic leaf. More generally,  any submanifold $N$ containing $m\in \mathcal{O}$ and satisfying the transversality condition \eqref{eq:slice} will be cosymplectic in a neighborhood of $m$, see Exercise~\ref{exer:slice}.

We outline the proof of the uniqueness of the transverse structure.

\begin{itemize}
\item Since $m_0$ and $m_1$ belong to the same symplectic leaf, they are connected by a composition of Hamiltonian flows as in \eqref{eq:flows}; this composition defines a Poisson automorphism of $M$ taking $m_0$ to $m_1$, and as a consequence we can reduce the problem to the case where $m_0=m_1$.  So we henceforth assume that we have a Weinstein splitting chart $S\times N$ around a point $m=(0,0)$ and $N' \subseteq S\times N$ is a cosymplectic submanifold, with induced Poisson structure denoted by $\pi_{N'}$, containing $m$ and satisfying 
$$
T_mS\oplus T_m N'=T_mM.
$$

\item By the previous transversality condition at $m$, the map $N'\to N$ given by restriction of the natural projection $S\times N \to N$ is a local diffeomorphism  around $m$. By shrinking $N$, we can assume that $N'$ projects diffeomorphically onto $N$, and hence there is a smooth map $\phi: N \to S$ such that
$$
N'= \{(\phi(y),y)\, | \, y\in N\}.
$$
Let $\pi'$ be the Poisson structure on $N$ corresponding to $\pi_{N'}$ via the projection $N'\stackrel{\sim}{\to} N$, or, equivalently, such that the inverse diffeomorphism 
$$
\iota_\phi: N \stackrel{\sim}{\to} N', \quad \;\; \iota_\phi(y)=(\phi(y),y),
$$ 
is a Poisson map. 
The remainder of the proof consists in showing that $\pi_N$ and $\pi'$ are locally isomorphic around $m$. We start by checking that they are gauge equivalent (see Exercise~\ref{exer:Bfield}).

\item The symplectic leaves of $(N',\pi_{N'})$ are obtained by intersecting $N'$ with the leaves of $(S,\pi_S)\times (N,\pi_N)$ (Exercise~\ref{exer:cosympleaves}). So each leaf is 
given by
$$
(S\times \mathcal{O})\cap N'=\{(\phi(y),y)\,|\, y\in \mathcal{O}\},
$$
with symplectic form given by the pullback of $\omega_S \times \omega_{\mathcal{O}}$, where  $(\mathcal{O},\omega_\mathcal{O})$ is a symplectic leaf of $(N,\pi_N)$ and $\omega_S=\pi_S^{-1}$. 
We know that the leaves of $\pi'$ correspond to the leaves of $\pi_{N'}$ under the Poisson diffeomorphism $\iota_\phi: (N,\pi')\to (N',\pi_{N'})$, so they 
are given by 
$$
\iota_\phi^{-1}((S\times \mathcal{O})\cap N')=\mathcal{O},
$$
with symplectic form given by the pullback of $\omega_S\times \omega_\mathcal{O}$ by $\iota_\phi: \mathcal{O}\to S\times \mathcal{O}$,
$$
\iota^*\phi^*\omega_S + \omega_\mathcal{O},
$$
where $\iota: \mathcal{O} \hookrightarrow N$ is the inclusion.
Following Exercise~\ref{exer:Bfield}, we conclude that $\pi'$ is the gauge transformation of  $\pi_N$ by the closed 2-form $B=\phi^*\omega_S$,
$$
(\pi')^\sharp = \pi_N^\sharp\circ (\Id + B^\flat\circ \pi_N^\sharp)^{-1}.
$$

\item In a neighborhood admitting Weinstein's splitting coordinates, the symplectic form $\omega_S$ is exact, and the local isomorphism between $\pi_N$ and $\pi'$ around $m$ can be obtained by Moser's trick (Exercise~\ref{exer:poissonmoser}), as explained in the next exercise.

\end{itemize}

\begin{exer}
Let $\pi$ be a Poisson structure on $N$ vanishing at $y_0\in N$.
Let $\pi'$ be a Poisson structure gauge equivalent to $\pi$ via an exact 2-form $B=-d\alpha$. The following items show that there are open neighborhoods $U$ and $U'$ of $y_0$ and a Poisson diffeomorphism $(U,\pi)\stackrel{\sim}{\to} (U',\pi')$.
\begin{itemize}
\item[(a)] Let $B_t=tB$. Show that there is an open neighborhood of $y_0$ where $(\Id + B_t^\flat\circ \pi^\sharp)$ is invertible for all $t\in [0,1]$, so there we can consider the family of Poisson structures $\pi_t:= \pi_{B_t}$ (note that they all vanish at $y_0$).
\item[(b)] Let $X_t=\pi_t^\sharp(\alpha)$. Show that there is a neighborhood of $y_0$ where the flow $\phi_t$ of $X_t$ is defined for all $t\in [0,1]$. By Exercise~\ref{exer:poissonmoser},
$$
(\phi_t)_*\pi_t=\pi_0=\pi,
$$
so $\phi_1$ gives the desired Poisson diffeomorphism.
\end{itemize}
\end{exer}

At any point $m$ of a Poisson manifold $(M,\pi)$ there is a natural Lie algebra structure on the conormal space to the symplectic leaf $\mathcal{O}$ through $m$.
Denoting by $\nu_m=T_mM/T_m\mathcal{O}$ the normal space  to $\mathcal{O}$ at $m$, the Lie bracket on the conormal space
$$
\nu_m^*= (T_mM/T_m\mathcal{O})^* = \mathrm{Ann}(T_m\mathcal{O}) = \Ker(\pi^\sharp_m)
$$
is defined as follows: for $\alpha, \beta \in \mathrm{Ann}(T_m\mathcal{O})$,  take functions $f,g \in C^\infty(M)$ such that
$\alpha=df|_m$ and $\beta=dg|_m$ and set
$$
[\alpha,\beta]:= d (\{f,g\})|_m.
$$
Note that $d (\{f,g\})|_m \in \mathrm{Ann}(T_m\mathcal{O})$ as a consequence of the Jacobi identity:
$$
d (\{f,g\})|_m(X_h)= dg|_m(X_{\{h,f\}}) - df|_m(X_{\{h,g\}})=0
$$
for any $h\in C^\infty(M)$.

\begin{exer} Check that, by the Leibniz identity, the value of $[\alpha,\beta]$ is well defined, independent of the choices of $f$ and $g$. 
\end{exer}

The Lie algebra $(\nu_m^*, [\cdot,\cdot])$ is called the {\em transverse Lie algebra}, or the {\em isotropy Lie algebra} of $\pi$ at $m$. 
The corresponding linear Poisson structure on the normal space $\nu_m$ is closely related to the transverse Poisson structure. To see that, we consider
 Weinstein's splitting coordinates centered at $m$, with respect to which the transverse Poisson structure on $N=\{(y_j)\}$ is written as 
$$
\pi_N= \sum_{i<j}\varphi_{ij}(y)\frac{\partial}{\partial y_i}\wedge \frac{\partial}{\partial y_j},
$$ 
with $\varphi_{ij}(0)=0$. Then the Lie bracket on $\nu_m^* = T_0^*N$ satisfies
$$
[dy_i,dy_j]=d(\{y_i,y_j\})|_0 = \sum_k \frac{\partial \varphi_{ij}}{\partial y_k}(0)dy_k,
$$
i.e., the structural constants in this basis are $c_{ijk}= \frac{\partial \varphi_{ij}}{\partial y_k}(0)$. The corresponding linear Poisson structure on $\nu_m=T_0N$ is given by
$$
\sum_{i<j} \left ( \sum_k  \frac{\partial \varphi_{ij}}{\partial y_k}(0) y_k \right) \frac{\partial}{\partial y_i}\wedge \frac{\partial}{\partial y_j},
$$
which is the first-order term of the Taylor expansion of $\pi_N$ around $m$. So the linear Poisson structure dual to the transverse Lie algebra at a given point is the {\em linear approximation} to the transverse Poisson structure at that point.

A consequence of 
 Weinstein's splitting theorem is that, in the local study of Poisson manifolds,  it suffices to consider Poisson structures in the neighborhoods of points where they vanish. Given a Poisson manifold $(M,\pi)$ and $m\in M$ such that $\pi_m=0$, a natural question is whether 
$\pi$ is locally isomorphic around $m$ to its linear approximation on $\nu_m=T_mM$. In the affirmative case, $\pi$ is said to be {\em linearizable} at $m$. A simple example of a Poisson structure that is not linearizable is $(x^2+y^2)\partial_x\wedge \partial_y$ at the origin: its linear approximation is the zero Poisson structure, while $\pi$ is nondegenerate away from the origin.

\begin{exer}
 Show that, in dimension 2, any Poisson structure whose isotropy Lie algebra at a vanishing point is nonabelian is linearizable at that point (use Exercise~\ref{exer:log2dim}).
\end{exer}

The systematic study of  linearization of Poisson structures began with Weinstein in \cite{We83} and led to an active line of research in Poisson geometry; see \cite[Ch.~4]{DZ} for many of the important results and references. A guiding question has been that of characterizing Lie algebras $\mathfrak{g}$ with the property that any Poisson structure having $\mathfrak{g}$ as its isotropy Lie algebra at a vanishing point is linearizable at that point. Lie algebras with this property are called {\em nondegenerate}, and {\em degenerate} otherwise. One of the main theorems in this direction, obtained by Conn in \cite{Conn}, is that any semisimple Lie algebra of compact type\footnote{I.e., a Lie algebra with negative definite Killing form, or equivalently such that all of its integrating Lie groups are compact.} is nondegenerate (as conjectured in \cite[$\S$ 6]{We83}). Conn's proof of this linearization result was analytic and relied on Nash-Moser type techniques; a soft geometric proof was later obtained in \cite{CrFeLin}.
There is also a ``semi-local'' version of the linearization problem in the neighborhood of a general symplectic leaf (rather than just vanishing points of the Poisson structure) as well as an extension of Conn's theorem to this setting \cite{CrMaLin} (see \cite{FeMaLin} for further generalizations and references).

The next exercise shows that $\mathfrak{sl}(2,\mathbb{R})$ is degenerate by giving an explicit example of a nonlinearizable Poisson structure with $\mathfrak{sl}(2,\mathbb{R})$ as its isotropy Lie algebra \cite[Prop.~6.3]{We83}.


\begin{exer}
Consider the linear Poisson structure on $\mathbb{R}^3$ dual to $\mathfrak{sl}(2,\mathbb{R})$,
$$
\pi_1 = -z\partial_x\wedge \partial_y + x \partial_y\wedge \partial_z + y \partial_z\wedge \partial_x = X\wedge Y,
$$
where $X=x \partial_y- y \partial_x$ and $Y=\partial_z + \frac{z}{x^2+y^2}(x\partial_x+y\partial_y)$.
Let $f\in C^\infty(\mathbb{R})$ be such that $f(t)=0$ for $t\leq 0$, and $f(t)>0$ for $t>0$, consider the vector field 
$$
Z=\frac{f(x^2+y^2-z^2)}{{x^2+y^2}}(x\partial_x + y\partial_y),
$$
and let $\pi= (X-Z)\wedge Y$.
\begin{itemize}
\item[(a)] Show that $[X,Y]=0$ and $[Z,Y]=0$, and conclude that $[\pi,\pi]= 0$ (see Exercise~\ref{exer:ld}).
\item[(b)] Verify that the linear approximation of $\pi$ at the origin is $\pi_1$.
\item[(c)] Notice that $X-Z$  is the Hamiltonian vector field of $h(x,y,z)=-z$ with respect to $\pi$.
Check that its integral curves in the region $\{ x^2+y^2-z^2>0\}$ spiral in towards the cone $\{ x^2+y^2-z^2 = 0\}$. By comparison with the symplectic leaves of $\pi_1$  (see Exercise~\ref{exer:g*leaves}, part (b)), conclude that $\pi$ is not linearizable at the origin. 
\end{itemize}

\end{exer}

The linearization problem has been extensively studied in the realm of Poisson Lie groups.
Most compact Lie groups with the standard Poisson structure turn out to be non-linearizable at the group unit; in the case of a compact simple Lie group, the only exception is $SU(2)$ \cite{CGR}. By contrast, the dual Poisson Lie group $K^*$ to any compact, semisimple Lie group $K$ is globally linearizable, i.e., it is globally Poisson diffeomorphic to its linear approximation on $T_eK^*=\mathfrak{k}^*$ \cite{GinzWein}. See \cite{AMLin} for further results and references.


\section{Further topics} \label{sec:topics}


\subsection{Symplectic realizations and dual pairs} \label{subsec:sympreal}
As previously seen, many interesting examples of Poisson manifolds arise as quotients of symplectic manifolds. Given the degenerate nature of Poisson structures and their singular symplectic foliations, it is often convenient to study Poisson manifolds by looking at symplectic manifolds that realize them as quotients.

Given a symplectic manifold $(S,\omega)$ and a surjective submersion $\varphi: S\to M$, a natural question is whether the symplectic structure on $S$ descends to a Poisson structure on $M$, or in other words if one can endow $M$ with a Poisson structure for which $\varphi$ is a Poisson map. The conditions for this to happen are described in Exercise~\ref{exer:lib1}, from where the so-called {\em Libermann's lemma} directly follows:
{\em If $M$ carries a Poisson structure such that $\varphi$ is a Poisson map then 
the symplectic orthogonal distribution to the subbundle $\mathcal{V}=\Ker(d\varphi)\subseteq TM$ is integrable, and the converse holds provided the fibers of $\varphi$ are connected.} 

The necessity of the connectedness assumption on $\varphi$-fibers is illustrated by the next exercise.

\begin{exer}
Consider the 2-sphere $\mathbb{S}^2$ equipped with a symplectic structure and the natural projection $\varphi: \mathbb{S}^2\to \mathbb{R}P^2$. Verify that $\varphi$ cannot be made into a Poisson map.
\end{exer}

A Poisson map $S\to M$ from a symplectic manifold $S$ to a Poisson manifold $M$ is called a {\em symplectic realization} of $M$. We call a symplectic realization {\em full} when it is a surjective submersion. We list below some examples that we have already encountered.


\begin{example}\label{ex:realiz}
\begin{itemize}
\item[(a)]  Inclusions of symplectic leaves into Poisson manifolds are symplectic realizations, generally not full. 

\item[(b)] The cotangent bundle projection $T^*M\to M$ is a full symplectic realization of the zero Poisson structure on $M$.

\item[(c)] As explained in $\S$ \ref{subsec:linear1}, symplectic realizations $\mu: S\to \mathfrak{g}^*$ of the dual of a Lie algebra $\mathfrak{g}$ are the same as moment maps of Hamiltonian $\mathfrak{g}$-actions on symplectic manifolds. In case of a free $\mathfrak{g}$-action, $\mu$ is a submersion, and hence a full symplectic realization of the open subset $\mu(S)\subseteq \mathfrak{g}^*$.

\item[(d)] When a Lie group $G$ acts freely and properly on a symplectic manifold $(S,\omega)$ by symplectomorphisms, the quotient map $S\to S/G$ is a full symplectic realization of the induced Poisson structure on $S/G$ (see $\S$ \ref{subsec:symmetry}). In particular, if  $\mathfrak{g}$ is the Lie algebra of $G$, then $T^*G$ is a full symplectic realization of $\mathfrak{g}^*=T^*G/G$ (see Exercise~\ref{exer:liepois}). 
\end{itemize}
\hfill $\diamond$
\end{example}

\begin{exer}\label{exer:leafcorrespond}
Let $\varphi: S\to M$ be a full symplectic realization with connected fibers. Let $\mathcal{V}=\ker(d\varphi)$, and consider the (smooth) distribution $\mathcal{V}+\mathcal{V}^\omega \subseteq TS$. Show that if $\mathcal{O}$ is a leaf of the symplectic foliation of $\pi$, then $\varphi^{-1}(\mathcal{O})$ is a leaf of $\mathcal{V}+\mathcal{V}^\omega$ (in particular, $\mathcal{V}+\mathcal{V}^\omega$ is integrable). Moreover, this map induces a homeomorphism between the leaf spaces of the corresponding foliations on $M$ and $S$. (C.f. $\S$ \ref{subsub:quot}.)
\end{exer}

Let $\varphi: S \to M$ be a full symplectic realization, and suppose that the foliation tangent to $\mathcal{V}^\omega$ is {\em simple}, in the sense that its leaves are given by fibers of a surjective submersion 
$\varphi': S\to M'$. By Libermann's lemma, $M'$ carries a Poisson structure such that $\varphi'$ is also a full symplectic realization.

Given a symplectic manifold $S$ and Poisson manifolds $M_1$ and $M_2$, a diagram
\begin{equation}\label{eq:dual}
 \xymatrix {
 & \ar[dl]_{\varphi_1} S \ar[dr]^{\varphi_2} &\\
M_1 &  & M_2}
\end{equation}
of full symplectic realizations satisfying $(\mathrm{ker}(d\varphi_1))^\omega=\mathrm{ker}(d\varphi_2)$ is called a {\em dual pair}.

For a key example, consider a free and proper Hamiltonian $G$-action on a symplectic manifold $S$, with moment map $\mu: S\to \mathfrak{g}^*$. Then (see Exercise~\ref{exer:reduction}) $S/G \leftarrow S \stackrel{\mu}{\to} \mu(\mathfrak{g}^*)$ is a dual pair.

Poisson manifolds $M_1$ and $M_2$ fitting into a dual pair \eqref{eq:dual} turn out to share many features.
\begin{itemize}
\item If $\varphi_1$- and $\varphi_2$-fibers are connected, then there is an induced bijection between symplectic leaves of $M_1$ and $M_2$ (in fact, a homeomorphism of leaf spaces) with respect to which $\mathcal{O}_1\hookrightarrow M_1$ corresponds to $\mathcal{O}_2\hookrightarrow M_2$ if and only if
$$
\varphi_1^{-1}(\mathcal{O}_1)=\varphi_2^{-1}(\mathcal{O}_2).
$$
This follows from Exercise~\ref{exer:leafcorrespond}.
In particular, there is an identification of Casimir functions on $M_1$ and $M_2$ (i.e., $H^0_{\pi}(M_1)\cong H^0_\pi(M_2)$).

\item For any $x \in S$, the transverse Poisson structures at $\varphi_1(x)$ and $\varphi_2(x)$ are anti-isomorphic (i.e., isomorphic up to a sign); see Exercise \ref{exer:dualprestransv} for a proof.

\item If $\varphi_1$- and $\varphi_2$-fibers are 1-connected (i.e., connected and simply-connected), then there is an induced isomorphism $H^1_{\pi}(M_1)\cong H^1_\pi(M_2)$ preserving modular classes \cite{Cr,Ginz,GinzLu}.

\end{itemize}

A dual pair $M_1 \stackrel{\varphi_1}{\leftarrow} S \stackrel{\varphi_2}{\to} {M}_2$  with the additional requirements that $\varphi_1$ and $\varphi_2$ have 1-connected fibers and   are complete Poisson maps (see Exercise \ref{exer:complete})  defines a {\em Morita equivalence} between $M_1$ and $\overline{M}_2$ \cite{XuMorita}. The notion of Morita equivalence of Poisson manifolds bears close resemblance with the original concept  of Morita equivalence in algebra, see e.g. \cite{BuWeMorita}; following this analogy,  the set of self-Morita equivalences of a Poisson manifold is called its {\em Picard group} \cite{BuFe,BuWePic}. Various examples of Morita equivalent Poisson manifolds and Picard groups are described in \cite{BuFe,BuWeMorita}; see \cite{BuWePic,Joel,RadShl} for the study of Morita equivalence and Picard groups of log-symplectic manifolds.


The fundamental question as to whether any Poisson manifold admits a full symplectic realization was answered in the affirmative by Weinstein \cite{CDW,We83} and Karasev \cite{Karas}.
The original proof consisted in gluing local symplectic realizations using a uniqueness property. More direct, global arguments to prove this result were given in \cite{CrMaReal} (see also \cite{FrejlichMarcut}).

One can take as the starting point to construct  a full symplectic realization $S\to M$ of an arbitrary Poisson manifold $M$ the fact that the zero Poisson structure on $M$ has a canonical symplectic realization  given by the cotangent bundle projection $q: (T^*M,\omega_{can}) \to M$. For a given Poisson structure $\pi$ on $M$, the intuitive idea is to obtain a symplectic realization by ``deforming'' $\omega_{can}$ using the flow of a suitably chosen vector field $Y$ on $T^*M$. This vector field is defined with the aid of 
the horizontal lift $\mathrm{hor}: q^*TM=T^*M\times_M TM \to T(T^*M)$ of a linear connection on $TM$ by 
$$
Y|_\xi=\mathrm{hor}(\xi,\pi^\sharp(\xi)), \qquad \xi \in T^*M.
$$ 
(The vector field $Y$ is more conceptually understood  as a {\em Poisson spray} \cite[$\S$ 1]{CrMaReal}, i.e., the Poisson-geometric analogue of a classical geodesic spray, see \cite[Chp.~11]{CFMbook}.)

Denoting the local flow of $Y$ by $\phi_t$, the fact that $Y$ vanishes on the zero section of $T^*M$ ensures that $\phi_t$ is defined for $t\in [0,1]$ in a sufficiently small neighborhood of $M\subseteq T^*M$. The main result \cite{CrMaReal} is that {\em there is a neighborhood $U \subseteq T^*M$ of the zero section  where the 2-form
\begin{equation}\label{eq:real}
\omega = \int_0^1 (\phi_t)^* \omega_{can} dt
\end{equation}
is symplectic and the restriction of the projection $q|_U: (U,\omega)\to (M,\pi)$ is a (full) symplectic realization.}

One can refine the problem of existence of full symplectic realizations of Poisson manifolds by the additional requirement that the symplectic realization be {\em complete}. 
For instance, the realizations of Example~\ref{ex:realiz} (d) are all complete (Exercise~\ref{exer:completequot}). As it turns out \cite{CF2}, it is {\em not} the case that any Poisson manifold admits a complete full symplectic realization. 

To hint at what is involved in this problem, consider the case of the dual of a Lie algebra. As mentioned in Example~\ref{ex:realiz} (d), any Lie group $G$ integrating a Lie algebra $\mathfrak{g}$ gives rise to a natural complete full symplectic realization of $M=\mathfrak{g}^*$ by $T^*G$. 
More generally, the existence of a complete full symplectic realization of a Poisson manifold is equivalent to the existence of a certain global object ``integrating'' the Poisson  manifold (discussed in $\S$ \ref{subsec:sympgrp} below); the issue is that, as opposed to Lie algebras, not every Poisson manifold is ``integrable'' in this sense. 

The fact that complete full symplectic realizations do not always exist implies that the notion of Morita equivalence mentioned above is not defined for all Poisson manifolds; on the other hand, among ``integrable'' Poisson manifolds, Morita equivalence is indeed an equivalence relation \cite{XuMorita} (see also \cite{BuWeMorita}).

\subsection{Symplectic groupoids}\label{subsec:sympgrp}

The full symplectic realization of a Poisson manifold $(M,\pi)$ defined by the symplectic form \eqref{eq:real} on a neighborhood $U\subseteq T^*M$ of the zero section via the cotangent bundle projection $\mathsf{t}:= q|_U: U\to  M$ comes naturally equipped with more structure:

\begin{itemize}
\item the map $\mathsf{s}:= q\circ \phi_1: U\to M$ is anti-Poisson, and  $M\stackrel{\mathsf{t}}{\leftarrow} U \stackrel{\mathsf{s}}{\rightarrow} \overline{M}$ is a dual pair;
\item the zero-section $M\hookrightarrow U$ is a lagrangian submanifold.
\end{itemize}
There is, in addition, a group-like multiplication on $U$, and understanding how these structures fit together leads to the notion of (local) {\em symplectic groupoid} \cite[Thm.~1.4]{CDW}.

Recall that a {\em Lie groupoid} consists of manifolds $\grd$ and $M$ equipped with structure maps given by surjective submersions $\tar,\sour:\grd\to M$ (called {\em
source and target maps}), a smooth partial {\em multiplication} map 
$$
m: \grd^{(2)}\to \grd, \qquad  (g,h)\mapsto gh,
$$
defined on the submanifold $\grd^{(2)}=\{(g,h)\,|\,\sour(g)=\tar(h)\}\subseteq \grd \times \grd$ of composable pairs, a diffeomorphism $i: \grd\to \grd$, $g\mapsto g^{-1}$, called {\em
inversion} map, and  an embedding $\varepsilon: M\to \grd$, $x\mapsto 1_x$, called {\em unit}
map, satisfying 
\begin{itemize}
\item[(i)]  {$\sour(gh)=\sour(h)$,
$\tar(gh)=\tar(g)$}, 
\item[(ii)] {$(gh)k=g(hk)$},
\item[(iii)] {$\sour(1_x)=\tar(1_x)=x$}, and
{$g 1_{\sour(g)}= 1_{\tar(g)}g = g$},
\item[(iv)] {$\sour(g^{-1})=\tar(g)$, $\tar(g^{-1})=\sour(g)$} and {$g^{-1}g=
1_{\sour(g)}$, $g g^{-1}= 1_{\tar(g)}$}.
\end{itemize}
We usually identify $M$ with its image in $\grd$ under $\varepsilon$.
It is convenient to think of an element $g\in \grd$ as an ``arrow'' $\tar(g) \stackrel{g}{\leftarrow} \sour(g)$ from its source to its target, so that e.g. the composition law (i) for the multiplication looks like
$$
m(z \stackrel{g}{\leftarrow} y, y \stackrel{h}{\leftarrow} x)= z \stackrel{gh}{\leftarrow} x.
$$
We use the notation $\grd \rightrightarrows M$ for a Lie groupoid. Morphisms of Lie groupoids are smooth maps preserving all groupoid structural maps.

A Lie group $G$ is the same as a Lie groupoid for which $M$ is a point, $G \rightrightarrows \{*\}$. Any manifold $M$ can be viewed as a Lie groupoid $M\rightrightarrows M$ with source and target maps being the identity. We will give more examples below.

A Lie groupoid $\grd\rightrightarrows M$ defines an equivalence relation on $M$ where $x\sim y$ if there exists an arrow $y \stackrel{g}{\leftarrow} x$. The equivalence class of $x\in M$,
$$
\mathcal{O}_x = \{y \in M \,|\, \exists\;  y \stackrel{g}{\leftarrow} x\}= \tar(\sour^{-1}(x)) \subseteq M,
$$
is called the {\em orbit} of $x$, and it is a (immersed) submanifold of $M$. To each $x\in M$ there also corresponds a Lie group 
$$
\grd_x = \{ g \in \grd \;|\;  x \stackrel{g}{\leftarrow} x \} = \tar^{-1}(x)\cap \sour^{-1}(x)\subseteq \grd,
$$
called the {\em isotropy group} at $x$. 
The decomposition of $M$ into orbits and the collection of isotropy groups are essential ingredients of a Lie groupoid.

\begin{exer}\label{exer:graph}
On a Lie groupoid $\grd\rightrightarrows M$, consider the graph of the multiplication map,
$$
\Gamma_m:=\{(g,h, gh),\; (g,h) \in \grd^{(2)} \}.
$$
Check that $\dim(\Gamma_m)= \dim(\grd^{(2)})= 2 \dim(\grd) - \dim(M)$.
\end{exer}

\begin{example}\label{ex:liegrds}
\begin{itemize}
\item[(a)] 
Given a manifold $M$, the {\em pair
groupoid} $M\times M\rightrightarrows M$ has arrows $y \stackrel{(y,x)}{\leftarrow} x$ and multiplication given by $(z,y)(y,x)= (z,x)$. This Lie groupoid has $M$ itself as the only orbit and trivial isotropy groups. A closely related Lie groupoid is the {\em fundamental groupoid} of $M$, $\Pi(M)\rightrightarrows M$, whose arrows are homotopy classes of
paths on $M$,
$$
\gamma(1) \stackrel{[\gamma]}{\leftarrow}{\gamma(0)},
$$
with multiplication given by concatenation. 
The
orbits of $\Pi(M)$ are the connected components
of $M$, and its isotropy group at $x$ is the fundamental group of the connected component containing $x$. The map $\Pi(M)\to M\times M$, $[\gamma]\mapsto (\gamma(1),\gamma(0))$ is a groupoid morphism, which is an isomorphism when
$M$ is 1-connected (i.e., connected with trivial fundamental group).

\item[(b)] A Lie groupoid for which source and target maps coincide is a {\em bundle of Lie groups}. Any vector bundle is a particular example of a Lie groupoid of this kind, with multiplication given by fiberwise addition.

\item[(c)] An action of a Lie group $G$ on a manifold $M$ gives rise to an {\em action groupoid} $G\times M \rightrightarrows M$, with arrows 
$$
\sigma x \stackrel{(\sigma,x)}{\leftarrow} x
$$ 
and multiplication $(\sigma',y)(\sigma, x)=(\sigma'\sigma, x)$.
The orbits and isotropy groups of the action groupoid
coincide with those of the action.
\end{itemize}
\hfill $\diamond$
\end{example}

On a Lie groupoid $\grd \rightrightarrows M$,  right and left multiplication by an arrow $y \stackrel{g}{\leftarrow} x$ define diffeomorphisms
$$
r_g:  \sour^{-1}(y) \to \sour^{-1}(x),  \qquad l_g: \tar^{-1}(x)\to \tar^{-1}(y).
$$

\begin{exer}\label{exer:grdid}
Use  properties (i), (iii) and (iv) of Lie groupoids to verify the following identities.
\begin{itemize}

\item[(a)] For $V\in T_g\grd$, 
$$
V = (dm)_{(\tar(g),g)}((d\tar)_g(V),V)= (dm)_{(g,\sour(g))}(V,(d\sour)_g(V)).
$$
In particular, $X=dm_{(x,x)}(X,X)$ for $X\in T_xM$.

\item[(b)] For $u\in \ker(d\sour)_y$, 
$$
(dr_g)_y(u)=dm_{(y,g)}(u,0), \quad \mbox{and} \quad (d\tar)_g((dr_g)_y(u))=(d\tar)_y(u).
$$

\item[(c)] For $V\in T_g\grd$, \; 
$
(d\tar)_g(V) = dm_{(g,g^{-1})}(V, (di)_{g}(V)).
$
\end{itemize}

\end{exer}

To define symplectic groupoids, we need to consider differential forms on Lie groupoids $\grd\rightrightarrows M$ suitably compatible with the groupoid structure. 
A $2$-form $\omega\in \Omega^2(\grd)$ is called {\em multiplicative} if
\begin{equation}\label{eq:mult}
m^*\omega = \pr_1^*\omega + \pr_2^* \omega,
\end{equation}
where $\pr_1, \pr_2: \grd^{(2)}\to \grd$ are the natural projections. (The same definition applies to differential forms of any degree.) A {\em symplectic groupoid} is a Lie groupoid equipped with a multiplicative symplectic form.

\begin{exer}\label{exer:omegaid}
Let $\omega$ be a multiplicative 2-form on $\grd\rightrightarrows M$. 
\begin{itemize}
\item[(a)] Use Exercise~\ref{exer:grdid}, (a) and (b), to show that, for $y\stackrel{g}{\leftarrow} x$, $u\in (\ker(d\sour))_y$, and $V\in T_g\grd$,
$$
\omega((dr_g)_y(u),V)=\omega(u,(d\tar)_g(V)).
$$
Using the isomorphism $(dr_g)_y: \ker(d\sour)_y \stackrel{\sim}{\to} \ker(d\sour)_g$, conclude that 
$$
\ker(d\sour)\subseteq \ker(d\tar)^\omega.
$$
\item[(b)] For $X, Y\in T_xM$, check that $\omega(X,Y)=0$ (see Exercise~\ref{exer:grdid} (a)). 

\item[(c)] Use the previous item and Exercise~\ref{exer:grdid} (c) to show that, for $V,W\in T_g\grd$,
$$
\omega(V,W) + \omega((di)_g(V),(di)_g(W))=0.
$$
\end{itemize}
\end{exer}

On a symplectic groupoid $(\grd, \omega)$, the compatibility condition \eqref{eq:mult} is equivalent to requiring that the graph of the multiplication map, $\Gamma_m:=\{(g,h, gh),\; (g,h) \in \grd^{(2)} \}$, is an isotropic submanifold of $(\grd,\omega)\times(\grd,\omega)\times (\grd,-\omega)$. On the other hand, by Exercise~\ref{exer:omegaid} (b), the image of $M\stackrel{\varepsilon}{\hookrightarrow} \grd$ is isotropic. Therefore (see Exercise~\ref{exer:graph})
$$
\dim(\Gamma_m)=2 \dim(\grd)-\dim(M)\leq \frac{3}{2} \dim(\grd), \qquad \dim(M)\leq \frac{1}{2}\dim(\grd).
$$
Combining these two inequalities, we have that
$$
\dim(M)=\frac{1}{2} \dim(\grd), \quad \mbox{ and } \quad \dim(\Gamma_m)=\frac{3}{2}\dim(\grd).
$$
As a consequence, we see that {\em a symplectic form $\omega$ on a Lie groupoid $\grd$ is multiplicative if and only if $\Gamma_m$ is a lagrangian submanifold of $(\grd,\omega)\times(\grd,\omega)\times (\grd,-\omega)$}; this last condition is the original definition of symplectic groupoid in \cite{We87}.

The following are some important properties of a symplectic groupoid $(\grd,\omega)$.

\begin{itemize}
\item[(SG1)] {\em The submanifold of units $M\hookrightarrow \grd$ is lagrangian}, as shown above.
\item[(SG2)] {\em The inversion map satisfies $i^*\omega=-\omega$}, by Exercise~\ref{exer:omegaid} (c).
\item[(SG3)] $\ker(d\sour)=\ker(d\tar)^\omega$ (by Exercise~\ref{exer:omegaid} (a), since the subbundles $\ker(d\sour)$ and $\ker(d\tar)^\omega$ have the same rank).  
\end{itemize}

Another key property, linking symplectic groupoids to Poisson structures, will follow from the next exercise (see (SG4) below).

\begin{exer}\label{exer:unitpois} Let $(\grd,\omega)$ be a symplectic groupoid, and
consider the vector bundle map 
$$
\mu_\omega: \ker(d\sour)|_M \to T^*M,  \;\;\;\;\; \mu_\omega(u) =\omega(u,\cdot)|_{TM}.
$$
\begin{itemize}
\item[(a)] Check that $\mu_\omega$ is an isomorphism. 
\item[(b)] Let $y\stackrel{g}{\leftarrow} x$,  $u\in \ker(d\sour)|_y$ and $f\in C^\infty(M)$  such that $df|_y= \mu_\omega(u)$. Conclude from Exercise~\ref{exer:omegaid} (a) that
$$
X_{\tar^*f}|_g= (dr_g)_y(u).
$$

\item[(c)] For $f_1, f_2 \in C^\infty(M)$, show that
$$
\{\tar^*f_1,\tar^*f_2\}(g)= \omega(u_2, (d\tar)_y(u_1)),
$$
where $y=\tar(g)$, $u_i\in (\ker(d\sour))_y$, and $\mu_\omega(u_i)=df_i|_y$, $i=1,2$.
Note in particular that the function $\{\tar^*f_1,\tar^*f_2\}\in C^\infty(\grd)$ is constant along $\tar$-fibers.
\end{itemize}

\end{exer}

A consequence of the previous exercise (part (c)) is that, on a symplectic groupoid, 
$$
\{\tar^*C^\infty(M),\tar^*C^\infty(M)\}\subseteq \tar^*C^\infty(M),
$$
and therefore $M$ inherits a Poisson structure such that $\tar$ is a Poisson map.
(In the special case of symplectic groupoids with $\tar$-connected fibers, this property can also be deduced from (SG3) and Libermann's lemma, see $\S$ \ref{subsec:sympreal}).

\begin{exer}\label{exer:targetcomplete}
Let $(\grd,\omega)$ be a symplectic groupoid and $f\in C^\infty(M)$.
It follows from Exercise~\ref{exer:unitpois} (b) that the Hamiltonian vector field $X_{\tar^*f}$ is tangent to the $\sour$-fibers, and,  for $y\stackrel{g}{\leftarrow} x$,
$$
(r_g)_*(X_{\tar^*f }|_{\sour^{-1}(y)}) = X_{\tar^*f}|_{\sour^{-1}(x)}.
$$
Use this last property to show that the Poisson map $\tar: \grd \to M$ is complete.
\end{exer}

In conclusion, for a symplectic groupoid $(\grd,\omega)$,
\smallskip
\begin{itemize}
\item[(SG4)]
{\em there exists a (unique) Poisson structure on $M$ making $\tar: \grd\to M$ into a complete, full symplectic realization.}
\end{itemize}

Note that, by the discussion in $\S$ \ref{subsub:quot} (see \eqref{eq:diag}), the symplectic leaves of the Poisson structure on $M$ are the connected components of the orbits of $\grd\rightrightarrows M$.


\begin{example}\label{ex:sympgrd}
\begin{itemize}

\item[(a)] Given a symplectic manifold $(S,\omega)$, the pair groupoid $S\times S \rightrightarrows S$ is a symplectic groupoid with respect to $\pr_1^*\omega -\pr_2^*\omega$. The induced Poisson structure on $S$ is the one defined by $\omega$.
Since the morphism $\Pi(S)\to S\times S$ is a local diffeomorphism, the pullback of 
$\pr_1^*\omega -\pr_2^*\omega$ makes $\Pi(S)$ into a symplectic groupoid over $S$.

\item[(b)] Let $E\to M$ be a vector bundle, viewed as a Lie groupoid with multiplication given by fiberwise addition (see Example~\ref{ex:liegrds} (b)). A symplectic form $\omega\in \Omega^2(E)$ is multiplicative if and only if it defines a linear Poisson structure, in the sense of $\S$ \ref{subsec:linear2}.
It follows from Exercise~\ref{exer:linarcanonical} that $(E,\omega)$ is isomorphic to $T^*M$
with its canonical symplectic form $\omega_{can}$. The induced Poisson structure on $M$ is trivial.

\item[(c)] For a Lie group $G$, the action groupoid $(G\times \mathfrak{g}^*)\rightrightarrows \mathfrak{g}^*$ (with respect to the coadjoint action) becomes a symplectic groupoid when equipped with the symplectic form obtained from the identification $T^*G\simeq G\times \mathfrak{g}^*$ via left (or right) translations.  The induced Poisson structure on $\mathfrak{g}^*$ is the linear Poisson structure dual to the Lie algebra $\mathfrak{g}$ (up to a sign, depending on the chosen identification). 
\end{itemize}
\hfill $\diamond$
\end{example}

A Poisson manifold is called {\em integrable} if it is isomorphic to the manifold of units of a symplectic groupoid; in such a case one refers to the symplectic groupoid  as an {\em integration} of the Poisson manifold. The examples above show that Poisson structures that are nondegenerate (symplectic), trivial or linear are integrable, and also that integrations are not unique.  On the other hand,  as first observed by Weinstein
\cite{We87}, there are Poisson manifolds that fail to be integrable; a concrete example is 
the Poisson structure on $M=\mathbb{S}^2\times \mathbb{R}$ defined by the family of symplectic structures $\omega_t= (1+t^2)\omega_{area}$ on $\mathbb{S}^2$, as in Example~\ref{ex:family} (see \cite{CatFel,CF2,We87}).
This phenomenon can be understood in light of the broader integrability problem for Lie algebroids \cite{CF}, briefly recalled in $\S$ \ref{subsec:cotangent} below.

We mention other important classes of integrable Poisson manifolds.

\begin{itemize}
    \item If $M$ is an integrable Poisson manifold endowed with a free and proper action of a Lie group $G$ by Poisson automorphisms, then the Poisson manifold $M/G$ (see $\S$ \ref{subsec:symmetry}) is also integrable. As shown in \cite{FOR}, $M$ admits an integration by a symplectic groupoid $\grd$ carrying a Hamiltonian $G$-action lifting the $G$-action on $M$, and its symplectic reduction at zero
    is a symplectic groupoid over $M/G$. The simplest instance of this construction is when $M$ has the trivial Poisson structure and $\grd=T^*M$; in this case, it is well known that the canonical lift of  the $G$-action on  $M$ to $T^*M$ is Hamiltonian (see Exercise~\ref{exer:linarcanonical}) and its reduction at zero is $T^*(M/G)$.

    \item Any Poisson Lie group $(G,\pi)$ is integrable \cite{LuWedouble}. The  symplectic groupoids integrating $G$ described in \cite{LuWedouble} have the property of carrying a second groupoid structure making them into a symplectic groupoid over the dual Poisson Lie group $G^*$, 
$$
\begin{matrix}
(\grd,\omega) & \rightrightarrows & G^* \\ \downdownarrows &  & \downdownarrows\\ G & \rightrightarrows & \ast.
\end{matrix}
$$
These are examples of {\em symplectic double groupoids}. The simplest case is given by $(T^*G,\omega_{can})$, which can be regarded as a symplectic groupoid over $G$ (with the trivial Poisson structure) and over $\mathfrak{g}^*$, see Example~\ref{ex:sympgrd}, (b) and (c).

More generally, any Poisson homogeneous space is integrable, see \cite{BIL}.  
    \smallskip
    \item Log-symplectic Poisson manifolds are integrable. Their symplectic groupoids were constructed and classified in \cite{GuaLi} by means of blow-up and gluing operations.
\end{itemize}

For an integrable Poisson manifold $(M,\pi)$, there is a canonical
symplectic groupoid $\Pi(M,\pi)\rightrightarrows M$ that generalizes the fundamental groupoid of a symplectic manifold (Example~\ref{ex:sympgrd} (a)), obtained as a quotient of ``cotangent paths'' on $M$ by ``cotangent homotopies'' \cite{CatFel,CF2}. We refer to $\Pi(M,\pi)$ as the 
{\em Poisson fundamental groupoid} of $(M,\pi)$. 
For a general Poisson manifold $(M,\pi)$, the construction of $\Pi(M,\pi)$ produces a {\em topological} groupoid which is always smooth on an open neighborhood of the units $M\subseteq \Pi(M,\pi)$, yielding a ``local'' symplectic groupoid (its global smoothness occurs if and only if $M$ is integrable \cite{CF2}).

We outline the key points relating symplectic groupoids to the study of symplectic realizations in $\S$ \ref{subsec:sympreal}.
The general symplectic realization of a Poisson manifold $(M,\pi)$ defined by the 2-form \eqref{eq:real} on a neighborhood of the zero section of $T^*M$  is isomorphic to the symplectic realization given by the target map on a smooth neighborhood of the units in $\Pi(M,\pi)$.  On the other hand, by (SG4) any integrable Poisson manifold has a complete, full symplectic realization given by the target map of a symplectic groupoid. More generally, as proven in \cite{CF2}, {\em the existence of a complete, full symplectic realization of a Poisson manifold is  equivalent to the existence of a symplectic groupoid that integrates it} (see also \cite{daniel}).

\begin{remark}[More on symplectic groupoids]
\
\begin{itemize}
\item[$\diamond$] Symplectic groupoids were originally introduced as part of a  quantization scheme for Poisson manifolds \cite{WeGQNC,WeXuext,Karas} (see \cite{Hawkins,bonechi} for more recent developments) based on the fact that groupoids are rich sources of noncommutative algebras via their convolution algebras.

\item [$\diamond$] 
The Poisson fundamental groupoid of a Poisson manifold $(M,\pi)$ can be regarded as an infinite-dimensional symplectic quotient, and it has an interpretation in terms of (classical) topological field theory as the phase space of the so-called {\em Poisson sigma model} with target $M$ \cite{CatFel}.

\item[$\diamond$] Symplectic groupoids and their actions lead to a broad extension of the theory of Hamiltonian actions on symplectic manifolds, which allows symplectic realizations to be regarded as ``moment maps''  \cite{MiWe} (see also \cite[$\S$ 14.4]{CFMbook}). 
Symplectic groupoid actions are also key elements in the theory of Morita equivalence of Poisson manifolds \cite{XuMoritaSymp,XuMorita} (see also \cite{BuWeMorita}). 

\item[$\diamond$] A symplectic groupoid gives rise to a dual pair
\begin{equation}\label{eq:unitpic}
\xymatrix {
 & \ar[dl]_{\tar} \grd \ar[dr]^{\sour} &\\
M &  & \overline{M},}
\end{equation}
see properties (SG2), (SG3), and (SG4).
Moreover, if $\grd$ has 1-connected $\sour$-fibers, this dual pair defines a self Morita equivalence of the Poisson manifold $M$, which is the unit element of its Picard group.

\item
    [$\diamond$] Symplectic groupoids serve as a tool e.g. for
    the calculation of Poisson cohomology \cite{XuCoh} and linearization problems \cite{CrFeLin}. Compactness properties of symplectic groupoids have many relevant effects on their underlying Poisson manifolds; their systematic study has recently begun in \cite{PMCT1}.

\end{itemize}
\hfill $\diamond$
\end{remark}

\begin{remark}\label{eq:cour}
{\em Poisson groupoids} \cite{WeinCoiso} provide a common framework for
symplectic groupoids and Poisson Lie groups ($\S$ \ref{subsec:PL}). The extension of   Lie bialgebras  and their Drinfeld doubles (see Remark~\ref{rem:double}) to this context leads to the important notions of {\em Lie bialgebroids} and {\em Courant algebroids} \cite{LWX}. \hfill $\diamond$
\end{remark}



\subsection{The Lie algebroid of a Poisson manifold}\label{subsec:cotangent}

Lie algebroids (see $\S$ \ref{subsec:linear2}) are infinitesimal versions of Lie groupoids, extending the way Lie algebras arise from Lie groups. A key fact explaining the link between Poisson manifolds and symplectic groupoids is that any Poisson manifold has an associated Lie algebroid. From this perspective, the correspondence between Poisson manifolds and symplectic groupoids is analogous to that of Lie algebras and Lie groups.

The Lie algebroid of a Lie groupoid $\grd\rightrightarrows M$ is defined by as follows\footnote{An alternative convention reverses the roles of source and target maps, and uses left instead of right-invariant vector fields.}:
\begin{itemize}
\item As a vector bundle, $A= \mathrm{ker}(d\sour)|_M \to M$;
\item The anchor map is the restriction of $d\tar: T\grd\to TM$ to $A\subseteq T\grd|_M$,
$$
\rho = d\tar |_A : A \to TM.
$$
\item The Lie bracket on the space of sections $\Gamma(A)$ arises from
its identification with vector fields $X$ on $\grd$ that are tangent to $\sour$-fibers and right-invariant, i.e., $d r_g (X_h) = X_{hg}$.
\end{itemize}

A Lie algebroid $A\to M$ is {\em integrable} if it is isomorphic to the Lie algebroid of a Lie groupoid $\grd\rightrightarrows M$, as above. In this case, we will say that $\grd$ {\em integrates}  $A$.

Unlike Lie algebras, which can be always integrated to Lie groups, Lie algebroids are not necessarily integrable (see \cite{CF} for a description of obstructions to integrability and examples). Nonetheless, other fundamental results in classical Lie theory are still valid in this broader context; e.g., the construction of the Lie algebroid of a Lie groupoid gives rise to a {\em Lie functor} from Lie groupoids to Lie algebroids that defines an {\em equivalence of categories between source-simply connected Lie groupoids\footnote{We say that a Lie groupoid is source-simply connected if its $\sour$-fibers are 1-connected.} and integrable Lie algebroids}.
(In particular, an integrable Lie algebroid has unique source-simply connected integration, up to isomorphism.) 

Given a Poisson manifold $(M,\pi)$, its cotangent bundle $A=T^*M$ carries a Lie algebroid structure with anchor map $\pi^\sharp: T^*M \to TM$ and Lie bracket on $\Gamma(A)=\Omega^1(M)$ given by
$$
[\alpha,\beta]_\pi:= \Lie_{\pi^\sharp(\alpha)}\beta - \Lie_{\pi^\sharp(\beta)}\alpha - d (\pi(\alpha,\beta)),
$$
so  that $[df,dg]_\pi=d\{f,g\}$. 

\begin{exer}
Show that a Lie algebroid structure on $T^*M$, with anchor $\rho: T^*M\to TM$ and bracket $[\cdot,\cdot]$ on $\Omega^1(M)$, is induced by a Poisson structure on $M$ as above if and only if $\rho^*=-\rho$ and $[\cdot,\cdot]$ preserves the subspace $\Omega^1_{cl}(M)$ of closed 1-forms. (As a key step, verify that if $\Lambda: \Omega^1(M)\times \Omega^1(M)\to \Omega^1(M)$ is $C^\infty(M)$-bilinear, skew-symmetric, and $\Lambda(\Omega^1_{cl}(M),\Omega^1_{cl}(M))\subseteq \Omega^1_{cl}(M)$, then $\Lambda= 0$.
)
\end{exer}

\begin{exer} Given a Lie algebroid $A\to M$, a subbundle $B\to N$ is a {\em Lie subalgebroid} if $\rho(B)\subseteq TN$ and $[u,v]|_N\in \Gamma(B)$ whenever $u|_N$, $v|_N \in \Gamma(B)$.

 Let $(M,\pi)$ be a Poisson manifold, and consider the corresponding Lie algebroid $T^*M$. Show that a submanifold $N\hookrightarrow M$ is  coisotropic if and only if $\mathrm{Ann}(TN)\subseteq T^*M$ is a Lie subalgebroid. Hence there is a natural bijective correspondence between coisotropic submanifolds of $M$ and lagrangian Lie subalgebroids of $T^*M$ (i.e., Lie subalgebroids that are lagrangian submanifolds with respect to the canonical symplectic form). 
 \end{exer}

\begin{remark}\label{rem:Liealgprop}
Many properties of Poisson manifolds can be understood in terms of general features of Lie algebroids. 
For a Lie algebroid $A\to M$, 
\begin{itemize}
\item the distribution  $\rho(A)\subseteq TM$ given by the image of the anchor map is  integrable; the leaves of the corresponding (singular) foliation on $M$ are called {\em orbits} of $A$; 
\item for each $m\in M$, $\Ker(\rho|_m) \subseteq A|_m$ acquires a Lie algebra structure by restriction of the Lie bracket on $\Gamma(A)$, called the {\em isotropy Lie algebra} at $m$.
\end{itemize}
When $A=T^*M$ is the cotangent Lie algebroid of a Poisson structure $\pi$ on $M$, orbits of $A$ are the leaves of the symplectic foliation (see $\S$ \ref{subsec:integ}), and the isotropy Lie algebra of $A$ at $m$ is the transverse Lie algebra of $\pi$ at $m$ (see $\S$ \ref{subsec:transverse}). The Lie algebroid cohomology of $A$ (see $\S$ \ref{subsec:log}) agrees with the Poisson cohomology of $(M,\pi)$. \hfill $\diamond$
\end{remark}

\begin{remark} For a Poisson manifold $(M,\pi)$, the Lie algebroid structure on $T^*M$  gives rise to a linear Poisson structure $\pi_{TM}$ on $TM$ (see
$\S$ \ref{subsec:linear2}) uniquely characterized by the property that
$$
\{df,dg\}_{TM}= d\{f,g\},
$$
for all $f,g \in C^\infty(M)$, where here we view $df$, $dg$ as (linear) functions on $TM$. The Poisson structure $\pi_{TM}$ on $TM$ is called the {\em tangent lift} of $\pi$.

\begin{exer}
Check that if $(M,\pi)$ is symplectic, then so is its tangent lift $\pi_{TM}$ on $TM$. How is this symplectic structure on $TM$ related to the canonical symplectic structure on $T^*M$? 
\end{exer}
\hfill $\diamond$
\end{remark}

Symplectic groupoids are closely related to integrations of cotangent Lie algebroids of Poisson manifolds.

First, suppose that $(M,\pi)$ is the Poisson manifold of units of a symplectic groupoid $(\grd,\omega)$, and let $A$ be the Lie algebroid of $\grd$. The next exercise shows that the Lie algebroid $T^*M$ is naturally identified with $A$.

\begin{exer} Following Exercise~\ref{exer:unitpois}, consider the vector bundle isomorphism
$$
\mu_\omega: A \to T^*M,  \;\;\;\;\; \mu_\omega(u) =\omega(u,\cdot)|_{TM},
$$
and note that the Poisson structure $\pi$ on $M$ satisfies $\pi^\sharp=\rho \circ \mu_{\omega}^{-1}: T^*M\to TM$.
Verify that $\mu_\omega$ is an isomorphism of Lie algebroids, i.e., it intertwines anchor maps and preserves brackets.
\end{exer}

 Therefore any symplectic groupoid integrates the cotangent Lie algebroid of the Poisson structure on its units.  In the converse direction, if $M$ is a Poisson manifold with integrable cotangent Lie algebroid, then any source-simply-connected Lie groupoid integrating $T^*M$ is naturally a symplectic groupoid \cite[Thm.~5.2]{MacXu} (but this may not be the case for other integrations, see Remark~\ref{rem:nonint} below); moreover, the Poisson fundamental groupoid $\Pi(M,\pi)$ gives a concrete model for source-simply-connected integrations of $T^*M$ \cite{CF2}. 
 
 We conclude that {\em a Poisson manifold is integrable (in the sense of $\S$ \ref{subsec:sympgrp}) if and only if its cotangent Lie algebroid is integrable}.
 
 \begin{remark} \label{rem:nonint} 
Consider the $3$-sphere $\mathbb{S}^3$ equipped with the zero Poisson structure. As a Poisson manifold, it is integrated by $T^*\mathbb{S}^3=\mathbb{S}^3\times \mathbb{R}^3$, viewed as a symplectic groupoid with respect to fiberwise addition and canonical symplectic structure (see Example~\ref{ex:sympgrd} (b)).
In this case the cotangent Lie algebroid of $\mathbb{S}^3$ has trivial anchor and bracket. The Lie groupoid  $\mathbb{S}^3\times \mathbb{T}^3$ (viewed as a bundle of abelian Lie groups over $\mathbb{S}^3$) is another integration of this Lie algebroid, but it  cannot be a symplectic groupoid (since the manifold $\mathbb{S}^3\times \mathbb{T}^3$ carries no symplectic form). \hfill $\diamond$
\end{remark}


\subsection{Dirac structures}\label{subsec:dirac} 
On a symplectic manifold $M$, any submanifold $N\hookrightarrow M$ inherits a closed 2-form (usually not symplectic) via pullback of the ambient symplectic form. Dirac structures are geometric objects that give a similar description of the intrinsic geometry inherited by submanifolds of Poisson manifolds. The slogan one should have in mind is that  {\em  Dirac structures are to Poisson structures what closed 2-forms are to symplectic structures}. In what follows, we will also use the term ``presymplectic'' to refer to closed 2-forms.

Dirac structures have several applications to Poisson geometry and related areas, some of which are mentioned at the end of this subsection.  Here we will mostly focus on their role in the study of submanifolds of Poisson manifolds. Besides Courant's original paper \cite{courant}, general references include \cite{bursztynDirac} and \cite[Chp.~7]{CFMbook}.

A peculiarity of Dirac structures is that, in contrast with Poisson and symplectic structures, their general definition is not in terms of tensors on a manifold $M$, but rather as special types of vector subbundles of the direct sum 
$$
\TM:= TM\oplus T^*M.
$$
The vector bundle $\TM$ carries a natural nondegenerate, symmetric, fiberwise bilinear pairing $\SP{\cdot,\cdot}$ given by
$$
\SP{(X,\alpha),(Y,\beta)}=\beta(X)+\alpha(Y).
$$
With respect to this pairing, a subbundle $L\subseteq \TM$ is {\em lagrangian} (i.e., $L=L^\perp$) if and only if $\SP{\cdot,\cdot}$ vanishes on $L$ and $\mathrm{rank}(L)=\dim(M)$. As shown in the next exercise, 2-forms and bivector fields on $M$ can be characterized by appropriate lagrangian subbundles of $\TM$.

\begin{exer} \label{exer:kerL}
 Verify the following.
\begin{itemize}
\item[(a)] For a 2-form $\omega\in \Omega^2(M)$, the graph of $\omega^\flat:TM\to T^*M$ is a lagrangian subbundle of $\TM$, and this establishes a bijection between 2-forms on $M$ and lagrangian subbundles of $\TM$ that intersect $0\oplus T^*M$ trivially.
\item[(b)] The map taking a bivector field $\pi\in \mathfrak{X}^2(M)$ to the graph of $\pi^\sharp:T^*M\to TM$ is a bijection between bivector fields on $M$ and lagrangian subbundles of $\TM$ intersecting $TM\oplus 0$ trivially.
\end{itemize}
\end{exer}

To identify 2-forms that are closed and bivector fields that are Poisson, one considers the {\em Courant-Dorfman bracket} $\Cour{\cdot,\cdot}$ on $\Gamma(\TM)=\mathfrak{X}^1(M)\oplus \Omega^1(M)$ \cite{courant,Dorfman},
\begin{equation}\label{eq:CDbrk}
\Cour{(X,\alpha),(Y,\beta)} = ([X,Y],\Lie_X\beta - i_Y
d\alpha).
\end{equation}

A {\em Dirac structure} on $M$ \cite{courant,CouWe} is a lagrangian subbundle $L\subset \mathbb{T}M$ satisfying 
$$
\Cour{\Gamma(L),\Gamma(L)}\subset \Gamma(L),
$$
referred to as the {\em integrability condition}. The pair $(M,L)$ is a {\em Dirac manifold}.

\begin{exer}
Given a lagrangian subbundle $L\subseteq \TM$, verify that the operation
$$
(e_1,e_2,e_3)\mapsto \SP{\Cour{e_1,e_2},e_3}, \quad \text{ for }\; e_1,e_2, e_3 \in \Gamma(L),
$$
is skew-symmetric and $C^\infty(M)$-trilinear, so it defines an element $\Upsilon_L\in \Gamma(\wedge^3L^*)$, called the {\em Courant tensor} of $L$. Moreover, the integrability of $L$ is equivalent to the vanishing condition $\Upsilon_L=0$. 
\end{exer}


\begin{exer}
\begin{itemize}

\item[(a)] Let $L\subseteq \TM$ be given by the graph of $\omega^\flat$, for a 2-form $\omega\in \Omega^2(M)$. Given $e_i=(X_i, \omega^\flat(X_i)) \in \Gamma(L)$, $i=1,2,3$, check that
$\Upsilon_L(e_1,e_2,e_3)= d\omega(X_1,X_2,X_3)$. Therefore $L$ is integrable (i.e., it is a Dirac structure) if and only $d\omega=0$.

\item[(b)] Let $L\subseteq \TM$ be the graph of $\pi^\sharp$, for a bivector field $\pi\in \mathfrak{X}^2(M)$. For $e_i=(X_{f_i},df_i) \in \Gamma(L)$, $i=1,2,3$, check that
$\Upsilon_L(e_1,e_2,e_3)= \mathrm{Jac}(f_1,f_2,f_3)$. Hence $L$ is integrable if and only if $\pi$ is Poisson.

\item[(c)] Any subbundle $D\subseteq TM$ gives rise to a lagrangian subbundle $L=D\oplus \mathrm{Ann}(D)\subseteq \TM$. Show that $L$ is a Dirac structure if and only if $D$ is involutive (and hence the tangent distribution to a regular foliation on $M$).
\end{itemize}

\end{exer}

We define the {\em kernel} distribution of a Dirac structure $L$ on $M$ by 
$$
\ker(L)= L\cap (TM\oplus 0) \subseteq TM,
$$
so that Poisson structures on $M$ are equivalent to Dirac structures with trivial kernel, see Exercise~\ref{exer:kerL} (b).  (For Dirac structures defined by closed 2-forms, those with trivial kernel correspond to symplectic forms.) As a consequence, a general Dirac structure $L$ defines a Poisson structure on the (possibly empty) open subset of $M$ where its kernel is trivial.

The next exercise illustrates how ``singular'' Poisson structure may be encoded by (smooth!) Dirac structures.

\begin{exer} On $M= \mathbb{R}^3 = \{ (x,y,z) \}$, check that the subbundle of $\TM$ spanned by the sections $(\frac{\partial}{\partial y}, z dx)$, $(\frac{\partial}{\partial x}, -z dy)$, $(0,dz)$  is a Dirac structure. Verify that this Dirac structure has trivial kernel on the open subset $\{(x,y,z)\,|\, z\neq 0\}$, where it corresponds to the Poisson structure $\pi =
\frac{1}{z}\frac{\partial}{\partial x}\wedge
\frac{\partial}{\partial y}$.
\end{exer}

The exercise below shows how the abelian group of closed 2-forms on a manifold acts on the set of Dirac structures by means of an operation called {\em gauge transformation} \cite{SW}.

\begin{exer}\label{exer:Bfield2}
For a 2-form $B\in \Omega^2(M)$, consider the operation 
$$
\tau_B: \mathbb{T}M \to \mathbb{T}M, \;\; (X,\alpha)\mapsto (X, i_XB + \alpha).
$$ 
\begin{itemize}
\item[(a)] Check that $\tau_B$ preserves the natural symmetric pairing on $\mathbb{T}M$ and that it 
preserves the Courant-Dorfman bracket if and only if $dB=0$. In this case, if $L$ is a Dirac structure, so is $\tau_B(L)$. 

\item[(b)] Suppose that $\pi$ is a Poisson structure on $M$ and $L_\pi =\mathrm{graph}(\pi^\sharp)$ is the corresponding Dirac structure. For a closed 2-form $B$, show that the Dirac structure $\tau_B(L_\pi)$ is again given by a Poisson structure if and only if $(\Id + B^\flat \circ \pi^\sharp): T^*M \to T^*M$ is an isomorphism, in which case $\tau_B(L_\pi) =  L_{\pi_B}$ (see Exercise~\ref{exer:Bfield}).
\end{itemize}
\end{exer}

\begin{remark}[More on the Courant-Dorfman bracket]
The Courant-Dorfman bracket \eqref{eq:CDbrk} is an interesting object in its own right.
The following are its main properties:
\begin{enumerate}[(i)]
\item $\Cour{e_1,e_2}+\Cour{e_2,e_1}= d \SP{e_1,e_2}$,
\item $\Cour{e_1,\Cour{e_2,e_3}}=\Cour{\Cour{e_1,e_2},e_3} + \Cour{e_2,\Cour{e_1,e_3}}$,
\item $\Cour{e_1,fe_2}= (\mathcal{L}_{X_1}f)e_2 +f \Cour{e_1,e_2}$,
\item$\mathcal{L}_{X_1}\SP{e_2,e_3}=\SP{\Cour{e_1,e_2},e_3} + \SP{e_2,\Cour{e_1,e_3}}$,
\end{enumerate}
for $e_1,e_2,e_3\in \Gamma(\TM)$, and $e_1=(X_1,\alpha_1)$. The Courant-Dorfman bracket is not skew-symmetric, hence not a Lie bracket, but we see that it satisfies versions of the Jacobi and Leibniz identities (as well as a natural compatibility with the pairing). This is an example of a {\em Courant algebroid}, called the {\em standard Courant algebroid}. One can alternatively consider the skew-symmetrization of this bracket (which is Courant's original bracket in \cite{courant}), but this has the effect of introducing ``anomalies'' in properties (ii), (iii) and (iv) (that can be suitably understood in terms of ``Lie 2-algebras'' \cite{Roy,RoyWe}).
By (i) the Courant-Dorfman bracket coincides with its skew symmetrization on sections of lagrangian subbundles $L\subseteq \TM$, so the integrability condition of Dirac structures can be phrased using either one of the brackets.
\hfill $\diamond$
\end{remark}

We mention some of the key features of Dirac manifolds extending properties of Poisson manifolds.

\begin{itemize}
\item {\em (Lie algebroid)} For a Dirac structure $L$ on $M$, the vector bundle $L\to M$ acquires a Lie-algebroid structure with anchor given by the natural projection $L\to TM$ and bracket on $\Gamma(L)$ given by the restriction of the Courant-Dorfman bracket. When $L$ corresponds to a Poisson structure, the projection $\TM\to T^*M$ induces a Lie algebroid isomorphism $L\stackrel{\sim}{\to} T^*M$.

\item {\em (Characteristic distribution and presymplectic foliation)} The image of a lagrangian subbundle $L\subseteq \TM$ under the projection $\pr_{TM}: \TM\to TM$ defines a distribution  
$$
R:=\pr_{TM}(L)\subseteq TM
$$ 
that generalizes the characteristic distribution of a bivector field,  see $\S$ \ref{subsec:char}; additionally, over each point $x\in M$, $L$ defines a skew-symmetric bilinear form $\Omega_x$ on $R|_x$  by
$$
\Omega_x(X,Y)= \alpha(Y),
$$
for $X,Y\in R|_x$ and $\alpha$ any element in $T_x^*M$ such that $(X,\alpha)\in L|_x$.
\begin{exer}
Verify that $\mathrm{Ann}(R)=L\cap (0\oplus T^*M)$, and check that $\Omega_x$ is well defined. 
\end{exer}
When $L$ is a Dirac structure, the distribution $R$ is integrable and each leaf $\mathcal{O}\hookrightarrow M$,
$$
T\mathcal{O}=R|_\mathcal{O},
$$
carries a closed 2-form $\omega_\mathcal{O}$ defined  pointwise by the 2-forms $\Omega_x$ on $R|_x$, $x\in \mathcal{O}$.
The collection of all such presymplectic leaves defines the {\em presymplectic foliation} of $L$, which uniquely characterizes $L$.
Note that $\omega_\mathcal{O}$ is nondegenerate at $x\in M$ if and only if $\ker(L)|_x=\{0\}$, so 
the presymplectic foliation of $L$ has symplectic leaves if and only if $L$ has trivial kernel at all points, which happens if and only if $L$ is defined by a Poisson structure.

\item {\em (Integration)} The ``differentiation/integration'' type of correspondence relating Poisson manifolds and symplectic groupoids also extends to Dirac manifolds \cite{BCWZ}. The global objects arising in this more general setting are known as {\em presymplectic groupoids} (or {\em quasi-symplectic groupoids} \cite{XuQuasi}) and consist of Lie groupoids $\grd\rightrightarrows M$ carrying a multiplicative, closed 2-form $\omega$ satisfying a weaker version of nondegeneracy.

\end{itemize}

\begin{exer}
Verify that two Dirac structures on a manifold $M$ related by a gauge transformation $\tau_B$ (as in Exercise~\ref{exer:Bfield2}) have the same leaves, and the corresponding closed 2-forms on a given leaf $\mathcal{O}\hookrightarrow M$ differ by the pullback of $B$ to $\mathcal{O}$ (cf. Exercise~\ref{exer:Bfield}).
\end{exer}

The next exercise shows that any Dirac manifold gives rise to a Poisson algebra; building on this fact, the subsequent exercise indicates that a Dirac structure on a manifold $M$ can be also thought of as a Poisson structure on the leaf space of a foliation on $M$.

\begin{exer} \label{exer:adm}
For a Dirac manifold $(M,L)$, we say that a vector field $X$ is a {\em Hamiltonian vector field} for a function $f\in C^\infty(M)$ if $(X,df)\in \Gamma(L)$; functions that admit a Hamiltonian vector field are called {\em admissible} (note that it may happen that a function admits many or no Hamiltonian vector fields). Show that the set of admissible functions, denoted by $C^\infty_{adm}(M)$, is a Poisson algebra with bracket given by
$$
\{f,g\}=dg(X),
$$
where $X$ is any Hamiltonian vector field of $f$.
\end{exer}

\begin{exer}\label{exer:diracquot}
Let $(M,L)$ be a Dirac manifold, and let $\varphi: M \to B$ be a surjective submersion with connected fibers and such that $\ker(d\varphi)=\ker(L)$. Show that $\varphi^*: C^\infty(B)\to C^\infty(M)$ is an isomorphism onto $C^\infty_{adm}(M)$, and therefore (by Exercise \ref{exer:adm}) $B$ acquires a Poisson structure.
\end{exer}


Regarding ``morphisms'' of Dirac manifolds, there are two natural ways
in which Dirac manifolds can be related by a map, one extending the notion of Poisson map, and the other generalizing the pullback operation on 2-forms.

For Dirac manifolds $(M_1,L_1)$ and $(M_2,L_2)$, we say that a map $\varphi: M_1\to M_2$ is a {\em forward Dirac map} (or a {\em Dirac map}, for simplicity) if 
$$
L_2|_{\varphi(x)}=\{(d_x\varphi(X),\beta)  \;|\; (X,(d_x\varphi)^*\beta) \in L_1|_x\} \qquad \forall x \in M_1.
$$
We also write $L_2=\varphi_!L_1$ to denote that $L_1$ and $L_2$ are related as above.

\begin{exer} \label{exer:diracquot2}
In the context of Exercise~\ref{exer:diracquot}, check that the Poisson structure on $B$, denoted by $\pi$, satisfies $\varphi_! L=L_\pi$.
\end{exer}

\begin{exer}\label{exer:diracmap}
\begin{itemize}
\item[(a)] For Dirac structures $L_1$ and $L_2$ defined by Poisson structures, check that $\varphi: M_1\to M_2$ is a Dirac map if and only if it is a Poisson map.
\item[(b)] For Dirac structures $L_1$ and $L_2$ defined by closed 2-forms $\omega_1$ and $\omega_2$, show that $\varphi: M_1\to M_2$ is a Dirac map if and only if $\varphi$ is a submersion and $\varphi^*\omega_2=\omega_1$.
\end{itemize}
\end{exer}

As for pullbacks, consider a Dirac manifold $(M,L)$ and a smooth map $\varphi:N\to M$. The {\em backward image} of $L$ under $\varphi$ is defined by
$$
\varphi^!L:= \{(X,(d\varphi)^*\beta)\,|\, (d\varphi(X),\beta)\in L\}\subset TN\oplus T^*N.
$$
Over each point $y\in N$, the backward image of $L$ is a lagrangian subspace of $\mathbb{T}N|_y$, but it may happen that this family of lagrangian subspaces does not fit into a smooth vector subbundle of $\mathbb{T}N$ (see Exercise~\ref{exer:poissonpb} below); but if that is the case, then $\varphi^!L$ indeed defines a Dirac structure on $N$  \cite[Prop.~1.10]{bursztynDirac}, that we call the {\em pullback} of $L$.   

We collect some properties of pullbacks of Dirac structures and Dirac maps in the following three exercises.

\begin{exer}\label{exer:charLN}
Check that the characteristic distribution of the pullback $\varphi^!L$ is $(d\varphi)^{-1}(R)$, where $R$ is the characteristic distribution of $L$; moreover, the pointwise 2-forms on $(d\varphi)^{-1}(R)$ are  given by the pullbacks of those in $R$.
\end{exer}

\begin{exer} \label{exer:2formpb}
Let $\varphi:N\to M$ be a smooth map.
\begin{itemize}
\item[(a)] Check that $\varphi^!TM=TN$, and $\varphi^! T^*M= \ker(d\varphi)\oplus \mathrm{Ann}(\ker(d\varphi))$.
\item[(b)] Let $L$ be a Dirac structure on $M$, and let $\omega\in \Omega^2(M)$ be a closed 2-form. Check that 
$$
\varphi^!(\tau_\omega(L))=\tau_{\varphi^*\omega} (\varphi^!L).
$$
(In particular, $\varphi^!L_\omega= \varphi^!(\tau_\omega(TM))= \tau_{\varphi^*\omega}(TN)= L_{\varphi^*\omega}$, where $L_\omega$ is the Dirac structure defined by $\omega$.) 
\end{itemize}
\end{exer}

\begin{exer}\label{exer:funct}
Consider Dirac manifolds $(M_1,L_1)$ and $(M_2,L_2)$, and a smooth map $\varphi: M_1\to M_2$. 
\begin{itemize}
\item[(a)] Suppose that $\varphi_!L_1=L_2$. Check that $\varphi^!L_2=L_1$ if and only if $\ker(d\varphi)\subseteq \ker(L_1)$.
\item[(b)] Suppose that $\varphi^!L_2=L_1$. Check that $\varphi_!L_1=L_2$ if and only if  $R_2 |_{\varphi(x)}\subseteq d\varphi(T_xM_1)$ for all $x\in M_1$ (here $R_2$ is the projection of $L_2$ on $TM$).
\item[(c)] If $\varphi$ is a diffeomorphism, check that $\varphi_!L_1=L_2$ if and only if $L_1=\varphi^!L_2$.
\end{itemize} 
\end{exer}

The next exercise provides a simple example where the backward image of a Dirac structure is not smooth.

\begin{exer}\label{exer:poissonpb}
Let $L$ be the Dirac structure in $\mathbb{R}^2$ defined by the Poisson structure $x\partial_x\wedge\partial_y$.
Show that the backward image of $L$ under the inclusion $\varphi: \mathbb{R}\to \mathbb{R}^2$, $x\mapsto (x,0)$, is given by $T_x\mathbb{R}\oplus 0 \subseteq \mathbb{T}_x\mathbb{R}$ for $x\neq 0$, and for $x=0$ it is $0\oplus T^*_0\mathbb{R}$.
\end{exer}

A sufficient (but not necessary) condition to ensure the smoothness of the family $\varphi^!L$ is that \cite[Prop.~1.10]{bursztynDirac} 
\begin{itemize}
\item the  subspaces
$\ker((d_x\varphi)^*)\cap L|_{\varphi(x)}\subseteq T^*M|_{\varphi(x)}$, for $x\in N$, have the same dimension, or equivalently (by considering their annihilators) that the subspaces $\mathrm{Im}(d_x\varphi)+R|_{\varphi(x)}\subseteq TM|_{\varphi(x)}$ have the same dimension. 
\end{itemize}
We will refer to this property as the {\em co-regularity condition},
following the terminology in \cite{Geudens}.
In particular,  Dirac structures can be always pulled back by submersions (or, more generally, by maps that are transverse to their characteristic distributions).

Whenever the co-regularity condition holds, {\em the presymplectic leaves of the Dirac structure $\varphi^!L$ coincide with the connected components of 
$\varphi^{-1}(\mathcal{O})$, equipped with the 2-form $\varphi^*\omega_{\mathcal{O}}$, for each presymplectic leaf $(\mathcal{O},\omega_\mathcal{O})$ of $L$} \cite[Thm.~7.33]{CFMbook}. (In spite of what Exercise~\ref{exer:charLN} may suggest, this description of presymplectic leaves is not valid in general, without additional ``cleanness'' conditions, see \cite{CorFrejMart} for examples.)

\begin{exer} 
\begin{itemize}
\item[(a)] Let $\varphi: M\to B$ be a surjective submersion with connected fibers. Verify that the assignment $\pi\mapsto \varphi^!L_\pi$ establishes a bijection between Poisson
structures on $B$ and Dirac structures $L$ on $M$ satisfying $\ker(d\varphi)=\ker(L)$ (cf. Exercises~\ref{exer:diracquot} and \ref{exer:diracquot2}).
\item[(b)] Suppose that $M$ carries a $G$-action that is free and proper, and $\varphi: M\to B=M/G$ is the quotient map. Show that $\pi\mapsto \varphi^!L_\pi$ yields a one-to-one correspondence between Poisson structures on $B$ and $G$-invariant Dirac structures $L$ on $M$ for which $\ker(L)$ coincides with the distribution tangent to $G$-orbits (here $G$-invariance means that $\sigma^!L=L$, or equivalently $\sigma_!L=L$, for all $\sigma\in G$).
\end{itemize}

\end{exer}

The next exercise uses pullbacks and gauge transformations to give a Dirac-geometric characterization of dual pairs \cite{FrejlichMarcut}.

\begin{exer}\label{exer:dualdirac}
Consider a symplectic manifold $(S,\omega)$, Poisson manifolds $(M_1,\pi_1)$ and $(M_2,\pi_2)$, and let $\varphi_i: S\to M_i$, $i=1,2$, be surjective submersions.
Show that $M_1\stackrel{\varphi_1}{\leftarrow} S \stackrel{\varphi_2}{\to} \overline{M}_2$ is a dual pair (see \eqref{eq:dual}) if and only if $\dim(S)=\dim(M_1)+\dim(M_2)$ and
$$
\varphi_1^!(L_{\pi_1})=\tau_\omega(\varphi_2^!L_{\pi_2}).
$$
\end{exer}

We now suppose that  $N\stackrel{\iota}{\hookrightarrow} M$ is a submanifold. We know from the discussion so far that, for a Dirac structure $L$ on $M$,
the backward image $\iota^!L$ defines a Dirac structure on $N$ whenever it is smooth.
The co-regularity condition assuring  the smoothness of $\iota^!L$ in this case is that $\mathrm{Ann}(TN)\cap L|_N$ has constant rank \cite[$\S$ 3.1]{courant}, or equivalently that
\begin{equation*}\label{eq:condpb}
  TN+R|_N\subseteq TM|_N \; \mbox{ has constant rank.}
\end{equation*}
We refer to a submanifold satisfying these conditions as {\em co-regular}. 
If $TN+ R|_N=TM|_N$, we call the submanifold $N$ a {\em transversal}.

\begin{exer}\label{exer:kerpb} 
Check that the kernel of the pullback Dirac structure $\iota^!L$ on $N$ is the distribution given by the projection on $TN$ of $(TN\oplus \mathrm{Ann}(TN))\cap L$.
\end{exer}

\begin{exer}\label{exer:diractransv} Let $(M,L)$ be a Dirac manifold, and suppose that $\iota: N\hookrightarrow M$ is a submanifold containing $x\in M$ and such that $T_xM = R|_x \oplus T_xN = T_x\mathcal{O}\oplus T_xN$, where $\mathcal{O}$ is the presymplectic leaf through $x$. Show 
that there exists a neighborhood of $x$ in $N$ where $\iota^!L$ is given by a Poisson structure. 
\end{exer}

In the special case where $M$ is a Poisson manifold, the Poisson structure of the previous exercise is just the transverse Poisson structure of $\S$ \ref{subsec:transverse}, as given in Exercise~\ref{exer:slice} (compare Exercises~\ref{exer:charpiN} and \ref{exer:charLN}).
A uniqueness theorem for these more general ``transverse Poisson structures'' is proven in \cite[Thm.~4.5]{DufourWade}, see also \cite[$\S$ 5.6]{BLM}.

The next exercise shows that transverse Poisson structures are preserved under pullbacks by submersions. This fact will be used in Exercise~\ref{exer:dualprestransv} to prove the invariance of transverse Poisson structures in dual pairs stated in $\S$ \ref{subsec:sympreal}.

\begin{exer}\label{exer:transvpres}  Let $(M,L)$ be a Dirac manifold, $\varphi: S \to  M$ a submersion, and consider the pullback Dirac structure $\varphi^!L$ on $S$. Pick $m\in S$, let $\mathcal{O}$ be the presymplectic leaf through $x=\varphi(m)$ and $\mathcal{O}':=\varphi^{-1}(\mathcal{O})$. Let $ N'\hookrightarrow S$ be a submanifold containing $m$ such that
$T_mS= T_m(\mathcal{O}')\oplus T_mN'$. Check that, by shrinking $N'$ and using Exercise~\ref{exer:diractransv}, one can assume that
\begin{itemize}
\item $N'$ inherits a Poisson structure by the pullback of $\varphi^!L$;
\item $\varphi$ restricts to an embedding $N'\to M$, with image denoted by $N\subseteq M$;
\item $N$ inherits a Poisson structure by the pullback of $L$.
\end{itemize}
Verify that $\varphi$ restricts to a Poisson isomorphism $N'\to N$.
\end{exer}

\begin{exer}\label{exer:dualprestransv}
Let  $(M_1,\pi_1) \stackrel{\varphi_1}{\leftarrow} (S,\omega) \stackrel{\varphi_2}{\to} ({M}_2,-\pi_2)$ be a  dual pair (see \eqref{eq:dual}). For any $m\in S$, show that the transverse Poisson structures at $\varphi_1(m)$ and $\varphi_2(m)$ are isomorphic. The following are the key steps:
\begin{itemize}
\item Use the Darboux theorem to find a submanifold $N'\hookrightarrow S$ such that $T_mS=T_mN' \oplus R'|_m$, where $R'= \ker(d\varphi_1) + \ker(d\varphi_2)$, that is isotropic.
\item By Exercises \ref{exer:2formpb} (b) and \ref{exer:dualdirac}, $\varphi_1^! L_{\pi_1}$ and $\varphi_2^! L_{\pi_2}$ have the same pullback on $N'$.
\item By shrinking $N'$, this pullback is a Poisson structure that is isomorphic to the transverse Poisson structures at $\varphi_1(m)$ and $\varphi_2(m)$, by Exercise~\ref{exer:transvpres}.
\end{itemize}
\end{exer}

Let us now specialize to  the case where $(M,\pi)$ is a Poisson manifold and  $\iota: N \hookrightarrow M$ is a submanifold. As we have seen, modulo smoothness issues (cf. Exercise \ref{exer:poissonpb}), $N$ generally carries a Dirac structure $\iota^!L_\pi$ via pullback of the ambient Poisson structure. 

Recalling that $TN^\pi=\pi^\sharp(\mathrm{Ann}(TN))$, note that
\begin{itemize}
\item $\ker(\iota^!L_\pi)=TN\cap TN^\pi$ (see Exercise \ref{exer:kerpb});
\item $N$ is co-regular if and only if $\mathrm{Ann}(TN)\cap \ker(\pi^\sharp)$ has constant rank, which is equivalent to the condition that
$$
TN^\pi  \mbox{ has constant rank},
$$
since $\mathrm{Ann}(TN)\cap \ker(\pi^\sharp)\hookrightarrow \mathrm{Ann}(TN)\twoheadrightarrow TN^\pi$ is exact.
\end{itemize}

A submanifold $\iota: N \hookrightarrow M$ for which the backward image $\iota^!L_\pi$ is smooth (hence Dirac) and has trivial kernel ($TN\cap TN^\pi=0$) is called a {\em Poisson-Dirac submanifold} of $M$. Equivalently, a Poisson-Dirac submanifold $N$ is characterized by the fact that it carries a Poisson structure $\pi_N$ uniquely determined by the condition
\begin{equation*}\label{eq:PDcond}
L_{\pi_N}=\iota^!L_\pi.
\end{equation*}

\begin{example}\ \label{example:PDsubm}
\begin{itemize}
\item[(a)] A submanifold $N\stackrel{\iota}{\hookrightarrow} M$ is Poisson if and only if $TN^\pi=0$ (see Exercise~\ref{exer:tangent} (a)),
so Poisson submanifolds are co-regular Poisson-Dirac submanifolds.
Recall that, being a Poisson submanifold, $N$ carries a Poisson structure $\pi_N$ for which the inclusion is a Poisson map, which means that $L_\pi = \iota_! L_{\pi_N}$ (Exercise~\ref{exer:diracmap} (a)). But then $\pi_N$ satisfies $L_{\pi_N}=\iota^!L_\pi$ by Exercise~\ref{exer:funct} (a).

\item[(b)] Since cosymplectic submanifolds satisfy $TM|_N=TN\oplus TN^\pi$, they are co-regular Poisson-Dirac submanifolds, and the Poisson structure $\pi_N$ that $N$
inherits as a cosymplectic submanifold indeed satisfies $L_{\pi_N}=\iota^!L_\pi$ (compare its characterization in Exercise~\ref{exer:charpiN} with Exercise~\ref{exer:charLN}). 
Note that, by Exercise~\ref{exer:cosympR}, cosymplectic submanifolds are the same as Poisson-Dirac submanifolds that are also transversals\footnote{This fact suggests ``Poisson-Dirac transversals'', rather than ``Poisson transversals'', as a more accurate alternative terminology for ``cosymplectic submanifolds''.}.

\item[(c)] Poisson-Dirac submanifolds of a symplectic manifold are the same as symplectic submanifolds (see Exercise~\ref{exer:2formpb}). Note that in a symplectic manifold any submanifold is co-regular.
\end{itemize}
\hfill $\diamond$
\end{example}


It is shown in \cite{Geudens} that any co-regular Poisson-Dirac submanifold of $M$ is a Poisson submanifold in a cosymplectic submanifold of $M$ (see also \cite[Thm.~8.44]{CFMbook}).

On a co-regular Poisson-Dirac submanifold $N\hookrightarrow M$, symplectic leaves are given by the connected components of  the intersections of $N$ with the symplectic leaves in $M$, equipped with the pullback of the symplectic forms. The description of the induced Poisson bracket on $N$ in the next exercise is entirely analogous to the one for cosymplectic submanifolds in Exercise~\ref{exer:cosymPoiss}.

\begin{exer} For a co-regular Poisson-Dirac submanifold $N\hookrightarrow M$ (embedded, for simplicity), verify that the induced Poisson bracket $\{\cdot,\cdot\}_N$ on $N$ is given by
$$
\{f,g\}_N= \{\widehat{f},\widehat{g}\}|_N,
$$
where $\widehat{f}, \widehat{g} \in C^\infty(M)$ are extensions of $f, g \in C^\infty(N)$ satisfying $d\widehat{f}|_{TN^\pi} = d\widehat{g}|_{TN^\pi}=0$.
\end{exer}

For a detailed account of Poisson-Dirac submanifolds with different types of regularity conditions, see \cite{CorFrejMart}. For a broader discussion including other types of submanifolds, see e.g. \cite[Chp.~8]{CFMbook} and \cite{Zambon}, as well as \cite{Geudens}.

Besides offering an effective framework for the study of submanifolds of Poisson manifolds, Dirac structures
are essential in several other aspects of Poisson geometry, especially due to their flexible pullbacks and enlarged group of symmetries (see Exercise~\ref{exer:Bfield2}). They can be used to prove, clarify, and often extend, several fundamental results in Poisson geometry, see \cite{MeinDirac}. Those include the classification of Poisson homogeneous spaces \cite{LWX,MeinLu} (see Remark~\ref{rem:double}), the integration of Poisson homogeneous spaces by symplectic groupoids \cite{BIL}, the construction of symplectic realizations \cite{FrejlichMarcut},  linearization of Poisson-Lie groups \cite{AMLin}, and normal forms of Poisson structures around transversals (extending Weinstein's splitting theorem) \cite{BLM,FrejMarcNormal}.

We have only discussed Dirac structures in the ``standard'' Courant algebroid $\TM$.
Many situations of interest involve other types of Courant algebroids, such as ``exact'' ones, where the Courant-Dorfman bracket \eqref{eq:CDbrk} is twisted by a closed 3-form \cite{SW}. In this extended setting Dirac structures and Dirac maps provide a geometric framework for quasi-Hamiltonian geometry \cite{AMM} (see \cite{ABM,BCWZ}) and become natural ingredients in the construction of Poisson structures on moduli spaces of flat bundles over surfaces and their variants \cite{LiBland}. 
Dirac structures are also studied in complexified Courant algebroids, especially
in connection with generalized complex structures \cite{Gua}, see \cite{AgRu}.


\end{document}